\documentclass[11pt,a4paper,reqno]{amsart}
\usepackage{amssymb,cite,mathtools}
\usepackage{array,placeins}
\usepackage[T1]{fontenc}
\usepackage{lmodern,microtype}
\usepackage[hidelinks]{hyperref}
\ifx\pdfpagewidth\undefined
  \AtBeginDocument{}
\else
  \pdfpagewidth=\paperwidth
  \pdfpageheight=\paperheight
\fi
\makeatletter
\renewcommand*\subjclass[2][2020]{%
  \def\@subjclass{#2}%
  \def\subjclassname{\textup{#1} Mathematics Subject Classification}}
\newdimen\paper@needed
\newdimen\paper@remaining
\newcommand{\Needspace}[1]{%
  \par\penalty-100\begingroup
  \paper@needed=#1\relax
  \paper@remaining=\pagegoal
  \advance\paper@remaining by -\pagetotal
  \ifdim\paper@needed>\paper@remaining
    \ifdim\paper@remaining>0pt \vfil\fi
    \break
  \fi
  \endgroup}
\renewenvironment{enumerate}[1][]{%
  \ifnum\@enumdepth>3 \@toodeep\else
    \advance\@enumdepth\@ne
    \edef\@enumctr{enum\romannumeral\the\@enumdepth}%
    \expandafter\edef\csname the\@enumctr\endcsname{%
      \noexpand\roman{\@enumctr}}%
    \list{\textup{(\csname the\@enumctr\endcsname)}}{%
      \usecounter{\@enumctr}%
      \setlength{\leftmargin}{2em}%
      \setlength{\labelwidth}{\leftmargin}%
      \addtolength{\labelwidth}{-\labelsep}%
      \setlength{\itemindent}{0pt}%
      \setlength{\topsep}{4pt}%
      \setlength{\itemsep}{3pt}%
      \setlength{\parsep}{0pt}%
      \def\makelabel##1{\hss\llap{##1}}}%
  \fi
}{\endlist}
\renewcommand\subsection{\@startsection{subsection}{2}\z@
  {1.1\baselineskip plus .3\baselineskip minus .2\baselineskip}
  {.35\baselineskip}{\normalfont\bfseries}}
\makeatother
\numberwithin{equation}{section}
\newtheorem{theorem}{Theorem}[section]
\newtheorem{proposition}[theorem]{Proposition}
\newtheorem{lemma}[theorem]{Lemma}
\newtheorem{corollary}[theorem]{Corollary}
\theoremstyle{definition}
\newtheorem{definition}[theorem]{Definition}
\theoremstyle{remark}
\newtheorem{remark}[theorem]{Remark}
\newcommand{\R}{\mathbb R}
\newcommand{\dd}{\,\mathrm d}
\newcommand{\im}{\operatorname{Im}}

\newcommand{\Ec}{E_{\mathrm c}}
\newcommand{\norm}[1]{\left\|#1\right\|}
\newcommand{\avgT}[1]{\frac1T\int_0^T #1\dd t}
\newcommand{\proofstep}[2]{%
  \par\addvspace{.5\baselineskip}\Needspace{3\baselineskip}%
  \noindent\textbf{Step #1. #2.}\par\nobreak\smallskip
  \noindent\ignorespaces}
\allowdisplaybreaks[2]
\title[Threshold scattering for the planar cubic--quintic NLS]{Nonradial threshold scattering for the\newline planar cubic--quintic NLS}
\author[Bo Yang]{Bo Yang}
\address{School of Mathematics and Statistics, Huazhong University of Science and Technology, Wuhan 430074, Hubei, P. R. China; Hubei Key Laboratory of Engineering Modeling and Scientific Computing, Huazhong University of Science and Technology, Wuhan 430074, Hubei, P. R. China}
\email{boyang\_0717@hust.edu.cn}
\author[Lei Zhang]{Lei Zhang}
\address{School of Mathematics and Statistics, Huazhong University of Science and Technology, Wuhan 430074, Hubei, P. R. China; Hubei Key Laboratory of Engineering Modeling and Scientific Computing, Huazhong University of Science and Technology, Wuhan 430074, Hubei, P. R. China}
\email{lei\_zhang@hust.edu.cn}
\author[Bin Liu]{Bin Liu$^*$}
\address{School of Mathematics and Statistics, Huazhong University of Science and Technology, Wuhan 430074, Hubei, P. R. China; Hubei Key Laboratory of Engineering Modeling and Scientific Computing, Huazhong University of Science and Technology, Wuhan 430074, Hubei, P. R. China}
\email{binliu@mail.hust.edu.cn}
\thanks{$^*$Corresponding author.}
\thanks{This work was partially supported by the National Key Research and
Development Program of China (Grant No.~2023YFC2206100),
the Fundamental Research Funds for the Central Universities
(Grant No.~2026BRSXB003), and the National Natural Science Foundation
of China (Grant No.~12231008).}
\date{}
\hypersetup{
  pdftitle={Nonradial threshold scattering for the planar cubic-quintic NLS},
  pdfauthor={Bo Yang, Lei Zhang, and Bin Liu},
  pdfsubject={Threshold scattering for the planar cubic-quintic nonlinear Schrodinger equation},
  pdfkeywords={cubic-quintic nonlinear Schrodinger equation, threshold scattering, nonradial solutions, concentration compactness, interaction Morawetz estimates, propagation estimates}
}
\begin{document}
\begin{abstract}
We study the planar cubic--quintic nonlinear Schr\"odinger equation
\[
 i\partial_tu+\Delta u=-|u|^2u+|u|^4u,
 \qquad (t,x)\in\mathbb R\times\mathbb R^2,
\]
with initial data $u_0\in H^1_x(\mathbb R^2)$. We prove scattering in both time directions
when the mass of \(u_0\) equals the mass of the ground state for the focusing cubic equation.
This resolves the nonradial threshold scattering problem posed by Murphy
[\emph{Comm.\ Partial Differential Equations} \textbf{46} (2021), 2213--2234].
Together with the known subthreshold scattering theory and the existence
of standing waves at every larger mass, this identifies the sharp mass
range in which every $H^1$ solution scatters. The proof combines localized
interaction Morawetz estimates with propagation estimates on moving
halfspaces to control exterior kinetic energy.
\end{abstract}
\maketitle

\section{Introduction}
\label{sec:introduction}
We consider the Cauchy problem for the planar cubic--quintic nonlinear
Schr\"odinger equation (NLS)
\begin{equation}\label{eq:cq}
 i\partial_tu+\Delta u=F(u),\qquad
 F(z):=-|z|^2z+|z|^4z,\qquad (t,x)\in\R\times\R^2,
\end{equation}
with $u(0)=u_0\in H^1_x(\R^2)$. The cubic term is focusing, and the
quintic term is defocusing. The competing nonlinearities are motivated
by mean-field models of Bose gases with attractive two-body and repulsive three-body
interactions \cite{NguyenCQ}.

Unlike the classical NLS, equation~\eqref{eq:cq} has no nontrivial scaling symmetry.
This prevents us from normalizing the spatial
scale while keeping the equation fixed and complicates the compactness
argument. However, \eqref{eq:cq} is invariant
under phase rotations, translations in time and space, and orthogonal
transformations of the spatial variables, as well as under the
Galilean transformation
\[
 u(t,x)\longmapsto u_a(t,x)
 :=e^{i(a\cdot x-|a|^2t)}u(t,x-2at),\qquad a\in\R^2,
\]
and under time reversal combined with complex conjugation,
$u(t,x)\mapsto\overline{u(-t,x)}$.
For $f\in H^1_x(\R^2)$, define the mass, energy, and momentum by
\begin{equation}\label{eq:intro_conserved_quantities}
 \begin{aligned}
 M(f)&:=\int_{\R^2}|f(x)|^2\,dx,\\
 E(f)&:=\int_{\R^2}\left(\frac12|\nabla f(x)|^2
                 -\frac14|f(x)|^4+\frac16|f(x)|^6\right)\,dx,\\
 P(f)&:=\operatorname{Im}\int_{\R^2}\overline{f(x)}\nabla f(x)\,dx.
 \end{aligned}
\end{equation}
These quantities are conserved along $H^1$ solutions.

Consider the elliptic equation
\begin{equation}\label{eq:threshold}
 -\Delta Q+Q-Q^3=0,\qquad x\in\R^2.
\end{equation}
This equation admits a unique positive radial solution
$Q\in H^1_x(\R^2)$ \cite{Weinstein83,McLeodSerrin}.
Set $m_Q:=M(Q)$. The function $Q$ attains equality in the sharp
Gagliardo--Nirenberg inequality \cite{Weinstein83},
\begin{equation}\label{eq:GN}
 \norm f_{L^4_x(\R^2)}^4
 \leq 2\frac{M(f)}{m_Q}\norm{\nabla f}_{L^2_x(\R^2)}^2
 \qquad(f\in H^1_x(\R^2)).
\end{equation}
It also generates the standing wave $v(t,x)=e^{it}Q(x)$ for the
focusing cubic equation
\begin{equation}\label{eq:cubic_nls}
 i\partial_tv+\Delta v=-|v|^2v,\qquad (t,x)\in\R\times\R^2.
\end{equation}
We call Q the ground state of \eqref{eq:cubic_nls}.
For equation \eqref{eq:cubic_nls}, every $L^2$ solution of mass less than $m_Q$
scatters in both time directions \cite{DodsonSub}, while $e^{it}Q$
is a nonscattering solution at mass $m_Q$.
Hence $m_Q$ is the sharp mass threshold for scattering of all
solutions of the focusing cubic equation.

The same mass is relevant to \eqref{eq:cq}. Under the mass-preserving
rescaling $u(t,x)=h^{-1}v(t/h^2,x/h)$, the cubic coefficient remains
unchanged and the quintic coefficient becomes $h^{-2}$.
Thus the focusing cubic equation describes the large-scale limit,
where its ground-state standing wave arises as a limiting profile.
For \eqref{eq:cq}, Cheng \cite{Cheng} proved scattering for any
$H^1$ data with $M(u_0)<m_Q$. At each prescribed mass $m>m_Q$,
Carles and Sparber \cite{CS} constructed a standing wave, so a scattering
statement for all data of that mass is impossible.
This leaves the following endpoint question:
\emph{Does every $H^1$ solution of \eqref{eq:cq} with $M(u_0)=m_Q$
scatter in both time directions without additional symmetry or moment
assumptions?}

\subsection{Related work}
\noindent\textit{The classical focusing NLS.}
The classical threshold theory concerns the single-power equation
\begin{equation}\label{eq:classical_focusing_nls}
 i\partial_tv+\Delta v=-|v|^\alpha v,
 \qquad (t,x)\in\R\times\R^d,
 \qquad d\geq1,\quad \alpha>0.
\end{equation}
In the mass-critical case $\alpha=4/d$, scattering below the
ground-state mass is known for arbitrary $L^2(\R^d)$ data in all
dimensions \cite{KTV,DodsonSub}. At the threshold, Merle \cite{Merle}
classified finite-time blowup solutions in $H^1$, and Dodson
\cite{Dodson} classified nonscattering $L^2$ solutions for
$2\leq d\leq15$.
For the three-dimensional cubic equation, scattering below both
the ground-state mass--energy and mass--gradient thresholds was
proved for radial $H^1$ data by Holmer and Roudenko
\cite{HolmerRoudenko}, and for general $H^1$ data by Duyckaerts,
Holmer, and Roudenko \cite{DHR}.
In the energy-critical case $\alpha=4/(d-2)$, scattering below
both the ground-state energy and gradient norm was established
by Kenig and Merle \cite{KM} for radial $\dot H^1$ data in
$d=3,4,5$, and subsequently for general data in $d\geq5$
\cite{KillipVisan} and $d=4$ \cite{Dodson4D}.
The works \cite{DuyckaertsRoudenko,DuyckaertsMerle} exhibit special
threshold solutions that scatter in one time direction and approach
a solitary wave in the other, so the endpoint requires a separate
analysis.

\medskip\noindent\textit{The cubic--quintic NLS.}
We next consider the competing-power equation
\begin{equation}\label{eq:cq_dimensions}
 i\partial_tu+\Delta u=-|u|^2u+|u|^4u,
 \qquad (t,x)\in\R\times\R^d,\qquad d\geq2.
\end{equation}
In three dimensions, Zhang \cite{ZhangCQ} proved global
well-posedness for arbitrary $H^1$ data.
Killip, Oh, Pocovnicu, and Visan \cite{KOPV} established scattering
in a variationally defined mass--energy region, without radial or
finite-variance assumptions. Killip, Murphy, and Visan \cite{KMV}
extended this region across the part of its boundary not attained
by solitons. Ardila and Murphy \cite{ArdilaMurphy} classified
solutions with nonnegative virial on a specified positive-energy
soliton portion of the original boundary.

In two dimensions, the subthreshold theorem of Cheng \cite[Theorem~1.2]{Cheng}
requires only $u_0\in H^1_x(\R^2)$ and $M(u_0)<m_Q$.
At $M(u_0)=m_Q$, Carles and Sparber \cite[Theorem~1.8]{CS}
proved scattering for data in the conformal space
\[
 \Sigma:=\{f\in H^1_x(\R^2):xf\in L^2_x(\R^2)\},
\]
using a pseudoconformal estimate. Their result holds for
$M(u_0)\leq m_Q$ and gives convergence of $e^{-it\Delta}u(t)$
in $\Sigma$ as $t\to\pm\infty$.
Murphy \cite[Theorem~1.1]{Murphy} proved threshold scattering for
radial $H^1$ data without a finite-variance assumption.

The superthreshold standing waves in \cite[Theorem~1.9]{CS}
have also been studied through constrained energy minimization.
Li and Wang \cite{LiWangCQ} obtained the optimal rate at which
the gradient energy of normalized ground states vanishes as their
mass decreases to $m_Q$, together with the limiting profile.
Jiang, Wang, Xiao, Zhang, and Zhu \cite{JiangEtAlCQ} proved uniqueness, up to translation, of the positive minimizer of the energy at each prescribed mass $m>m_Q$.

These two-dimensional results cover scattering below $m_Q$, standing
waves above $m_Q$, and threshold scattering under radiality or
finite-variance assumptions. They leave open the threshold problem
for general $H^1$ data, explicitly posed by Murphy \cite{Murphy}.
We therefore consider \eqref{eq:cq} with $M(u_0)=m_Q$ without
either assumption and prove scattering in both time directions.
Combined with the subthreshold scattering theory and the existence
of standing waves at every larger mass, our result identifies
$M(u_0)\leq m_Q$ as the sharp mass range in which every $H^1$
solution scatters.

\begin{table}[!htbp]
\caption{Scattering and standing waves for the planar cubic--quintic NLS}
\label{tab:planar_cq_results}
\centering\footnotesize
\renewcommand{\arraystretch}{1.2}
\setlength{\tabcolsep}{4pt}
\begin{tabular}{@{}>{\raggedright\arraybackslash}p{.41\textwidth}>{\raggedright\arraybackslash}p{.28\textwidth}>{\raggedright\arraybackslash}p{.27\textwidth}@{}}
\hline
Assumptions & Conclusion & Reference\\
\hline
$M(u_0)<m_Q$
 & $H^1$ scattering & Cheng \cite{Cheng}\\
$u_0\in\Sigma$, $M(u_0)\leq m_Q$
 & Scattering in $\Sigma$ & Carles--Sparber \cite{CS}\\
$u_0\in H^1_{\rm rad}$, $M(u_0)=m_Q$
 & $H^1$ scattering & Murphy \cite{Murphy}\\
Each prescribed mass $m>m_Q$
 & Existence of a standing wave & Carles--Sparber \cite{CS}\\
\hline
$M(u_0)=m_Q$
 & $H^1$ scattering & Theorem~\ref{thm:threshold_scattering}\\
\hline
\end{tabular}
\end{table}
\FloatBarrier

\subsection{The main result}
We first define strong $H^1$ solutions of \eqref{eq:cq}.
\begin{definition}[Strong $H^1$ solution]
\label{def:strong_H1_solution}
Let $I\subset\R$ be an interval containing $t_0$, and let
$f\in H^1_x(\R^2)$. A strong $H^1$ solution of \eqref{eq:cq} on $I$
with initial data $u(t_0)=f$ is a function
$u\in C(I;H^1_x(\R^2))$ such that, for every compact interval $J\subset I$,
\[
 u\in L_t^4(J;W_x^{1,4}(\R^2))
       \cap L_t^8(J;W_x^{1,8/3}(\R^2)),
\]
and
\begin{equation}\label{eq:appendix_duhamel}
 u(t)=e^{i(t-t_0)\Delta}f
      -i\int_{t_0}^t e^{i(t-\tau)\Delta}F(u(\tau))\,d\tau,
 \qquad t\in I.
\end{equation}
The identity holds in $H^1_x(\R^2)$.
\end{definition}
The estimate $E(f)\geq\frac12\|\nabla f\|_2^2-\frac{3}{32}M(f)$,
together with the local theory, gives a unique global strong
solution of \eqref{eq:cq} for any initial data in $H^1_x(\R^2)$,
with a uniform $H^1$ bound \cite{CS}.
Our main result establishes scattering for these solutions at the
ground-state mass, without radial or finite-variance assumptions.

\Needspace{9\baselineskip}
\begin{theorem}[Threshold scattering]
\label{thm:threshold_scattering}
Let $u_0\in H^1_x(\R^2)$ satisfy $M(u_0)=m_Q$, and let $u$ be the
solution of~\eqref{eq:cq} with $u(0)=u_0$. Then
$u\in C(\R;H^1_x(\R^2))$, and there exist $u_\pm\in H^1_x(\R^2)$ such that
\begin{equation}\label{eq:scattering_conclusion}
 \|u\|_{L^4_{t,x}(\R\times\R^2)}<\infty,\qquad
 \|u(t)-e^{it\Delta}u_\pm\|_{H^1}\longrightarrow0
 \quad(t\to\pm\infty).
\end{equation}
\end{theorem}
The theorem requires neither radial symmetry nor finite variance.
Together with Cheng's result \cite{Cheng}, it proves scattering for any
$H^1$ data with $M(u_0)\leq m_Q$. Since a standing wave exists at
each mass $m>m_Q$, the mass condition is optimal for a statement
that applies to all initial data of the prescribed mass.

\subsection{Difficulties and strategy of the proof}
The first difficulty is that, at $M(f)=m_Q$, the right-hand side
of the coercive estimate
\[
 \frac12\|\nabla f\|_{L_x^2(\R^2)}^2
 -\frac14\|f\|_{L_x^4(\R^2)}^4
 \geq\frac12\left(1-\frac{M(f)}{m_Q}\right)
             \|\nabla f\|_{L_x^2(\R^2)}^2
\]
vanishes. Thus this estimate no longer controls the kinetic energy
with a strictly positive coercivity constant.
A second difficulty appears in the large-scale limit:
under $u(t,x)=h^{-1}v(t/h^2,x/h)$, equation~\eqref{eq:cq} becomes
\[
 i\partial_sv+\Delta_yv=-|v|^2v+h^{-2}|v|^4v.
\]
As $h\to\infty$, the limiting equation is \eqref{eq:cubic_nls}.
Its soliton can carry all the mass while an $L^2$-small remainder
retains kinetic energy, since
\[
 \|r_n\|_{L_x^2(\R^2)}\to0
 \quad\not\Rightarrow\quad
 \|\nabla r_n\|_{L_x^2(\R^2)}\to0.
\]
Here $r_n=u(t_n)-q_n$ is the remainder after subtracting the
cubic soliton profile $q_n$ along a sequence with
$\lambda(t_n)\to\infty$; see \eqref{eq:Qdefectform} in
Corollary~\ref{cor:Qdensity}.
Thus, We must control the exterior kinetic energy directly.

In the radial argument of Murphy \cite{Murphy}, radial Sobolev
estimates control the nonlinear tails, while the incoming/outgoing
decomposition and weighted radial Strichartz estimates yield
kinetic-energy localization. These tools require radial data.
For general data, mass compactness and the Gagliardo--Nirenberg
inequality control the nonlinear spatial tails, but do not control
the exterior gradient. However, localized interaction Morawetz estimates
give exterior spacetime bounds. We combine these bounds with
propagation between moving halfspaces to control the exterior
gradient through its Duhamel representation.

\medskip\noindent\textit{Compactness and the virial obstruction.}
Following the concentration--compactness approach of Kenig and
Merle \cite{KM} and the reduction in \cite{Murphy}, failure of the
theorem yields a global solution with
\begin{equation}\label{eq:minimal_mass_nonscattering}
 M(u)=m_Q,\qquad P(u)=0,\qquad
 \|u\|_{L^4_{t,x}([0,\infty)\times\R^2)}=\infty.
\end{equation}
Here time reversal with complex conjugation fixes the nonscattering direction, and a
Galilean transform gives zero momentum. We call $u$ a
\emph{minimal-mass nonscattering solution}.

Subthreshold scattering excludes splitting into several nonzero profiles.
A large-scale cubic limit cannot scatter in either time direction:
forward scattering would control a forward tail of $u$, and backward
scattering would bound its norm on $[0,T_n]$ for $T_n\to\infty$.
Dodson's classification \cite{Dodson} identifies this limit as a soliton.
Comparison in rescaled time and the center argument in \cite{DHR}
give a spatial scale and center satisfying
\[
 \lambda(t)=o(\sqrt t),\qquad c(t)=o(t).
\]
Together with mass compactness, these bounds imply
\[
 \int_{|x|>ct}|u(t,x)|^2\,dx\longrightarrow0
 \qquad(t\to\infty,\ c>0).
\]

Let $E>0$ be the conserved energy in this zero-momentum frame.
The localized virial argument and
Corollary~\ref{cor:ballistic_evacuation} give
\[
 \begin{aligned}
 \frac1T\int_0^T\|\nabla u(t)\|_{L_x^2(\R^2)}^2\,dt&\longrightarrow2E,\\
 \frac1T\int_0^T\int_{|x|\leq cT}|\nabla u(t,x)|^2\,dx\,dt
 &\longrightarrow0\qquad(c>0).
 \end{aligned}
\]
Consequently,
\begin{equation}\label{eq:intro_exterior_energy_average}
 \lim_{T\to\infty}\frac1T\int_0^T\int_{|x|>ct}
                  |\nabla u(t,x)|^2\,dx\,dt=2E
 \qquad(c>0).
\end{equation}
The remaining task is to prove decay of the kinetic energy in
these same exterior regions.

\medskip\noindent\textit{Exterior estimates on dyadic intervals.}
We localize the interaction estimates of
\cite{CGT,PlanchonVega} to halfspaces. Small exterior mass allows
us to absorb the focusing cubic pressure through nested cutoffs;
the quintic pressure has a favorable sign. For each spatial cutoff,
we partition $[T,2T]$ into suitable subintervals and apply Strichartz
estimates to the localized equation, using the truncated solution
at each left endpoint as initial data. This gives
\[
 \|u\|_{L_{t,x}^4(\{(t,x):T<t<2T,\ |x|>ct\})}
 +\|u\|_{L_{t,x}^8(\{(t,x):T<t<2T,\ |x|>ct\})}
 \longrightarrow0
 \quad(T\to\infty)
\]
for every $c>0$; see
Corollary~\ref{cor:cubic_window_unconditional} and
\eqref{eq:exterior_im_doubling}.
The corresponding derivative estimate is
\[
 \|\chi(x/T)F(u)\|_{L_t^{4/3}([T,2T];W_x^{1,4/3}(\R^2))}
 \longrightarrow0\qquad(T\to\infty),
\]
for each fixed exterior cutoff $\chi$ from
\eqref{eq:spatial_cutoff}; see
Corollary~\ref{cor:actual_exterior_cubic_N1}.
Both estimates are needed in the global gradient argument.

\medskip\noindent\textit{Directional propagation and summation in time.}
Dyadic decay alone does not give a finite global spacetime norm.
Inspired by the incoming/outgoing decomposition in \cite{Murphy},
we separate velocities relative to a moving line:
\[
 x\cdot e=at,\qquad |e|=1,\qquad a>0.
\]
The normal group velocity is $2\xi\cdot e$.
We use smooth Fourier cutoffs that add up to one to split the
solution into faster and slower components. The two velocity
ranges overlap in the transition region.
For the faster component we use the
reduced Duhamel formula
\[
 i\int_t^T e^{i(t-s)\Delta}F(u(s))\,ds
 \rightharpoonup u(t)
 \quad\text{in }H^1_x(\R^2)\qquad(T\to\infty),
\]
proved in Proposition~\ref{prop:reduced_future_duhamel}.
For the slower component we use \eqref{eq:appendix_duhamel}
from a fixed initial time. In each chosen time direction, the classical rays with the selected
normal velocities do not connect the support of the interior
forcing term to the exterior observation region.Thus the
phase derivative in the normal frequency stays away from zero,
allowing integration by parts to give decay.

To control the interior contribution to the exterior Duhamel term,
we establish two propagation estimates.
Lemma~\ref{lem:all_frequency_gap_crossing_S0} gives a
$\delta^{-4}$ bound for the speed gap $\delta=B-A$, uniform for
$T_0\geq1$, which allows iteration through nested cones with
shrinking gaps.
Lemma~\ref{lem:normal_crossing_full_frequency} controls one normal
derivative of the Duhamel term: the interior contribution decays
as $t^{-2}$ in $L_x^2$ when the forcing term is bounded in
$L_t^\infty L_x^2$, and belongs to $L_{t,x}^4$ when the forcing
term is bounded in $L_t^\infty L_{x_e}^2H_{x_\perp}^{1/4}$.
The localized strong limits in the latter estimate identify these
contributions in the weak reduced Duhamel formula.
These estimates concern the free Schr\"odinger group on $\R^2$
and are not specific to the cubic--quintic NLS.

For the exterior forcing term on the intervals
$I_j=[2^jT_0,2^{j+1}T_0]$, the dispersive estimate gives the factor
$2^{-|j-k|/4}$ when $|j-k|\geq2$, while Strichartz estimates
control the remaining pairs. We iterate through nested cones, first removing
the remainder in $\ell^\infty$ and then summing in $\ell^4$.
Proposition~\ref{prop:conditional_exterior_global_L4} thus gives
\[
 \sum_{j\geq0}
 \|u\|_{L_{t,x}^4(\{(t,x):t\in I_j,\ |x|>ct\})}^4
 =\int_{T_0}^\infty\int_{|x|>ct}|u(t,x)|^4\,dx\,dt<\infty
\]
for every $c>0$ and sufficiently large $T_0$.
The uniform $H^1$ bound allows us to include any finite initial
interval, so we may fix $T_0=1$ for all $c>0$.

\medskip\noindent\textit{Gradient bounds and the final contradiction.}
The global exterior spacetime bound and the dyadic derivative
estimate yield
\[
 \nabla u\in L_{t,x}^4(\{(t,x):t\geq T_0,\ |x|>ct\})
 \qquad(c>0)
\]
by Proposition~\ref{prop:gradient_from_exterior_tail}.
Combining the $\ell^2$ nonlinear coefficients with the resulting
$\ell^4$ derivative bounds, we use
$\ell^2\cdot\ell^4\subset\ell^{4/3}$ to obtain
\[
 \zeta F(u)\in
 L_t^{4/3}([T_0,\infty);W_x^{1,4/3}(\R^2))
\]
for the smooth exterior cutoffs $\zeta$ in
Corollary~\ref{cor:global_exterior_N1_from_gradient}.
The corresponding Duhamel tails tend to zero, and the slower free
component leaves the exterior region. Two orthogonal normal
directions on each of finitely many angular sectors then give
\[
 \|\nabla u(t)\|_{L_x^2(\{x:|x|>ct\})}\longrightarrow0
 \qquad(t\to\infty,\ c>0).
\]
The uniform $H^1$ bound makes its squared time average tend to
zero, contradicting \eqref{eq:intro_exterior_energy_average}.

\subsection{Organization of the paper}
Section~\ref{sec:prelim} proves the compactness and scale
estimates. Section~\ref{sec:rigidity} establishes the virial limits
and the reduced Duhamel formula. Sections~\ref{sec:exterior_windows}
and~\ref{sec:gradient_rigidity} prove the exterior spacetime and
gradient estimates, respectively. Section~\ref{subsec:normal_cone_rigidity_closure}
completes the proof of Theorem~\ref{thm:threshold_scattering}.
Appendix~\ref{app:local_theory} collects the analytic tools, and
Appendix~\ref{app:propagation_proofs} proves the propagation estimates.

\section{Preliminaries and threshold compactness}
\label{sec:prelim}

\subsection{Notation}
\label{subsec:notation}
For nonnegative quantities $X$ and $Y$, we write $X\lesssim Y$
if $X\leq CY$ for some $C>0$,
$X\gtrsim Y$ if $Y\lesssim X$, and $X\sim Y$ if both hold.

We use the unitary Fourier transform
$\widehat f(\xi):=(2\pi)^{-1}\int_{\R^2}e^{-ix\cdot\xi}f(x)\,dx$.
Fix a real even function $\chi\in C_c^\infty(\R^2)$ with
$0\leq\chi\leq1$, $\chi=1$ on $B(0,1)$, and $\chi=0$ outside $B(0,2)$.
For $N>0$, define
\[
 \widehat{P_{\leq N}f}(\xi):=\chi(\xi/N)\widehat f(\xi),
 \qquad P_{>N}:=I-P_{\leq N}.
\]
For $N\in2^{\mathbb Z}$, set $P_N:=P_{\leq N}-P_{\leq N/2}$.
Then $f=\sum_{N\in2^{\mathbb Z}}P_Nf$ in $L^2$ for $f\in L^2_x(\R^2)$.

Write $D:=|\nabla|$, $J^s:=\langle D\rangle^s=(1-\Delta)^{s/2}$,
and $\mathbb S^1:=\{e\in\R^2:|e|=1\}$.
For $e\in\mathbb S^1$ and a unit vector $e_\perp$ orthogonal to $e$, set
\begin{equation}\label{eq:directional_notation}
 \begin{gathered}
 x_e:=x\cdot e,\qquad y_e:=y\cdot e,\qquad
 x_\perp:=x\cdot e_\perp,\qquad \xi_e:=\xi\cdot e,\\
 \partial_e:=e\cdot\nabla,\qquad D_e:=-i\partial_e.
 \end{gathered}
\end{equation}
For a smooth vector field $X:\R^2\to\R^2$,
$\nabla X:=(\partial_jX_i)_{i,j=1}^2$ is its Jacobian matrix;
matrix norms are Euclidean operator norms.

For a time interval $I\subset\R$ and $0\leq s\leq1$, define
\[
 \begin{aligned}
 S^s(I)&:=L^\infty_t(I;H^s_x(\R^2))
          \cap L^4_t(I;W^{s,4}_x(\R^2))
          \cap L^8_t(I;W^{s,8/3}_x(\R^2)),\\
 N^s(I)&:=L^{4/3}_t(I;W^{s,4/3}_x(\R^2)).
 \end{aligned}
\]
We equip $S^s(I)$ with the sum norm. A Schr\"odinger-admissible pair
satisfies $2/q+2/r=1$, $2\leq q,r\leq\infty$, and $(q,r)\ne(2,\infty)$.
Each pair is fixed before a time limit is taken.
For the minimal-mass nonscattering solution $u$, set
\begin{equation}\label{eq:uniform_H1_bound}
 B_H:=\sup_{t\in\R}\|u(t)\|_{H^1_x}<\infty.
\end{equation}
Finiteness follows from mass and energy conservation and
\eqref{eq:coercive_energy}. For an arbitrary solution on an interval $I$,
$B_H$ denotes the bound $\sup_{t\in I}\|u(t)\|_{H^1}$.
The analytic estimates and local and stability theories used below
are collected in Appendix~\ref{app:local_theory}.

\subsection{Conservation laws and scattering criterion}
\label{subsec:conservation}
The mass, energy, and momentum in
\eqref{eq:intro_conserved_quantities} are conserved for $H^1$ solutions
by Lemma~\ref{lem:local_small_data}. The energy of \eqref{eq:cubic_nls} is
\[
 \Ec(f):=\frac12\norm{\nabla f}_{L^2_x(\R^2)}^2
          -\frac14\norm f_{L^4_x(\R^2)}^4
 \qquad(f\in H^1(\R^2)).
\]
When $M(u)=m_Q$, the sharp Gagliardo--Nirenberg inequality
\eqref{eq:GN} implies
\begin{equation}\label{eq:positive}
 E(u)=\Ec(u)+\frac16\norm u_{L^6_x(\R^2)}^6\geq\frac16\norm u_{L^6_x(\R^2)}^6>0.
\end{equation}
Moreover, H\"older's and Young's inequalities give
\begin{equation}\label{eq:coercive_energy}
 \norm u_{L^4_x(\R^2)}^4\leq\sqrt{M(u)}\norm u_{L^6_x(\R^2)}^3,\qquad
 E(u)\geq\frac12\norm{\nabla u}_{L^2_x(\R^2)}^2-\frac{3M(u)}{32}.
\end{equation}
Conservation of mass and energy then gives
$\norm{\nabla u(t)}_{L^2_x(\R^2)}^2\leq2E(u)+3M(u)/16$.
The continuation criterion in Lemma~\ref{lem:local_small_data}
gives global existence and a uniform $H^1$ bound.

The Galilean transform $u_a$ defined in the introduction satisfies,
for $a\in\R^2$,
\[
 P(u_a)=P(u)+aM(u),\qquad
 E(u_a)=E(u)+a\cdot P(u)+\frac12|a|^2M(u).
\]
When $M(u)=m_Q$, choosing $a=-P(u)/M(u)$ gives $P(u_a)=0$.
Since $M(u_a)=m_Q$,
\eqref{eq:GN} implies $E(u_a)=E(u)-|P(u)|^2/(2M(u))>0$.
In what follows, we use the normalized conserved energy $E:=E(u_a)$.

The scattering criterion in Lemma~\ref{lem:H1_scattering_criterion}
(Appendix~\ref{app:local_theory}) gives $H^1$ scattering in either
time direction from a uniform $H^1$ bound and a finite $L^4_{t,x}$ norm
in that direction.

\subsection{Nonlinear approximation for large-scale data}
\label{sec:compactness}

Proposition~\ref{prop:largescale} gives uniform $L^4_{t,x}$ bounds for
solutions of \eqref{eq:cq} with large-scale initial data. In the next
subsection, we use it to establish the Palais--Smale condition and
classify large-scale limits. The resulting uniform spatial localization
of mass and spatial-scale bounds enter the virial argument.

\begin{proposition}
\label{prop:largescale}
Suppose $h_n\to\infty$, $\sup_n|\xi_n|<\infty$,
$\sup_n\norm{f_n}_{H^1_x(\R^2)}<\infty$, and, for some $\phi\in L^2(\R^2)$,
\[
 f_n(x)=e^{i\gamma_n+ix\cdot\xi_n}h_n^{-1}
 \phi((x-x_n)/h_n)+o_{L^2_x(\R^2)}(1),
\]
where $\gamma_n\in\R$ are constant phase parameters and
$x_n,\xi_n\in\R^2$ are translation and frequency modulation parameters, respectively.
Let $w$ solve \eqref{eq:cubic_nls} with $w(0)=\phi$.
If $J\ni0$ and
$\norm w_{L^4_{s,y}(J\times\R^2)}<\infty$, then the solutions $v_n$
of \eqref{eq:cq} with $v_n(0)=f_n$ satisfy
\[
 \limsup_{n\to\infty}\norm{v_n}_{L^4_{t,x}(h_n^2J\times\R^2)}<\infty.
\]
\end{proposition}

\begin{proof}
Fix $0<\theta<1$, set $I_n:=h_n^2J$, and let $w_n$ solve
\eqref{eq:cubic_nls} with $w_n(0)=P_{\le h_n^\theta}\phi$.
The convergence $P_{\le h_n^\theta}\phi\to\phi$ in $L^2$ and
$L^2$ stability give uniform $L^4_{s,y}$ bounds for $w_n$ on $J$.
The persistence of regularity in Lemma~\ref{lem:persistence_regularity}
then gives, for all large $n$ and $0\le s\le1$,
\begin{equation}\label{eq:cubicpersistent}
 \norm{D^sw_n}_{L^\infty_s(J;L^2_y(\R^2))\cap L^4_{s,y}(J\times\R^2)
                         \cap L^8_s(J;L^{8/3}_y(\R^2))}
 \le C_J h_n^{s\theta}.
\end{equation}
Here $C_J$ depends only on $\norm\phi_{L^2_y(\R^2)}$ and
$\norm w_{L^4_{s,y}(J\times\R^2)}$, so the same bound holds on
compact subintervals of $J$ containing zero.

We next construct the approximate solution and check the hypotheses
of Lemma~\ref{lem:stability}. Define
\[
 \widetilde v_n(t,x):=e^{i\gamma_n+ix\cdot\xi_n-it|\xi_n|^2}h_n^{-1}
 w_n\!\left(t/h_n^2,(x-x_n-2t\xi_n)/h_n\right).
\]
It solves \eqref{eq:cubic_nls} in $(t,x)$, with equation error in \eqref{eq:cq}
\[
 e_n:=(i\partial_t+\Delta)\widetilde v_n-F(\widetilde v_n)=-|\widetilde v_n|^4\widetilde v_n.
\]
By the Sobolev estimate \eqref{eq:appendix_sobolev} and \eqref{eq:cubicpersistent},
\[
 \|w_n\|_{L^8_{s,y}(J\times\R^2)}
 \lesssim\|D^{1/2}w_n\|_{L^8_s(J;L^{8/3}_y(\R^2))}
 \lesssim_J h_n^{\theta/2}.
\]
H\"older's inequality and the fractional product estimate
\eqref{eq:appendix_fractional_product} then give, for $s=0,\frac12$,
\begin{align}
 &\norm{D^s\{h_n^{-5}(|w_n|^4w_n)
       (t/h_n^2,x/h_n)\}}_{L^{4/3}_{t,x}(I_n\times\R^2)}
 \notag\\
 &\qquad\lesssim h_n^{-2-s}\norm{w_n}_{L^8_{s,y}(J\times\R^2)}^4
                  \norm{D^sw_n}_{L^4_{s,y}(J\times\R^2)}
 \lesssim_J h_n^{-(2+s)(1-\theta)}. \label{eq:scaledquintic}
\end{align}
Spatial translations preserve these norms at each time, and
\[
 \norm{D^{1/2}(e^{ix\cdot\xi}g)}_{L^p_x(\R^2)}
 \lesssim\norm{D^{1/2}g}_{L^p_x(\R^2)}+|\xi|^{1/2}\norm g_{L^p_x(\R^2)},
 \qquad 1<p<\infty.
\]
Thus \eqref{eq:scaledquintic} and $\sup_n|\xi_n|<\infty$ imply
\begin{equation}\label{eq:boostedquintic}
 \norm{\langle D\rangle^{1/2}e_n}_{L^{4/3}_{t,x}(I_n\times\R^2)}
 \lesssim_{J,\sup_n|\xi_n|}
 h_n^{-2(1-\theta)}+h_n^{-\frac52(1-\theta)}\longrightarrow0.
\end{equation}
The same scaling and modulation estimates applied to
\eqref{eq:cubicpersistent} bound
$\norm{\widetilde v_n}_{L^\infty_t(I_n;H^1_x(\R^2))}$ and
$\norm{\langle D\rangle^{1/2}\widetilde v_n}_{L^4_{t,x}(I_n\times\R^2)}$
uniformly in $n$. Moreover,
\[
 \|f_n-\widetilde v_n(0)\|_{L^2_x(\R^2)}\longrightarrow0,
 \qquad
 \sup_n\|f_n-\widetilde v_n(0)\|_{H^1_x(\R^2)}<\infty.
\]
By the interpolation inequality \eqref{eq:appendix_bessel_interpolation} and
the Strichartz estimate \eqref{eq:appendix_strichartz},
\[
 \begin{aligned}
 \bigl\|e^{it\Delta}(f_n-\widetilde v_n(0))\bigr\|_{L^4_t(I_n;W^{1/2,4}_x(\R^2))}
 &\leq C\|f_n-\widetilde v_n(0)\|_{H^{1/2}_x(\R^2)}\\
 &\lesssim\|f_n-\widetilde v_n(0)\|_{L^2_x(\R^2)}^{1/2}
          \|f_n-\widetilde v_n(0)\|_{H^1_x(\R^2)}^{1/2}\longrightarrow0,
 \end{aligned}
\]
with $C$ independent of $I_n$.
Conservation of mass and energy, \eqref{eq:coercive_energy}, and
the uniform $H^1$ bound on $f_n$ give global solutions $v_n$ with
uniform $L^\infty_t(\R;H^1_x(\R^2))$ bounds.
Applying Lemma~\ref{lem:stability} on $I_n$, we obtain
\[
 \|v_n-\widetilde v_n\|_{L^4_{t,x}(I_n\times\R^2)}\longrightarrow0.
\]
The uniform $L^4_{t,x}(I_n\times\R^2)$ bound on
$\widetilde v_n$ completes the proof.
\end{proof}

\begin{remark}
\label{rem:largescale_scaling}
Here ``large-scale'' means $h_n\to\infty$ for the spatial scaling
parameters $h_n$ in Proposition~\ref{prop:largescale}.
Under the mass-preserving rescaling
$u(t,x)=h_n^{-1}v(t/h_n^2,x/h_n)$, equation~\eqref{eq:cq} becomes
\[
 i\partial_sv+\Delta_yv=-|v|^2v+h_n^{-2}|v|^4v.
\]
The coefficient $h_n^{-2}\to0$ makes the quintic term a perturbation
of \eqref{eq:cubic_nls}.
\end{remark}

\begin{remark}
\label{rem:largescale_convergence}
The solutions $v_n$ also have uniform $L^8_{t,x}(I_n\times\R^2)$ bounds,
and $\|v_n-\widetilde v_n\|_{L^\infty_t(I_n;L^2_x(\R^2))}\to0$.
For $z=v_n$ or $\widetilde v_n$, the Gagliardo--Nirenberg inequality \eqref{eq:appendix_GN_L8} gives
\[
 \|z\|_{L^8_{t,x}(I_n\times\R^2)}^8
 \leq C\|\nabla z\|_{L^\infty_t(I_n;L^2_x(\R^2))}^4
          \|z\|_{L^4_{t,x}(I_n\times\R^2)}^4
 \leq C.
\]
Put $z_n:=v_n-\widetilde v_n$. H\"older's inequality gives
\[
 \begin{aligned}
 &\|F(v_n)-F(\widetilde v_n)\|_{L^{4/3}_{t,x}(I_n\times\R^2)}\\
 &\quad\leq C\sum_{z\in\{v_n,\widetilde v_n\}}
 \Bigl(\|z\|_{L^4_{t,x}(I_n\times\R^2)}^2
       +\|z\|_{L^8_{t,x}(I_n\times\R^2)}^4\Bigr)
       \|z_n\|_{L^4_{t,x}(I_n\times\R^2)}\\
 &\quad\leq C\|z_n\|_{L^4_{t,x}(I_n\times\R^2)}.
 \end{aligned}
\]
The difference equation and the Strichartz estimate \eqref{eq:appendix_strichartz} therefore imply
\[
 \|z_n\|_{L^\infty_t(I_n;L^2_x(\R^2))}
 \leq C\Bigl(\|z_n(0)\|_{L^2_x(\R^2)}
       +\|e_n\|_{L^{4/3}_{t,x}(I_n\times\R^2)}
       +\|z_n\|_{L^4_{t,x}(I_n\times\R^2)}\Bigr)\longrightarrow0.
\]
\end{remark}

\subsection{The Palais--Smale condition and soliton limits}
\label{subsec:cubic_limits}

We give the compactness and rigidity arguments in forward time;
the backward-time statements follow by applying them to
$\overline{u(-t,x)}$.

To prove uniform spatial localization of mass for forward
nonscattering solutions of \eqref{eq:cq} with mass $m_Q$, we apply
the linear profile decomposition to $u(T_n)$, where $T_n\to\infty$.
If every nonzero profile had mass below $m_Q$, then $u$ would scatter,
a contradiction. Thus one profile has mass $m_Q$; the total mass
excludes other nonzero profiles and forces the remainder to vanish
in $L^2$. We then show that the normalized time parameter has a
finite limit, giving the Palais--Smale condition below.

The linear profile decomposition in $L^2(\R^2)$ was established
by Merle and Vega \cite{MerleVega}. For $H^1$-bounded sequences,
we combine the restrictions on scale and frequency modulation of Carles and Keraani
\cite{CarlesKeraani} with the frequency-truncated construction of
Cheng, Guo, and Zhao \cite{CGZ}.

\begin{lemma}[Linear profile decomposition]
\label{lem:linear_profile_decomposition}
Let $(f_n)$ be bounded in $H^1(\R^2)$ and fix $0<\theta<1$.
After passing to a subsequence, there exist profiles and parameters
such that, for every finite $J$,
\begin{equation}\label{eq:nonradialprofiles}
 f_n=\sum_{j\le J}a_n^j+r_n^J,\qquad
 a_n^j:=e^{i\gamma_n^j+ix\cdot\xi_n^j}
 e^{-it_n^j\Delta}
 \left[(h_n^j)^{-1}(P_n^j\phi^j)
              ((\,\cdot-x_n^j)/h_n^j)\right].
\end{equation}
Here $\phi^j\in L^2$, $r_n^J\in H^1$,
$\gamma_n^j,t_n^j\in\R$, $x_n^j,\xi_n^j\in\R^2$,
$h_n^j>0$, and $|\xi_n^j|\leq C_j$.
For each $j$, one of the following holds:
\begin{itemize}
 \item $h_n^j\equiv1$, $P_n^j:=\operatorname{Id}$,
       $\xi_n^j=0$, and $\phi^j\in H^1$;
 \item $h_n^j\to\infty$ and
       $P_n^j:=P_{\leq(h_n^j)^\theta}$.
\end{itemize}
For every $j$, the normalized time parameter
$\sigma_n^j:=-t_n^j/(h_n^j)^2$ converges in $[-\infty,\infty]$, and the parameters satisfy the orthogonality
\begin{equation}\label{eq:profile_parameter_orthogonality}
 \begin{aligned}
 &\frac{h_n^j}{h_n^k}+\frac{h_n^k}{h_n^j}
 +\frac{|t_n^j-t_n^k|}{(h_n^j)^2}
 +h_n^j|\xi_n^j-\xi_n^k|\\
 &\qquad+\frac{|x_n^j-x_n^k+2t_n^j(\xi_n^j-\xi_n^k)|}{h_n^j}
 \longrightarrow\infty\qquad(j\ne k).
 \end{aligned}
\end{equation}
The initial profiles are asymptotically orthogonal in $L^2$ and
in the gradient norm:
\begin{equation}\label{eq:profile_initial_orthogonality}
 \langle D^sa_n^j,D^sa_n^k\rangle_{L^2_x(\R^2)}\longrightarrow0
 \qquad(s=0,1,\ j\ne k).
\end{equation}
For every fixed $J$, we have the expansions
\begin{equation}\label{eq:profile_pythagorean}
 \|D^sf_n\|_{L^2_x(\R^2)}^2
 =\sum_{j\leq J}\|D^sa_n^j\|_{L^2_x(\R^2)}^2
   +\|D^sr_n^J\|_{L^2_x(\R^2)}^2+o_n(1),\qquad s=0,1.
\end{equation}
Moreover, setting $z_n^J(t):=e^{it\Delta}r_n^J$, we have
\begin{equation}\label{eq:linear_profile_remainder}
 \lim_{J\to\infty}\limsup_{n\to\infty}
 \left(\|z_n^J\|_{L^4_t(\R;W^{1/2,4}_x(\R^2))}
       +\|z_n^J\|_{L^8_t(\R;W^{1/2,8/3}_x(\R^2))}\right)=0.
\end{equation}
\end{lemma}

The Palais--Smale condition states that every sequence $T_n\to\infty$
has a subsequence along which $u(T_n)$ converges strongly in $L^2$
modulo phase rotations, translations, frequency modulations, and
mass-preserving dilations.

\begin{lemma}[Palais--Smale condition]
\label{lem:single_profile_compactness}
Let $u$ be an $H^1$ solution of \eqref{eq:cq} satisfying
\[
 M(u)=m_Q,\qquad
 \norm u_{L^4_{t,x}([0,\infty)\times\R^2)}=\infty.
\]
For every sequence $T_n\to\infty$, after passing to a subsequence,
there exist $\phi\in L^2(\R^2)$, $\gamma_n\in\R$,
$x_n,\xi_n\in\R^2$, and $h_n>0$, with
$\sup_n|\xi_n|<\infty$, such that
\begin{equation}\label{eq:singlefullprofile}
 u(T_n,x)=e^{i\gamma_n+ix\cdot\xi_n}h_n^{-1}
             \phi((x-x_n)/h_n)+o_{L^2_x(\R^2)}(1),\qquad M(\phi)=m_Q.
\end{equation}
Either $h_n\equiv1$, in which case $\xi_n=0$ and $\phi\in H^1$,
or $h_n\to\infty$.
\end{lemma}

\begin{proof}
By conservation of mass and energy and \eqref{eq:coercive_energy},
the solution is global and has the uniform $H^1$ bound
\eqref{eq:uniform_H1_bound}.

Fix a sequence $T_n\to\infty$ and set $f_n:=u(T_n)$.
Choose $0<\theta<1$ and apply
Lemma~\ref{lem:linear_profile_decomposition} to $(f_n)$.
The decomposition \eqref{eq:nonradialprofiles} satisfies
\[
 \sum_jM(\phi^j)\leq m_Q,\qquad
 \sup_J\limsup_{n\to\infty}\|r_n^J\|_{H^1_x(\R^2)}\leq B_H.
\]

\par\medskip\Needspace{3\baselineskip}\noindent\textbf{Only one nonzero profile.}\par\nobreak\smallskip
\noindent
We first show that at least one profile is nonzero. Otherwise,
all profiles are zero and \eqref{eq:linear_profile_remainder} gives
\[
 \|e^{it\Delta}f_n\|_{L^4_t(\R;W^{1/2,4}_x(\R^2))}\to0.
\]
Applying Lemma~\ref{lem:stability} with $\widetilde v\equiv0$,
$e\equiv0$, $t_0=0$, and initial data $f_n$ gives a finite
$L^4_{t,x}$ norm for $u$, which yields a contradiction.

Suppose that $M(\phi^j)<m_Q$ for every $j$, and let $U_n^j$ be
the solution of \eqref{eq:cq} with $U_n^j(0)=a_n^j$.
For each fixed $j$, we obtain global spacetime bounds uniform in $n$.
Write $\sigma^j:=\lim_n\sigma_n^j$.

\par\smallskip\Needspace{3\baselineskip}\noindent
\text{Case 1: $h_n^j\equiv1$.}\par\nobreak\smallskip\noindent
Let $v^j$ solve \eqref{eq:cq} with
$v^j(\sigma^j)=e^{i\sigma^j\Delta}\phi^j$ if $\sigma^j\in\R$,
or with $H^1$ scattering state $\phi^j$ as $t\to\sigma^j$ if
$\sigma^j=\pm\infty$. By the subthreshold scattering theorem
\cite[Theorem~1.2]{Cheng} and
Lemma~\ref{lem:H1_scattering_criterion}, we have
$\|v^j\|_{S^1(\R)}<\infty$.
Continuity in $H^1$ when $\sigma^j\in\R$, and scattering when
$\sigma^j=\pm\infty$, give
\[
 \|a_n^j-e^{i\gamma_n^j}v^j(\sigma_n^j,\,\cdot-x_n^j)\|_{H^1_x(\R^2)}
 =\|e^{i\sigma_n^j\Delta}\phi^j-v^j(\sigma_n^j)\|_{H^1_x(\R^2)}
 \longrightarrow0.
\]
Apply Lemma~\ref{lem:stability} with approximate solution
$e^{i\gamma_n^j}v^j(t+\sigma_n^j,x-x_n^j)$, zero equation error,
and initial datum $a_n^j$.
The Strichartz estimate \eqref{eq:appendix_strichartz} controls the
free evolution of the initial difference by its $H^1$ norm.
On each interval $[-T,T]$, the stability constants depend only on
the global bounds for $v^j$, and are independent of $n$ and $T$.
Uniqueness identifies the exact solution with $U_n^j$.
Letting $T\to\infty$ and then $n\to\infty$, we obtain
\begin{equation}\label{eq:unit_profile_convergence}
 \|U_n^j-e^{i\gamma_n^j}v^j(t+\sigma_n^j,x-x_n^j)\|_{L^4_{t,x}(\R\times\R^2)}
 \longrightarrow0.
\end{equation}

\par\smallskip\Needspace{3\baselineskip}\noindent
\text{Case 2: $h_n^j\to\infty$.}\par\nobreak\smallskip\noindent
If $\sigma^j\in\R$, let $w^j$ and $w_n^j$ solve
\eqref{eq:cubic_nls} with
\[
 w^j(\sigma^j)=e^{i\sigma^j\Delta}\phi^j,
 \qquad w_n^j(\sigma^j)=P_{\leq(h_n^j)^\theta}w^j(\sigma^j).
\]
If $\sigma^j=\pm\infty$, choose $w^j$ and $w_n^j$ with scattering
states $\phi^j$ and $P_n^j\phi^j$ in that time direction.
In both cases, Lemma~\ref{lem:cubic_subthreshold} and the $L^2$
stability theory for \eqref{eq:cubic_nls} give
$w_n^j\to w^j$ in $L^\infty_s(\R;L^2_y(\R^2))\cap L^4_{s,y}(\R\times\R^2)$.
Using continuity when $\sigma^j\in\R$, and scattering when
$\sigma^j=\pm\infty$, we obtain
\begin{equation}\label{eq:large_profile_initial_discrepancy}
\begin{aligned}
 \|w_n^j(\sigma_n^j)-e^{i\sigma_n^j\Delta}P_n^j\phi^j\|_{L^2_y(\R^2)}
 &\leq\|w_n^j-w^j\|_{L^\infty_s(\R;L^2_y(\R^2))}
       +\|w^j(\sigma_n^j)-e^{i\sigma_n^j\Delta}\phi^j\|_{L^2_y(\R^2)}\\
 &\quad+\|(1-P_n^j)\phi^j\|_{L^2_y(\R^2)}\longrightarrow0.
\end{aligned}
\end{equation}
Mass conservation, $L^2$ stability, and
Lemma~\ref{lem:persistence_regularity} give
\begin{equation}\label{eq:cubic_profile_bounds}
 \|w_n^j\|_{L^\infty_s(\R;L^2_y(\R^2))}
 +\|w_n^j\|_{L^4_{s,y}(\R\times\R^2)}
 +(h_n^j)^{-\theta}\|\nabla w_n^j\|_{L^\infty_s(\R;L^2_y(\R^2))}
 \lesssim_j1.
\end{equation}

To bound $U_n^j$ by Lemma~\ref{lem:stability}, define the approximate
solution obtained from \eqref{eq:cubic_nls}:
\begin{equation}\label{eq:large_profile_approximation}
 \widetilde U_n^j(t,x)
 :=e^{i\gamma_n^j+ix\cdot\xi_n^j-it|\xi_n^j|^2}(h_n^j)^{-1}
 w_n^j\!\left(\frac{t}{(h_n^j)^2}+\sigma_n^j,
              \frac{x-x_n^j-2t\xi_n^j}{h_n^j}\right).
\end{equation}
By \eqref{eq:large_profile_initial_discrepancy},
$a_n^j-\widetilde U_n^j(0)$ tends to zero in $L^2$. Scaling,
\eqref{eq:cubic_profile_bounds}, and $|\xi_n^j|\leq C_j$ give
\[
 \|\nabla(a_n^j-\widetilde U_n^j(0))\|_{L^2_x(\R^2)}
 \lesssim_j(h_n^j)^{\theta-1}
       +\|a_n^j-\widetilde U_n^j(0)\|_{L^2_x(\R^2)}\longrightarrow0.
\]
Thus the initial-data difference vanishes in $H^1$.
Applying \eqref{eq:cubicpersistent} and \eqref{eq:boostedquintic}
to $w_n^j$ on $\R$ gives
\[
\begin{aligned}
 &\|\widetilde U_n^j\|_{L^\infty_t(\R;H^1_x(\R^2))}
  +\|\widetilde U_n^j\|_{L^4_t(\R;W^{1/2,4}_x(\R^2))}\lesssim_j1,\\
 &\bigl\||\widetilde U_n^j|^4\widetilde U_n^j
       \bigr\|_{L^{4/3}_t(\R;W^{1/2,4/3}_x(\R^2))}\longrightarrow0.
\end{aligned}
\]
Since the equation error is $-|\widetilde U_n^j|^4\widetilde U_n^j$,
applying Lemma~\ref{lem:stability} on both time halflines gives
\begin{equation}\label{eq:large_profile_convergence}
 \|U_n^j-\widetilde U_n^j\|_{L^4_{t,x}(\R\times\R^2)}\to0,
 \qquad
 \limsup_{n\to\infty}\|U_n^j\|_{L^4_{t,x}(\R\times\R^2)}<\infty.
\end{equation}

In both cases, conservation of mass and energy,
\eqref{eq:coercive_energy}, and the uniform $L^4_{t,x}$ bound give,
by Lemma~\ref{lem:H1_scattering_criterion}, the first estimate in
\eqref{eq:nonlinear_profile_bounds}. The second follows from
\eqref{eq:cubic_profile_bounds} and
\eqref{eq:large_profile_approximation}:
\begin{equation}\label{eq:nonlinear_profile_bounds}
\begin{aligned}
 \limsup_{n\to\infty}\|U_n^j\|_{S^1(\R)}&<\infty,\\
 \limsup_{n\to\infty}\|\widetilde U_n^j\|_{L^\infty_t(\R;H^1_x(\R^2))}
 &<\infty\qquad(h_n^j\to\infty).
\end{aligned}
\end{equation}
To apply Lemma~\ref{lem:stability} with initial data $f_n$, we need
bounds independent of the number of nonlinear profiles.
After passing to a diagonal subsequence, set
$\alpha_j:=\lim_{n\to\infty}\norm{a_n^j}_{H^1_x(\R^2)}$.
Then \eqref{eq:profile_pythagorean} gives
$\sum_j\alpha_j^2\leq B_H^2$. Choose $J_0$ so that $\alpha_j$ is less
than half the threshold $\delta_0$ in Lemma~\ref{lem:local_small_data}
for $j>J_0$. For fixed $J>J_0$ and sufficiently large $n$,
the small-data correction estimate \eqref{eq:appendix_small_data_correction} gives
\begin{equation}\label{eq:small_profile_correction}
 \norm{U_n^j-e^{it\Delta}a_n^j}_{S^1(\R)}
 \lesssim\norm{a_n^j}_{H^1_x(\R^2)}^3\qquad(J_0<j\leq J).
\end{equation}
Define the approximate solution
\[
 V_n^J:=\sum_{j\leq J}U_n^j+e^{it\Delta}r_n^J,
 \qquad V_n^J(0)=f_n.
\]
For this fixed $J_0$, the identity
\[
 V_n^J=\sum_{j\le J_0}U_n^j+e^{it\Delta}r_n^{J_0}
       +\sum_{J_0<j\le J}(U_n^j-e^{it\Delta}a_n^j)
\]
and $\sum_{j>J_0}\alpha_j^3<\infty$ imply
\begin{equation}\label{eq:profile_sum_bound}
 \sup_J\limsup_{n\to\infty}\norm{V_n^J}_{S^1(\R)}<\infty.
\end{equation}
To estimate the equation error, set
\[
 e_n^J:=(i\partial_t+\Delta)V_n^J-F(V_n^J)
       =\sum_{j\leq J}F(U_n^j)-F(V_n^J).
\]
By H\"older's inequality,
\[
 \|F(v)\|_{N^1(\R)}
 \leq C\bigl(\|v\|_{L^4_{t,x}(\R\times\R^2)}^2+\|v\|_{L^8_{t,x}(\R\times\R^2)}^4\bigr)
             \|v\|_{L^4_t(\R;W^{1,4}_x(\R^2))}.
\]
Sobolev embedding controls the two spacetime norms by
$\|v\|_{S^1(\R)}$. The small-data estimate \eqref{eq:appendix_small_data} therefore gives
\begin{equation}\label{eq:profile_nonlinearity_tail}
 \sup_{J>J_0}\limsup_{n\to\infty}
 \sum_{J_0<j\leq J}\|F(U_n^j)\|_{N^1(\R)}
 \lesssim\sum_{j>J_0}(\alpha_j^3+\alpha_j^5).
\end{equation}
Combining \eqref{eq:profile_nonlinearity_tail} with
\eqref{eq:profile_sum_bound} and
\eqref{eq:nonlinear_profile_bounds} for $j\leq J_0$ yields
\begin{equation}\label{eq:profile_error_bound}
 \sup_J\limsup_{n\to\infty}\|e_n^J\|_{N^1(\R)}<\infty.
\end{equation}
We claim that
\begin{equation}\label{eq:profile_error_smallness}
 \lim_{J\to\infty}\limsup_{n\to\infty}
 \|e_n^J\|_{N^0(\R)}=0.
\end{equation}
For $j>J_0$, split the profile sum into its free evolution and
the corrections in \eqref{eq:small_profile_correction}.
The Strichartz estimate \eqref{eq:appendix_strichartz} and
orthogonality \eqref{eq:profile_initial_orthogonality} give
\begin{equation}\label{eq:nonlinear_profile_tail}
\begin{aligned}
 &\sup_{J>J_0}\limsup_{n\to\infty}
 \left\|\sum_{J_0<j\leq J}U_n^j\right\|_{S^1(\R)}\\
 &\quad\lesssim
   \left(\sum_{j>J_0}\alpha_j^2\right)^{1/2}
       +\sum_{j>J_0}\alpha_j^3
       \longrightarrow0\qquad(J_0\to\infty).
\end{aligned}
\end{equation}
For fixed $J_0$ and $J>J_0$, set
$S_n^{J_0}:=\sum_{j\leq J_0}U_n^j$ and write
\begin{equation}\label{eq:profile_error_decomposition}
\begin{aligned}
 e_n^J
 &=\left[\sum_{j\leq J_0}F(U_n^j)-F(S_n^{J_0})\right]
   +\sum_{J_0<j\leq J}F(U_n^j)\\
 &\quad+F(S_n^{J_0})-F(V_n^J).
\end{aligned}
\end{equation}
The second term is controlled by \eqref{eq:profile_nonlinearity_tail}.
Each term of $F(S_n^{J_0})-F(V_n^J)$ contains
$V_n^J-S_n^{J_0}=\sum_{J_0<j\leq J}U_n^j+e^{it\Delta}r_n^J$
or its conjugate, to be controlled by the two remainder estimates.
We first treat the bracketed term, a finite sum of cubic and
quintic mixed terms of the form
\[
 U_n^{j_1}\overline{U_n^{j_2}}U_n^{j_3},\qquad
 U_n^{j_1}\overline{U_n^{j_2}}U_n^{j_3}
             \overline{U_n^{j_4}}U_n^{j_5},
\]
where $j_\ell\leq J_0$ and the indices are not all equal.

For the cubic mixed terms, it suffices to prove
\begin{equation}\label{eq:profile_product_limit}
 \|U_n^jU_n^k\|_{L^2_{t,x}(\R\times\R^2)}\longrightarrow0
 \qquad(j\ne k),
\end{equation}
and control the remaining factors in $L^4_{t,x}(\R\times\R^2)$.
In $L^4_{s,y}(\R\times\R^2)$, approximate $v^j$ when $h_n^j\equiv1$
and $w^j$ when $h_n^j\to\infty$ by bounded smooth functions
$\Psi^j\in C_c^\infty(\R\times\R^2)$. Set
\[
 \Psi_n^j(t,x):=e^{i\gamma_n^j+ix\cdot\xi_n^j-it|\xi_n^j|^2}
 (h_n^j)^{-1}\Psi^j\!\left(\frac{t-t_n^j}{(h_n^j)^2},
                         \frac{x-x_n^j-2t\xi_n^j}{h_n^j}\right).
\]
By \eqref{eq:unit_profile_convergence}, the $L^4$ convergence
$w_n^j\to w^j$, and \eqref{eq:large_profile_convergence}, for every
$0<\eta<1$ we may choose these functions so that
\begin{equation}\label{eq:profile_smooth_approximation}
 \limsup_{n\to\infty}
 \|U_n^j-\Psi_n^j\|_{L^4_{t,x}(\R\times\R^2)}<\eta
\end{equation}
for all indices under consideration. It suffices to prove
\eqref{eq:profile_product_limit} on a subsequence. By diagonal
extraction over $j\ne k$, assume that each nonnegative term in
\eqref{eq:profile_parameter_orthogonality} converges in $[0,\infty]$,
and retain this subsequence below.
If $h_n^j/h_n^k\to0$ or $\infty$, boundedness and compact support give
\[
 \|\Psi_n^j\Psi_n^k\|_{L^2_{t,x}(\R\times\R^2)}^2
 \lesssim\left(\frac{\min(h_n^j,h_n^k)}
                        {\max(h_n^j,h_n^k)}\right)^2\longrightarrow0.
\]
Otherwise, $h_n^j/h_n^k\to\rho\in(0,\infty)$. If
$|t_n^j-t_n^k|/(h_n^j)^2\to\infty$, the time supports of
$\Psi_n^j$ and $\Psi_n^k$ are disjoint for all sufficiently large $n$.
If $|t_n^j-t_n^k|/(h_n^j)^2$ is bounded and
$h_n^j|\xi_n^j-\xi_n^k|\to\infty$, write
$t=t_n^j+(h_n^j)^2s$. On the common spatial support,
\[
 \left|\frac{x_n^j-x_n^k+2t_n^j(\xi_n^j-\xi_n^k)}{h_n^j}
       +2h_n^js(\xi_n^j-\xi_n^k)\right|\lesssim1.
\]
The set of such $s$ has length
$O((h_n^j|\xi_n^j-\xi_n^k|)^{-1})$, so the product norm tends to zero.
If both $|t_n^j-t_n^k|/(h_n^j)^2$ and
$h_n^j|\xi_n^j-\xi_n^k|$ are bounded,
\eqref{eq:profile_parameter_orthogonality} gives
\[
 \frac{|x_n^j-x_n^k+2t_n^j(\xi_n^j-\xi_n^k)|}{h_n^j}
 \longrightarrow\infty,
\]
and the spatial supports are disjoint for all sufficiently large $n$.
Thus, for each fixed choice of approximants,
\begin{equation}\label{eq:smooth_profile_product_limit}
 \|\Psi_n^j\Psi_n^k\|_{L^2_{t,x}(\R\times\R^2)}\longrightarrow0.
\end{equation}
By H\"older's inequality,
\begin{equation}\label{eq:profile_product_comparison}
\begin{aligned}
 \|U_n^jU_n^k\|_{L^2_{t,x}(\R\times\R^2)}
 &\leq C_{j,k}\Bigl(
       \|U_n^j-\Psi_n^j\|_{L^4_{t,x}(\R\times\R^2)}
       +\|U_n^k-\Psi_n^k\|_{L^4_{t,x}(\R\times\R^2)}\Bigr)\\
 &\quad+\|\Psi_n^j\Psi_n^k\|_{L^2_{t,x}(\R\times\R^2)}.
\end{aligned}
\end{equation}
By \eqref{eq:nonlinear_profile_bounds} and
\eqref{eq:profile_smooth_approximation}, $C_{j,k}$ is independent of
$\eta$ for sufficiently large $n$. Taking $\limsup_{n\to\infty}$ in
\eqref{eq:profile_product_comparison}, then using
\eqref{eq:profile_smooth_approximation} and
\eqref{eq:smooth_profile_product_limit}, gives the bound $C_{j,k}\eta$.
Letting $\eta\downarrow0$ proves \eqref{eq:profile_product_limit}.

For the quintic mixed terms, first suppose $h_n^j\to\infty$.
Scaling \eqref{eq:large_profile_approximation} and applying the
Gagliardo--Nirenberg inequality and \eqref{eq:cubic_profile_bounds} give
\begin{equation}\label{eq:large_profile_approximation_L8}
 \|\widetilde U_n^j\|_{L^8_{t,x}(\R\times\R^2)}^8
 =(h_n^j)^{-4}\|w_n^j\|_{L^8_{s,y}(\R\times\R^2)}^8
 \lesssim_j(h_n^j)^{-4(1-\theta)}\longrightarrow0.
\end{equation}
Set $z_n^j:=U_n^j-\widetilde U_n^j$. Its $L^4_{t,x}$ convergence
in \eqref{eq:large_profile_convergence} and uniform $H^1$ bound
from \eqref{eq:nonlinear_profile_bounds} imply
\begin{equation}\label{eq:large_profile_difference_L8}
 \|z_n^j\|_{L^8_{t,x}(\R\times\R^2)}^8
 \lesssim\|\nabla z_n^j\|_{L^\infty_t(\R;L^2_x(\R^2))}^4
          \|z_n^j\|_{L^4_{t,x}(\R\times\R^2)}^4\longrightarrow0.
\end{equation}
Hence \eqref{eq:large_profile_approximation_L8} and
\eqref{eq:large_profile_difference_L8} give
\begin{equation}\label{eq:large_profile_L8_smallness}
 \|U_n^j\|_{L^8_{t,x}(\R\times\R^2)}\longrightarrow0
 \qquad(h_n^j\to\infty).
\end{equation}
For a quintic mixed term containing this profile, put it in
$L^8_{t,x}$, one remaining factor in $L^4_{t,x}$, and the other
three in $L^8_{t,x}$. H\"older's inequality gives smallness in
$N^0(\R)$.
For $h_n^j\equiv1$, set
$z_n^j:=U_n^j-e^{i\gamma_n^j}v^j(t+\sigma_n^j,x-x_n^j)$.
By \eqref{eq:unit_profile_convergence},
\eqref{eq:nonlinear_profile_bounds}, and
\eqref{eq:appendix_GN_L8},
\[
 \|z_n^j\|_{L^8_{t,x}(\R\times\R^2)}^8
 \lesssim
 \|\nabla z_n^j\|_{L^\infty_t(\R;L^2_x(\R^2))}^4
 \|z_n^j\|_{L^4_{t,x}(\R\times\R^2)}^4\longrightarrow0.
\]
If $h_n^j=h_n^k\equiv1$, then $\xi_n^j=\xi_n^k=0$, and
\eqref{eq:profile_parameter_orthogonality} reduces to
$|t_n^j-t_n^k|+|x_n^j-x_n^k|\to\infty$.
Approximating $v^j$ and $v^k$ in $L^8_{t,x}$ by smooth compactly
supported functions, and using the preceding convergence, gives
\begin{equation}\label{eq:unit_profile_product_limit}
 \|U_n^jU_n^k\|_{L^4_{t,x}(\R\times\R^2)}\longrightarrow0
 \qquad(j\ne k,\ h_n^j=h_n^k\equiv1).
\end{equation}
The remaining three factors are controlled in $L^4_{t,x}$,
$L^8_{t,x}$, and $L^8_{t,x}$.
By \eqref{eq:linear_profile_remainder} and Sobolev embedding,
\begin{equation}\label{eq:free_remainder_spacetime_smallness}
 \lim_{J\to\infty}\limsup_{n\to\infty}
 \left(\|e^{it\Delta}r_n^J\|_{L^4_{t,x}(\R\times\R^2)}
       +\|e^{it\Delta}r_n^J\|_{L^8_{t,x}(\R\times\R^2)}\right)=0.
\end{equation}

We now prove \eqref{eq:profile_error_smallness}. For fixed $J_0$,
the mixed-term estimates \eqref{eq:profile_product_limit},
\eqref{eq:large_profile_L8_smallness}, and
\eqref{eq:unit_profile_product_limit} make the bracketed term in
\eqref{eq:profile_error_decomposition} tend to zero.
For the difference $F(S_n^{J_0})-F(V_n^J)$, use H\"older with the factor
$V_n^J-S_n^{J_0}$ in $L^4_{t,x}$ and the other factors in
$L^4_{t,x}$ for the cubic terms or $L^8_{t,x}$ for the quintic terms.
These other factors are uniformly bounded by
\eqref{eq:profile_sum_bound} and the Strichartz estimate \eqref{eq:appendix_strichartz}
applied to the free remainder. Thus \eqref{eq:nonlinear_profile_tail},
\eqref{eq:free_remainder_spacetime_smallness}, and
\eqref{eq:profile_nonlinearity_tail} control the two remaining terms.
Taking $\limsup_{n\to\infty}$ in
\eqref{eq:profile_error_decomposition}, then $J\to\infty$ and
$J_0\to\infty$ in this order, proves
\eqref{eq:profile_error_smallness}.
By the interpolation inequality \eqref{eq:appendix_bessel_interpolation} and H\"older's inequality,
\[
 \|e_n^J\|_{N^{1/2}(\R)}
 \lesssim\|e_n^J\|_{N^0(\R)}^{1/2}
             \|e_n^J\|_{N^1(\R)}^{1/2}.
\]
Combining \eqref{eq:profile_error_bound} and
\eqref{eq:profile_error_smallness}, we obtain
\[
 \lim_{J\to\infty}\limsup_{n\to\infty}
 \|\langle D\rangle^{1/2}e_n^J\|_{L^{4/3}_{t,x}(\R\times\R^2)}=0.
\]

Choose first $J$ and then $n$ sufficiently large. Since
$V_n^J(0)=f_n$, Lemma~\ref{lem:stability} gives a finite
$L^4_{t,x}$ norm for the solution with initial data $u(T_n)$.
Uniqueness and the local spacetime bound in
Lemma~\ref{lem:local_small_data} then give a finite forward
$L^4_{t,x}$ norm for $u$, a contradiction.

Thus one profile has mass $m_Q$, and \eqref{eq:profile_pythagorean}
excludes all other nonzero profiles. Relabeling it as the first gives
$\|r_n^1\|_{L^2_x(\R^2)}\to0$.

\par\medskip\Needspace{3\baselineskip}\noindent\textbf{Boundedness of the normalized time parameter.}\par\nobreak\smallskip
\noindent
Write $f_n=a_n+r_n^1$, omitting the index of the remaining
profile in \eqref{eq:nonradialprofiles}. Set $\sigma_n:=-t_n/h_n^2$.
The free evolution of $a_n$ is
\[
 e^{is\Delta}a_n(x)
 =e^{i\gamma_n+ix\cdot\xi_n-is|\xi_n|^2}h_n^{-1}
 [e^{i(s/h_n^2+\sigma_n)\Delta}P_n\phi]
       ((x-x_n-2s\xi_n)/h_n).
\]
If $\sigma_n\to+\infty$, the Strichartz estimate \eqref{eq:appendix_strichartz} gives
\begin{equation}\label{eq:profile_free_future_tail}
 \|e^{ir\Delta}P_n\phi\|_{L^4_{r,y}(\{r>\sigma_n\}\times\R^2)}
 \leq\|e^{ir\Delta}\phi\|_{L^4_{r,y}(\{r>\sigma_n\}\times\R^2)}
       +C\|(P_n-I)\phi\|_{L^2_y(\R^2)}\longrightarrow0.
\end{equation}
For $\sigma_n\to-\infty$, the analogue of
\eqref{eq:profile_free_future_tail} holds on $\{r<\sigma_n\}$. Scaling and the Strichartz estimate for
$r_n^1$ give $\|e^{is\Delta}f_n\|_{L^4_{s,x}(I\times\R^2)}\to0$,
where $I:=[0,\infty)$ or $I:=(-\infty,0]$, respectively.
Interpolation and the Strichartz estimate \eqref{eq:appendix_strichartz} with $s=1$ yield
\[
 \|\langle D\rangle^{1/2}e^{is\Delta}f_n\|_{L^4_{s,x}(I\times\R^2)}
 \leq C\|e^{is\Delta}f_n\|_{L^4_{s,x}(I\times\R^2)}^{1/2}
            \|f_n\|_{H^1_x(\R^2)}^{1/2}\longrightarrow0.
\]
Applying Lemma~\ref{lem:stability} with $\widetilde v\equiv0$,
$e\equiv0$, $t_0=0$, and initial data $f_n$,
we obtain $\|u\|_{L^4_{t,x}([T_n,\infty)\times\R^2)}<\infty$
if $\sigma_n\to+\infty$, and
$\sup_n\|u\|_{L^4_{t,x}([0,T_n]\times\R^2)}<\infty$
if $\sigma_n\to-\infty$.
Both conclusions contradict
$\|u\|_{L^4_{t,x}([0,T]\times\R^2)}\to\infty$ as $T\to\infty$.
Hence $\sigma_n\to\sigma\in\R$. Replacing $\phi$ by
$e^{i\sigma\Delta}\phi$ and using $P_n\phi\to\phi$ in $L^2$,
we obtain \eqref{eq:singlefullprofile}.
When $h_n\equiv1$, we have $\phi\in H^1$ and $\xi_n=0$.
\end{proof}

The Palais--Smale condition yields uniform spatial localization of
mass. When $h_n\to\infty$, Proposition~\ref{prop:largescale} and the
threshold classification for \eqref{eq:cubic_nls} identify the
limiting profile as a soliton.

\begin{proposition}
\label{prop:compactness}
Let $u$ be an $H^1$ solution of \eqref{eq:cq} satisfying
\[
 M(u)=m_Q,\qquad
 \norm u_{L^4_{t,x}([0,\infty)\times\R^2)}=\infty.
\]
For $t\geq0$, define the spatial scale function $\lambda(t)$ by
\begin{equation}\label{eq:halfwidth}
 \begin{gathered}
 \mathcal C_t(r):=\sup_{y\in\R^2}\int_{B(y,r)}|u(t,x)|^2\dd x,\\
 R(t):=\inf\{r>0:\mathcal C_t(r)\ge m_Q/2\},\qquad
 \lambda(t):=\max\{1,R(t)\}.
 \end{gathered}
\end{equation}
Then there exists a spatial center $c:[0,\infty)\to\R^2$ such that
\[
 \int_{B(c(t),\lambda(t))}|u(t,x)|^2\dd x\geq\frac{m_Q}{2}
 \qquad(t\geq0).
\]
For every $\epsilon>0$, there exists $A_\epsilon<\infty$ such that
\[
 \sup_{t\ge0}\int_{|x-c(t)|>A_\epsilon\lambda(t)}
 |u(t,x)|^2\dd x<\epsilon.
\]
Moreover,
\begin{equation}\label{eq:subparabolic}
 \lambda(t)=o(\sqrt t)\quad(t\to\infty).
\end{equation}
Define the soliton orbit by
\[
 \mathcal O_Q:=
 \{e^{i\gamma+ix\cdot\nu}\ell^{-1}Q((x-y)/\ell):
       \gamma\in\R,\ \nu,y\in\R^2,\ \ell>0\}.
\]
Then
\begin{equation}\label{eq:Qlarge}
 \lim_{L\to\infty}\
 \sup_{\{t\ge0:\lambda(t)\ge L\}}
 \operatorname{dist}_{L^2_x(\R^2)}(u(t),\mathcal O_Q)=0,
\end{equation}
where the supremum over the empty set is taken to be zero.
\end{proposition}

\begin{proof}
\proofstep{1}{Uniform spatial localization of mass}
First choose a measurable center. The ball mass
$(t,y,r)\mapsto\int_{B(y,r)}|u(t,x)|^2\,dx$ is continuous
and tends to zero as $|y|\to\infty$ for fixed $t$ and $r$.
Thus, for some $y_{t,r}\in\R^2$,
\[
 \mathcal C_t(r)=\int_{B(y_{t,r},r)}|u(t,x)|^2\,dx.
\]
The function $\mathcal C_t$ is continuous, tends to zero as
$r\downarrow0$, and tends to $m_Q$ as $r\to\infty$. Hence
$0<R(t)<\infty$ and $\mathcal C_t(\lambda(t))\geq m_Q/2$.
Both $R$ and $\lambda$ are Borel measurable. The set
\[
 K(t):=\left\{y\in\R^2:
       \int_{B(y,\lambda(t))}|u(t,x)|^2\,dx\geq m_Q/2\right\}
\]
is nonempty and compact. To select a point measurably, enumerate
$\mathbb Z^2$ as $k_1,k_2,\ldots$, put $D_j:=k_j+[0,1]^2$, and set
\[
 j_0(t):=\min\{j\geq1:D_j\cap K(t)\ne\varnothing\},
 \qquad C_0(t):=D_{j_0(t)}.
\]
Divide each selected square $C_m(t)$ into four closed squares of
half the side length, ordered from lower left to lower right, then
upper left to upper right. Let $C_{m+1}(t)$ be the first meeting
$K(t)$; such a square exists since $C_m(t)\cap K(t)\ne\varnothing$.
The squares are nested, with
\[
 \operatorname{diam}C_m(t)=\sqrt2\,2^{-m}\longrightarrow0,
 \qquad \bigcap_{m\geq0}C_m(t)=\{c(t)\}.
\]
Since $K(t)$ is closed and meets every $C_m(t)$, we have $c(t)\in K(t)$.

For any fixed square $D$ in this construction, continuity in the
center gives
\[
 \max_{y\in D}\int_{B(y,\lambda(t))}|u(t,x)|^2\,dx
 =\sup_{y\in D\cap\mathbb Q^2}
       \int_{B(y,\lambda(t))}|u(t,x)|^2\,dx.
\]
Since $\lambda$ is Borel, the right-hand side and the condition
$D\cap K(t)\ne\varnothing$ are Borel in $t$. Only countably many
fixed squares occur at each stage, so induction makes each
first-square choice Borel. Thus $c(t)$ is Borel as the limit of
the centers of $C_m(t)$.

Given $T_n\to\infty$, Lemma~\ref{lem:single_profile_compactness}
gives a subsequence, still denoted by $(T_n)$, satisfying
\eqref{eq:singlefullprofile}. Set
\[
 \Phi_n(x):=e^{i\gamma_n+ix\cdot\xi_n}h_n^{-1}
                  \phi((x-x_n)/h_n).
\]
Uniformly over measurable sets $\Omega\subset\R^2$,
\[
 \left|\int_\Omega|u(T_n,x)|^2\,dx
       -\int_\Omega|\Phi_n(x)|^2\,dx\right|
 \leq2\sqrt{m_Q}\,\|u(T_n)-\Phi_n\|_{L^2_x(\R^2)}=o(1).
\]
Since $|\phi|^2\in L^1(\R^2)$ and $M(\phi)=m_Q$, choose
$0<r_0<1<R_0$ such that
\[
 \sup_{y\in\R^2}\int_{B(y,r_0)}|\phi(x)|^2\,dx<\frac{m_Q}{4},
 \qquad
 \int_{B(0,R_0)}|\phi(x)|^2\,dx>\frac{3m_Q}{4}.
\]
For large $n$, every ball of radius $r_0h_n$ contains less than
$m_Q/2$ of the mass of $u(T_n)$, whereas $B(x_n,R_0h_n)$ contains
more than $m_Q/2$. Since $h_n\geq1$,
\[
 r_0h_n\leq R(T_n)\leq\lambda(T_n)\leq R_0h_n.
\]
The balls $B(x_n,R_0h_n)$ and $B(c(T_n),\lambda(T_n))$ intersect
because their masses sum to more than $m_Q$. Thus
$|c(T_n)-x_n|\leq\lambda(T_n)+R_0h_n\leq2R_0h_n$.
With $a=r_0$ and $A=2R_0$, we obtain
\begin{equation}\label{eq:widthcomparison}
 a\le\lambda(T_n)/h_n\le A,\qquad |c(T_n)-x_n|/h_n\le A.
\end{equation}
The tails of $\phi$ and \eqref{eq:widthcomparison} give uniform
spatial localization of mass on this subsequence. On bounded time
intervals, $L^2$ continuity gives the same localization and bounds
for $\lambda$ and $c$. Failure of uniform localization for $t\geq0$
would therefore yield $T_n\to\infty$ contradicting this subsequence
argument.

\proofstep{2}{Identification of the limiting profile}
If $\lambda(T_n)\to\infty$, local boundedness gives $T_n\to\infty$.
Pass to a subsequence satisfying \eqref{eq:singlefullprofile} and
\eqref{eq:widthcomparison}; the latter gives $h_n\to\infty$.
Let $w$ solve \eqref{eq:cubic_nls} with $w(0)=\phi$ on its
maximal lifespan. Forward scattering of $w$ would contradict the
hypothesis on $u$ by Proposition~\ref{prop:largescale} on $[0,\infty)$.
Similarly, $w$ cannot scatter backward. Its $L^4_{s,y}$ norm therefore
diverges in both time directions of its maximal lifespan.

Since $M(w)=m_Q$, Lemma~\ref{lem:cubic_threshold_classification}
applies. A pseudoconformal solution satisfies
$\norm{w(s)}_{L^4_y(\R^2)}^4=C(T-s)^{-2}$, so its $L^4_{s,y}$ norm
is finite toward $-\infty$ on the branch $s<T$, or toward
$+\infty$ on the time-reversed branch $s>T$.
Thus $w$ is a soliton and $\phi\in\mathcal O_Q$.
Absorbing its fixed symmetry parameters into
\eqref{eq:singlefullprofile} gives \eqref{eq:Qlarge} by contradiction
along sequences with $\lambda(T_n)\to\infty$. If $\ell$ is the scale
of $\phi$, the resulting $Q$ profile has scale
$\ell h_n\sim\lambda(T_n)$ by \eqref{eq:widthcomparison}.

\proofstep{3}{The growth of the spatial scale}
If \eqref{eq:subparabolic} fails, there exist $\epsilon>0$ and
$T_n\to\infty$ with $\lambda(T_n)\geq\epsilon\sqrt{T_n}$.
Step~2 gives a soliton $w$, and \eqref{eq:widthcomparison} gives
$T_n/h_n^2\leq A$ for some $A<\infty$.
Since $w$ has finite $L^4_{s,y}$ norm on $[-A-1,0]$, applying
Proposition~\ref{prop:largescale} on this interval and using
\[
 [-T_n,0]\subset h_n^2[-A-1,0],
\]
we obtain, by time-translation invariance of \eqref{eq:cq},
\[
 \sup_n\|u\|_{L^4_{t,x}([0,T_n]\times\R^2)}<\infty.
\]
This contradicts the hypothesis on $u$ and proves
\eqref{eq:subparabolic}.
\end{proof}

\begin{remark}\label{rem:spatial_width}
In \eqref{eq:halfwidth}, $R(t)$ is the smallest radius of a ball
containing at least half the mass, with its center allowed to vary.
Thus $\lambda(t)=\max\{1,R(t)\}$ measures the spatial scale of mass
concentration. Proposition~\ref{prop:compactness} shows that
$B(c(t),A_\epsilon\lambda(t))$ contains all but $\epsilon$ of the
mass, uniformly in $t\geq0$.
\end{remark}

\section{Virial estimates and time-averaged kinetic energy evacuation}\label{sec:rigidity}
In this section, we study the minimal-mass nonscattering solution
defined by \eqref{eq:minimal_mass_nonscattering}.
Using Proposition~\ref{prop:compactness} and localized virial
identities, we prove sublinear motion of its spatial center and
time-averaged decay of $\norm{u(t)}_{L^4_x(\R^2)}^4$ and $\norm{u(t)}_{L^6_x(\R^2)}^6$.
The total kinetic energy has time average $2E>0$, while its
time average inside any ball of radius $vT$ tends to zero for
every fixed $v>0$. We also describe the soliton approximation
and the kinetic energy retained by its remainder, before establishing
the reduced Duhamel formula used in the exterior estimates.

\subsection{Localized virial estimates}
We first control the motion of the spatial center, then use a
localized virial identity to obtain the time-averaged nonlinear
and kinetic estimates.

\begin{theorem}\label{thm:main}
Let $u$ be a minimal-mass nonscattering solution as in
\eqref{eq:minimal_mass_nonscattering}.
Let $\lambda(t)$ and $c(t)$ be the spatial scale function
and spatial center given by Proposition~\ref{prop:compactness}. Then $\lambda(t)=o(\sqrt t)$,
and there exists $C_0>0$ such that
\begin{equation}\label{eq:width}
 1\leq\lambda(t)\leq C_0\sqrt{1+t}\qquad(t\geq0).
\end{equation}
For every $\epsilon>0$, there exists $A_\epsilon\geq1$ such that
\begin{equation}\label{eq:MC}
 \sup_{t\geq0}\int_{|x-c(t)|>A_\epsilon\lambda(t)}|u(t,x)|^2\dd x
 \leq\epsilon.
\end{equation}
Moreover,
\begin{equation}\label{eq:center_sublinear}
 R_c(T):=\sup_{0\leq t\leq T}|c(t)|=o(T)\qquad(T\to\infty),
\end{equation}
and
\begin{align}
 \lim_{T\to\infty}\avgT{\norm{u(t)}_{L^6_x(\R^2)}^6}
 &=\lim_{T\to\infty}\avgT{\norm{u(t)}_{L^4_x(\R^2)}^4}=0,\label{eq:nlevac}\\
 \lim_{T\to\infty}\avgT{\norm{\nabla u(t)}_{L^2_x(\R^2)}^2}&=2E>0.\label{eq:totalavg}
\end{align}
\end{theorem}

\begin{proof}[Proof of Theorem~\ref{thm:main}]
By \eqref{eq:coercive_energy} and conservation of mass and energy,
$\sup_{t\geq0}\norm{u(t)}_{H^1}\leq B_H<\infty$, where $B_H$
is defined in \eqref{eq:uniform_H1_bound}.
Proposition~\ref{prop:compactness} gives \eqref{eq:width} and
\eqref{eq:MC}.

\proofstep{1}{Nonlinear tail estimates}
Fix $\epsilon>0$ and choose $A_\epsilon$ as in \eqref{eq:MC}.
Let $\zeta\in C^\infty(\R^2;[0,1])$ equal to 0 on $B(0,1)$ and 1 outside $B(0,2)$. Set
\[
 v(t,x):=\zeta\!\left(\frac{x-c(t)}{A_\epsilon\lambda(t)}\right)u(t,x).
\]
Using \eqref{eq:MC}, the bound $A_\epsilon\lambda(t)\geq1$,
and the Gagliardo--Nirenberg inequality, we obtain, for
$2<p<\infty$,
\[
 \norm{v(t)}_{L^2_x(\R^2)}^2\leq\epsilon,\qquad
 \norm{\nabla v(t)}_{L^2_x(\R^2)}\leq B_H+C\sqrt{m_Q},\qquad
 \norm{v(t)}_{L^p_x(\R^2)}^p\leq C_{p,B_H,m_Q}\epsilon.
\]
As $v(t,x)=u(t,x)$ for $|x-c(t)|\geq2A_\epsilon\lambda(t)$,
it follows that
\begin{equation}\label{eq:nonlinear_spatial_tails}
 \kappa_p(A):=\sup_{t\geq0}
 \int_{|x-c(t)|>A\lambda(t)}|u(t,x)|^p\dd x
 \longrightarrow0\qquad(A\to\infty).
\end{equation}

\proofstep{2}{Sublinear motion of the spatial center}
The local boundedness of $c(t)$ established in the
proof of Proposition~\ref{prop:compactness} gives $R_c(T)<\infty$
for every $T>0$. Choose $\theta\in C_c^\infty(\R^2;[0,1])$ equal to 1 on $B(0,1)$
and 0 outside $B(0,2)$, and define
\[
 h_R(x):=x\theta(x/R),\qquad
 Z_R(t):=\int h_R(x)|u(t,x)|^2\dd x.
\]
Since $P(u)=0$, the continuity equation gives
\[
 Z_R'(t)=2\int(\nabla h_R-I)\im(\overline u\nabla u)\dd x,\qquad
 |Z_R'(t)|\leq CB_H
 \left(\int_{|x|>R}|u(t,x)|^2\dd x\right)^{1/2},
\]
where $\nabla h_R=(\partial_j(h_R)_i)_{i,j=1}^2$ is the Jacobian
matrix of $h_R$ and $I$ is the $2\times2$ identity matrix.
For fixed $\epsilon>0$ and $T>0$, choose
\[
 R=2R_c(T)+2A_\epsilon C_0\sqrt{1+T}+2,
\]
then $B(c(t),A_\epsilon\lambda(t))\subset B(0,R)$ for every
$t\in[0,T]$. Integrating over \([0,T]\) and combining with \eqref{eq:MC}, we obtain
\begin{equation}\label{eq:truncated_center_motion}
 |Z_R(t)-Z_R(0)|\leq CB_H\sqrt{\epsilon}\,T
 \qquad(0\leq t\leq T).
\end{equation}
On the other hand, on $B(c(t),A_\epsilon\lambda(t))$ we have
$h_R(x)=x$. On its complement, $|h_R(x)-c(t)|\leq CR$.
Since $M(u)=m_Q$, splitting the integral below over these two
sets gives
\begin{equation}\label{eq:truncated_center_error}
 \begin{aligned}
 |Z_R(t)-m_Qc(t)|
 &=\left|\int(h_R(x)-c(t))|u(t,x)|^2\dd x\right|\\
 &\leq m_QA_\epsilon\lambda(t)+C\epsilon R
 \qquad(0\leq t\leq T).
 \end{aligned}
\end{equation}
Combining \eqref{eq:truncated_center_motion} and
\eqref{eq:truncated_center_error}, and taking the supremum over
$0\leq t\leq T$, we obtain
\[
 m_QR_c(T)\leq m_Q|c(0)|+CB_H\sqrt{\epsilon}\,T
      +C_\epsilon\sqrt{1+T}+C\epsilon R_c(T),
\]
where $C$ is independent of $\epsilon$ and $T$.
For sufficiently small $\epsilon$, absorb the last term and
divide by $T$. Letting first $T\to\infty$ and then
$\epsilon\downarrow0$ proves \eqref{eq:center_sublinear}.

\proofstep{3}{Localized virial identity}
Choose a smooth nonincreasing function $\rho:[0,\infty)\to[0,1]$
equal to 1 on $[0,1]$ and 0 on $[2,\infty)$. Define
\[
 \varphi:=\rho^2,\qquad
 \psi(r):=\frac1r\int_0^r\varphi(s)\dd s\quad(r>0),\qquad \psi(0)=1,
\]
and set
\[
 X_R(x):=\psi(|x|/R)x,\qquad
 \mathcal A_R(t):=\im\int\overline{u(t,x)}X_R(x)\cdot\nabla u(t,x)\dd x.
\]
For $a=\rho,\varphi,\psi$, write $a_R(x):=a(|x|/R)$.
The radial and tangential eigenvalues of $\nabla X_R$ are $\varphi_R$
and $\psi_R$, respectively. Since $\varphi$ is nonincreasing,
$0\leq\varphi_R\leq\psi_R\leq1$. In particular,
\begin{equation}\label{eq:virial_multiplier_bounds}
 \begin{gathered}
 \norm{\nabla X_R}_{L^\infty_x(\R^2)}\leq1,\qquad \norm{X_R}_{L^\infty_x(\R^2)}\leq2R,\\
 \operatorname{div}X_R=\varphi_R+\psi_R\leq2,\qquad
 \norm{\Delta\operatorname{div}X_R}_{L^\infty_x(\R^2)}\leq CR^{-2},\\
 |\mathcal A_R(t)|\leq CRB_H\sqrt{m_Q}.
 \end{gathered}
\end{equation}
Moreover, $X_R(x)=x$ and $\operatorname{div}X_R=2$ on $B(0,R)$.
For $x\ne0$, the radial derivative and tangential gradient are
\[
 \partial_ru:=\frac{x}{|x|}\cdot\nabla u,\qquad
 \not\nabla u:=\nabla u-\frac{x}{|x|}\partial_ru.
\]
The localized virial identity for \eqref{eq:cq} is
\begin{align}
 \mathcal A_R'={}&
 2\int\varphi_R|\partial_r u|^2+
 2\int\psi_R|\not\nabla u|^2
 -\frac12\int\Delta(\varphi_R+\psi_R)|u|^2\notag\\
 &-\frac12\int(\varphi_R+\psi_R)|u|^4
 +\frac23\int(\varphi_R+\psi_R)|u|^6.\label{eq:virialexact}
\end{align}
For smooth solutions, this follows by integration by parts
\cite{Cazenave,KM,Murphy}. For each $T_0>0$,
Lemma~\ref{lem:local_small_data} gives smooth solutions $u_m$ of
\eqref{eq:cq} with $u_m\to u$ in $C([0,T_0];H^1(\R^2))$.
Passing to the limit in the integrated identity gives
\eqref{eq:virialexact} for $u$.

Using $\psi_R\geq\varphi_R$, $\varphi_R+\psi_R\leq2$, and
$\varphi_R+\psi_R=2$ on $B(0,R)$, we deduce from
\eqref{eq:virialexact} that
\begin{align}
 \mathcal A_R'\geq{}&
 2\int\varphi_R|\nabla u|^2-\norm u_{L^4_x(\R^2)}^4+\frac43\norm u_{L^6_x(\R^2)}^6
 -Cm_QR^{-2}-C\int_{|x|>R}|u|^6.\label{eq:viriallower}
\end{align}
\proofstep{4}{Time-averaged limits}
It remains to prove \eqref{eq:nlevac} and \eqref{eq:totalavg}.
For $A\geq1$, set
\begin{equation}\label{eq:virial_radius}
 R_A(T):=1+R_c(T)+AC_0\sqrt{1+T}.
\end{equation}
By \eqref{eq:center_sublinear}, for fixed $A$ we have
$R_A(T)/T\to0$ and $R_A(T)\to\infty$ as $T\to\infty$.
For $0\leq t\leq T$, \eqref{eq:width} implies
\[
 \{|x|>R_A(T)\}\subset\{|x-c(t)|>A\lambda(t)\}.
\]
Since $M(\rho_Ru)\leq m_Q$, \eqref{eq:GN}
and integration by parts give
\[
 \norm{\rho_Ru}_{L^4_x(\R^2)}^4\leq2\norm{\nabla(\rho_Ru)}_{L^2_x(\R^2)}^2
 =2\int\varphi_R|\nabla u|^2-2\int\rho_R\Delta\rho_R|u|^2.
\]
Since $\rho_R=1$ on $B(0,R)$ and
$\norm{\rho_R\Delta\rho_R}_{L^\infty_x(\R^2)}\leq CR^{-2}$, it follows that
\begin{equation}\label{eq:localized_threshold_GN}
 2\int\varphi_R|\nabla u|^2-\norm u_{L^4_x(\R^2)}^4
 \geq-Cm_QR^{-2}-\int_{|x|>R}|u|^4.
\end{equation}
Inserting \eqref{eq:localized_threshold_GN} into
\eqref{eq:viriallower} and integrating with the fixed radius
$R=R_A(T)$ from \eqref{eq:virial_radius}, we obtain
\begin{equation}\label{eq:averaged_sextic_bound}
 \frac43\avgT{\norm u_{L^6_x(\R^2)}^6}
 \leq C\frac{R_A(T)}T+CR_A(T)^{-2}
      +C\kappa_4(A)+C\kappa_6(A).
\end{equation}
Letting $T\to\infty$ in \eqref{eq:averaged_sextic_bound} and
then $A\to\infty$, we obtain the sixth-power limit in
\eqref{eq:nlevac} by \eqref{eq:nonlinear_spatial_tails}.
By \eqref{eq:coercive_energy} and Cauchy--Schwarz in time,
\[
 \avgT{\norm{u(t)}_{L^4_x(\R^2)}^4}
 \leq\sqrt{m_Q}\left(\avgT{\norm{u(t)}_{L^6_x(\R^2)}^6}\right)^{1/2}
 \longrightarrow0.
\]
Energy conservation now gives
\[
 \avgT{\norm{\nabla u(t)}_{L^2_x(\R^2)}^2}
 =2E+\frac12\avgT{\norm{u(t)}_{L^4_x(\R^2)}^4}
     -\frac13\avgT{\norm{u(t)}_{L^6_x(\R^2)}^6}
 \longrightarrow2E,
\]
which proves \eqref{eq:totalavg}. The positivity of $E$ follows
from \eqref{eq:positive}.
\end{proof}

\subsection{Time-averaged kinetic energy evacuation}
We apply the localized virial identity to translated vector fields
to show that the time-averaged kinetic energy in $B(z,vT)$ tends
to zero uniformly in $z$.

\begin{corollary}\label{cor:ballistic_evacuation}
Let $u$ be a minimal-mass nonscattering solution.
Then, for every fixed $v>0$,
\begin{equation}\label{eq:ballistic_evacuation}
 \lim_{T\to\infty}\sup_{z\in\R^2}
 \frac1T\int_0^T\int_{|x-z|\leq vT}|\nabla u(t,x)|^2\dd x\dd t=0.
\end{equation}
The center $z$ may depend on $T$, but remains fixed during the time integration.
Moreover, as $T\to\infty$,
\begin{align}
 \sup_{0\leq t\leq T}\int_{|x|>vT}|u(t,x)|^2\dd x
 &\longrightarrow0,\label{eq:ballistic_mass}\\
 \frac1T\int_0^T\int_{|x|>vT}|\nabla u(t,x)|^2\dd x\dd t
 &\longrightarrow2E.\label{eq:ballistic_energy}
\end{align}
\end{corollary}
\begin{proof}
Fix $v,\epsilon>0$ and choose $A_\epsilon$ as in \eqref{eq:MC}.
For $R>0$ and $z\in\R^2$, define
\[
 \mathcal A_{R,z}(t):=
 \im\int\overline{u(t,x)}X_R(x-z)\cdot\nabla u(t,x)\dd x,
\]
where $X_R$ is the vector field in \eqref{eq:virialexact}.
The argument leading to \eqref{eq:virialexact}, with $X_R(x)$
replaced by $X_R(x-z)$, gives the translated virial identity.
Dropping the nonnegative sixth-power term yields
\begin{equation}\label{eq:translated_virial_lower}
 \mathcal A_{R,z}'(t)
 \geq2\int_{|x-z|\leq R}|\nabla u|^2-\norm u_{L^4_x(\R^2)}^4-Cm_QR^{-2}.
\end{equation}
Since $P(u)=0$, at each fixed time $t$ we have
\[
 \mathcal A_{R,z}(t)=\im\int\overline{u(t,x)}
 \bigl[X_R(x-z)-X_R(c(t)-z)\bigr]\cdot\nabla u(t,x)\dd x.
\]
By \eqref{eq:virial_multiplier_bounds},
\[
 |X_R(x-z)-X_R(c(t)-z)|\leq
 \begin{cases}
  A_\epsilon\lambda(t),
       &|x-c(t)|\leq A_\epsilon\lambda(t),\\
  4R,&|x-c(t)|>A_\epsilon\lambda(t).
 \end{cases}
\]
Splitting the integral over these two regions, Cauchy--Schwarz,
mass conservation, and \eqref{eq:MC} give
\begin{equation}\label{eq:translated_action_bound}
 \begin{aligned}
 |\mathcal A_{R,z}(t)|
 &\leq\norm{\nabla u(t)}_{L^2_x(\R^2)}\left[
 A_\epsilon\lambda(t)\norm{u(t)}_{L^2_x(\R^2)}
 +4R\left(\int_{|x-c(t)|>A_\epsilon\lambda(t)}|u(t,x)|^2\dd x\right)^{1/2}
 \right]\\
 &\leq B_HA_\epsilon\lambda(t)\sqrt{m_Q}+4B_HR\sqrt\epsilon,
 \end{aligned}
\end{equation}
uniformly in $z$. Integrating the virial inequality
\eqref{eq:translated_virial_lower} over $[0,T]$ with $R=vT$, we bound the endpoint terms $\mathcal A_{R,z}(T)$
and $\mathcal A_{R,z}(0)$ by \eqref{eq:translated_action_bound}.
Dividing by $2T$ and using \eqref{eq:width} and \eqref{eq:nlevac},
we obtain
\[
 \limsup_{T\to\infty}\sup_{z\in\R^2}
 \frac1T\int_0^T\int_{|x-z|\leq vT}|\nabla u(t,x)|^2\dd x\dd t
 \leq4B_Hv\sqrt\epsilon.
\]
Letting $\epsilon\downarrow0$ proves \eqref{eq:ballistic_evacuation}.

For fixed $\epsilon>0$, \eqref{eq:width} and \eqref{eq:center_sublinear}
imply $B(c(t),A_\epsilon\lambda(t))\subset B(0,vT)$ for all
$0\leq t\leq T$ once $T$ is sufficiently large.
Consequently, \eqref{eq:MC} gives
\[
 \sup_{0\leq t\leq T}\int_{|x|>vT}|u(t,x)|^2\dd x
 \leq\sup_{0\leq t\leq T}
 \int_{|x-c(t)|>A_\epsilon\lambda(t)}|u(t,x)|^2\dd x
 \leq\epsilon.
\]
Letting $T\to\infty$ and then $\epsilon\downarrow0$ proves
\eqref{eq:ballistic_mass}.
Finally, subtracting \eqref{eq:ballistic_evacuation} with $z=0$
from \eqref{eq:totalavg} gives \eqref{eq:ballistic_energy}.
\end{proof}

\subsection{Soliton approximation and the kinetic energy remainder}
We combine \eqref{eq:nlevac} with Proposition~\ref{prop:compactness}
to obtain approximation by the soliton orbit in time density and
identify the kinetic energy retained by the remainder.

\begin{corollary}\label{cor:Qdensity}
Let $u$ be a minimal-mass nonscattering solution.
Then, as $T\to\infty$, for every fixed $L<\infty$,
\begin{equation}\label{eq:widthdensity}
 \frac{|\{t\in[0,T]:\lambda(t)\leq L\}|}{T}\longrightarrow0.
\end{equation}
For every $\eta>0$,
\begin{equation}\label{eq:Qdensity}
 \frac{|\{t\in[0,T]:\operatorname{dist}_{L^2}(u(t),\mathcal O_Q)>\eta\}|}{T}
 \longrightarrow0,
\end{equation}
where $\mathcal O_Q$ is the soliton orbit defined in
Proposition~\ref{prop:compactness}.
If $\lambda(t_n)\to\infty$, there exist decompositions
\begin{equation}\label{eq:Qdefectform}
 \begin{gathered}
 u(t_n)=q_n+r_n,\qquad
 q_n(x):=e^{i\gamma_n+ix\cdot\nu_n}\ell_n^{-1}
 Q\left(\frac{x-y_n}{\ell_n}\right),\\
 \ell_n\sim\lambda(t_n),\qquad
 \norm{r_n}_{L^2_x(\R^2)}\to0,
 \end{gathered}
\end{equation}
The implicit constants are independent of $n$. Moreover,
\begin{equation}\label{eq:Qdefect}
 \nu_n\to0,\qquad\norm{\nabla q_n}_{L^2_x(\R^2)}\to0,\qquad
 \norm{\nabla r_n}_{L^2_x(\R^2)}^2\to2E>0.
\end{equation}
In particular, for every fixed $N>0$,
\[
 \norm{\nabla P_{\leq N}u(t_n)}_{L^2_x(\R^2)}\to0,\qquad
 \norm{\nabla P_{>N}u(t_n)}_{L^2_x(\R^2)}^2\to2E.
\]
\end{corollary}
\begin{proof}
By \eqref{eq:halfwidth} and the choice of $c(t)$ in
Proposition~\ref{prop:compactness}, the ball $B(c(t),\lambda(t))$
contains at least $m_Q/2$ of the mass. Cauchy--Schwarz on this ball gives
\[
 \begin{aligned}
 \frac{m_Q}{2}
 &\leq\int_{B(c(t),\lambda(t))}|u(t,x)|^2\dd x\\
 &\leq\sqrt\pi\lambda(t)
       \left(\int_{B(c(t),\lambda(t))}|u(t,x)|^4\dd x\right)^{1/2}
 \leq\sqrt\pi\lambda(t)\norm{u(t)}_{L^4_x(\R^2)}^2.
 \end{aligned}
\]
Hence
\begin{equation}\label{eq:halfwidth_L4_lower}
 \norm{u(t)}_{L^4_x(\R^2)}^4\geq m_Q^2/(4\pi\lambda(t)^2).
\end{equation}
By \eqref{eq:halfwidth_L4_lower} and \eqref{eq:nlevac}, we obtain
\eqref{eq:widthdensity}.
Combining \eqref{eq:widthdensity} with \eqref{eq:Qlarge}, we obtain
\eqref{eq:Qdensity}.
For a sequence with $\lambda(t_n)\to\infty$, \eqref{eq:Qlarge}
provides $q_n\in\mathcal O_Q$ with $\delta_n:=\norm{u(t_n)-q_n}_{L^2_x(\R^2)}\to0$.
The half-mass radius argument for \eqref{eq:widthcomparison}
gives $\ell_n\sim\lambda(t_n)$ for all sufficiently large $n$.
Enlarging the comparison constant to include the remaining
finitely many terms proves \eqref{eq:Qdefectform}.

Since $u(t_n)=q_n+r_n$, the momentum difference satisfies
\[
 \begin{aligned}
 P(u(t_n))-P(q_n)
 &=\im\int\bigl(\overline{r_n}\nabla u(t_n)
                    +\overline{q_n}\nabla r_n\bigr)\dd x\\
 &=\im\int\overline{r_n}\nabla\bigl(u(t_n)+q_n\bigr)\dd x.
 \end{aligned}
\]
Since $P(u(t_n))=0$ and $P(q_n)=m_Q\nu_n$, Cauchy--Schwarz gives
\[
 \begin{aligned}
 m_Q|\nu_n|
 &\leq\delta_n\bigl(\norm{\nabla u(t_n)}_{L^2_x(\R^2)}+\norm{\nabla q_n}_{L^2_x(\R^2)}\bigr)\\
 &\leq\delta_n\bigl(B_H+\sqrt{m_Q}|\nu_n|
                         +\ell_n^{-1}\norm{\nabla Q}_{L^2_x(\R^2)}\bigr).
 \end{aligned}
\]
For sufficiently large $n$, we absorb the $|\nu_n|$ term to obtain
$\nu_n\to0$. Since $\ell_n\to\infty$,
\begin{equation}\label{eq:soliton_gradient_decay}
 \norm{\nabla q_n}_{L^2_x(\R^2)}
 \leq\sqrt{m_Q}|\nu_n|+\ell_n^{-1}\norm{\nabla Q}_{L^2_x(\R^2)}
 \longrightarrow0.
\end{equation}
It follows from \eqref{eq:uniform_H1_bound} that
$\sup_n\norm{r_n}_{H^1}<\infty$.
For every $2<p<\infty$, scaling of $q_n$ and the
Gagliardo--Nirenberg inequality applied to $r_n$ give
\begin{equation}\label{eq:soliton_remainder_Lp_decay}
 \begin{aligned}
 \norm{q_n}_{L^p_x(\R^2)}&=\ell_n^{2/p-1}\norm Q_{L^p_x(\R^2)}\longrightarrow0,\\
 \norm{r_n}_{L^p_x(\R^2)}&\leq C_p\norm{r_n}_{L^2_x(\R^2)}^{2/p}
                      \norm{\nabla r_n}_{L^2_x(\R^2)}^{1-2/p}\longrightarrow0.
 \end{aligned}
\end{equation}
Thus $\norm{u(t_n)}_{L^p_x(\R^2)}\to0$. Taking $p=4,6$ in
\eqref{eq:soliton_remainder_Lp_decay} and using energy conservation, we obtain
\[
 \norm{\nabla u(t_n)}_{L^2_x(\R^2)}^2
 =2E+\frac12\norm{u(t_n)}_{L^4_x(\R^2)}^4-\frac13\norm{u(t_n)}_{L^6_x(\R^2)}^6
 \longrightarrow2E.
\]
Moreover, \eqref{eq:soliton_gradient_decay} implies
\[
 \left|\norm{\nabla r_n}_{L^2_x(\R^2)}^2-\norm{\nabla u(t_n)}_{L^2_x(\R^2)}^2\right|
 \leq\norm{\nabla q_n}_{L^2_x(\R^2)}
       \bigl(\norm{\nabla r_n}_{L^2_x(\R^2)}+\norm{\nabla u(t_n)}_{L^2_x(\R^2)}\bigr)
 \longrightarrow0.
\]
This proves \eqref{eq:Qdefect}. For fixed $N>0$,
\[
 \norm{\nabla P_{\leq N}u(t_n)}_{L^2_x(\R^2)}
 \leq\norm{\nabla q_n}_{L^2_x(\R^2)}+CN\delta_n\longrightarrow0.
\]
Since $\nabla P_{>N}u(t_n)-\nabla u(t_n)
=-\nabla P_{\leq N}u(t_n)$, the high-frequency limit follows.
\end{proof}

\subsection{Reduced Duhamel formula}
We now combine \eqref{eq:MC} with \eqref{eq:width} to show that the
free term in the Duhamel formula with terminal time $T$ converges
weakly to zero as $T\to\infty$. This gives the reduced Duhamel
formula used in the global exterior estimates of
Subsection~\ref{sec:dyadic_cone_completion}.

\begin{proposition}\label{prop:reduced_future_duhamel}
Let $u$ be a minimal-mass nonscattering solution.
Then, for every fixed $t\geq0$,
\begin{equation}\label{eq:linear_future_weak_zero}
 e^{i(t-T)\Delta}u(T)\rightharpoonup0
 \quad\hbox{in }H^1(\R^2)\qquad(T\to\infty).
\end{equation}
Consequently,
\begin{equation}\label{eq:reduced_future_duhamel}
 i\int_t^T e^{i(t-s)\Delta}F(u(s))\dd s
 \rightharpoonup u(t)
 \quad\hbox{in }H^1(\R^2)\qquad(T\to\infty).
\end{equation}
\end{proposition}

\begin{proof}
Fix $t\geq0$ and set $v_T:=e^{i(t-T)\Delta}u(T)$ for $T>t$.
Unitarity, conservation of mass, and \eqref{eq:uniform_H1_bound} give
$\norm{v_T}_{L^2_x(\R^2)}=\sqrt{m_Q}$ and $\norm{v_T}_{H^1}\leq B_H$.

We first prove that $v_T\rightharpoonup0$ in $L^2$.
Fix $g\in\mathcal S(\R^2)$ and $\varepsilon>0$, choose
$A_\varepsilon$ as in \eqref{eq:MC}, and set
$B_T:=B(c(T),A_\varepsilon\lambda(T))$.
Splitting $\langle u(T),e^{i(T-t)\Delta}g\rangle_{L^2}$ over
$B_T$ and its complement, Cauchy--Schwarz, \eqref{eq:MC}, and
the dispersive estimate yield
\[
 \begin{aligned}
 |\langle v_T,g\rangle_{L^2}|
 &\leq\sqrt{m_Q}\,|B_T|^{1/2}
          \norm{e^{i(T-t)\Delta}g}_{L^\infty_x(\R^2)}
       +\sqrt\varepsilon\norm{e^{i(T-t)\Delta}g}_{L^2_x(\R^2)}\\
 &\leq C\sqrt{m_Q}\,
       \frac{A_\varepsilon\lambda(T)}{T-t}\norm g_{L^1_x(\R^2)}
       +\sqrt\varepsilon\norm g_{L^2_x(\R^2)}.
 \end{aligned}
\]
By \eqref{eq:width},
$\lambda(T)/(T-t)\leq C_0\sqrt{1+T}/(T-t)\to0$.
Taking first $T\to\infty$ and then $\varepsilon\downarrow0$
gives $\langle v_T,g\rangle_{L^2}\to0$.
Since $\mathcal S(\R^2)$ is dense in $L^2$ and
$\norm{v_T}_{L^2_x(\R^2)}=\sqrt{m_Q}$, we obtain $v_T\rightharpoonup0$ in $L^2$.

To obtain weak convergence in $H^1$, take $\phi\in\mathcal S(\R^2)$.
Then
\[
 \langle v_T,\phi\rangle_{H^1}
 =\langle v_T,(1-\Delta)\phi\rangle_{L^2}\longrightarrow0.
\]
The bound $\norm{v_T}_{H^1}\leq B_H$ and density of
$\mathcal S(\R^2)$ in $H^1$ prove \eqref{eq:linear_future_weak_zero}.
Finally, for $T>t$, the Duhamel formula \eqref{eq:appendix_duhamel} gives
\[
 u(t)=e^{i(t-T)\Delta}u(T)
       +i\int_t^T e^{i(t-s)\Delta}F(u(s))\dd s
 \qquad\text{in }H^1(\R^2).
\]
Taking the weak $H^1$ limit and using
\eqref{eq:linear_future_weak_zero} proves \eqref{eq:reduced_future_duhamel}.
\end{proof}

\section{Exterior spacetime estimates}\label{sec:exterior_windows}

In this section, we establish exterior spacetime bounds for the
minimal-mass nonscattering solution, first on $[T,2T]$ and then
on the positive time half-line.

We begin with an interaction Morawetz estimate based on small
exterior mass. Together with the time-averaged kinetic energy
evacuation in \eqref{eq:ballistic_evacuation}, it gives decay of
the exterior Strichartz norms on $[T,2T]$ and of the localized
nonlinearity in $L_t^{4/3}W^{1,4/3}_x$ on these intervals.
Combining the exterior \(L_{t,x}^4\) and \(L_{t,x}^8\) decay estimates with two propagation estimates and the reduced Duhamel formula, we derive the global exterior \(L_{t,x}^4\) bound in Proposition~\ref{prop:conditional_exterior_global_L4}.

For estimates outside balls, we use fixed smooth cutoffs
$\chi\in C^\infty(\R^2;[0,1])$ with bounded derivatives of the
orders used below. For some $a>0$, they satisfy
\begin{equation}\label{eq:spatial_cutoff}
 \chi(x)=0\quad(|x|\leq a),\qquad \chi_T(x):=\chi(x/T).
\end{equation}

\subsection{Exterior interaction Morawetz estimate}
\label{subsec:exterior_im_finite}

We first bound the exterior $L_t^4L_x^8$ norm in terms of the
exterior mass. We use the two-dimensional correlation estimate
obtained by Colliander, Grillakis, and Tzirakis
\cite{CGT} and by Planchon and Vega \cite{PlanchonVega}.
The focusing cubic term is controlled by a finite sequence of
nested spatial cutoffs, while the quintic term has a nonnegative
contribution. Fix a nondecreasing smooth function
$\vartheta:\R\to[0,1]$ which is 0 on $(-\infty,0]$ and 1
on $[1,\infty)$.

\begin{proposition}[Exterior interaction Morawetz estimate]
\label{prop:exterior_im_finite}
Let $I=[t_0,t_1]$ have length $H$, and let
$u\in C(I;H^1(\R^2))$ solve \eqref{eq:cq} with
\[
 \sup_{t\in I}\|u(t)\|_{H^1}\leq B_H.
\]
Fix a unit vector $e\in\R^2$, $b\in\R$, and $L\geq4$. Set
\[
 \delta:=\sup_{t\in I}\|u(t)\|_{L_x^2(\{x:x\cdot e>b\})},\qquad
 \chi_0(x):=\vartheta\left(\frac{x\cdot e-b-L}{L}\right),\qquad
 w_0:=\chi_0u.
\]
There exists a constant $\delta_*>0$ such that, if
$\delta\leq\delta_*$, then
\begin{equation}\label{eq:exterior_im_finite}
 \int_I\bigl\|D^{1/2}|w_0(t)|^2\bigr\|_{L_x^2(\R^2)}^2\,dt
 \leq C_{B_H}\left(\delta^3+
          \frac{\delta^2 H}{L}+\frac{\delta H}{L^2}\right).
\end{equation}
Consequently,
\begin{equation}\label{eq:exterior_im_L4L8}
 \int_I\|u(t)\|_{L_x^8(\{x:x\cdot e>b+2L\})}^4\,dt
 \leq C_{B_H}\left(\delta^3+
          \frac{\delta^2 H}{L}+\frac{\delta H}{L^2}\right)
\end{equation}
and
\begin{equation}\label{eq:exterior_im_quintic_critical}
 \int_I\int_{x\cdot e>b+2L}|u(t,x)|^8\,dx\,dt
 \leq C_{B_H}\left(\delta^4+
          \frac{\delta^3 H}{L}+\frac{\delta^2 H}{L^2}\right).
\end{equation}
\end{proposition}

\begin{proof}
The case $\delta=0$ is immediate, so assume $\delta>0$.
We first derive the estimate for smooth solutions and justify
the passage to $H^1$ solutions in Step~4.

\proofstep{1}{A density estimate}
For every nonnegative $f\in L_x^1(\R^2)\cap\dot H^{1/2}(\R^2)$,
we claim that
\begin{equation}\label{eq:im_pressure_absorption}
 \int f(x)^2\bigl(|x|^{-1}*f\bigr)(x)\,dx
 \leq C\|f\|_{L_x^1(\R^2)}\bigl\|D^{1/2}f\bigr\|_{L_x^2(\R^2)}^2.
\end{equation}
Indeed, the Hardy--Littlewood--Sobolev inequality
bounds the left-hand side by
\[
 C\|f\|_{L_x^{4/3}(\R^2)}\|f\|_{L_x^{8/3}(\R^2)}^2.
\]
Interpolation gives
\begin{equation}\label{eq:im_density_interpolation}
 \|f\|_{L_x^{4/3}(\R^2)}\leq\|f\|_{L_x^1(\R^2)}^{2/3}\|f\|_{L_x^4(\R^2)}^{1/3},\qquad
 \|f\|_{L_x^{8/3}(\R^2)}^2\leq\|f\|_{L_x^1(\R^2)}^{1/3}\|f\|_{L_x^4(\R^2)}^{5/3}.
\end{equation}
Combining \eqref{eq:im_density_interpolation} with
$\|f\|_{L_x^4(\R^2)}\lesssim\|D^{1/2}f\|_{L_x^2(\R^2)}$ from
the Sobolev estimate \eqref{eq:appendix_sobolev} proves
\eqref{eq:im_pressure_absorption}.

\proofstep{2}{The localized interaction identity}
We next localize the interaction identity and estimate the terms
created by the spatial cutoff. For $\ell\geq1$ and $r>b$, set
\[
 \eta(x):=\vartheta\left(\frac{x\cdot e-r}{\ell}\right),\quad
 w:=\eta u,\quad \rho:=|u|^2,\quad \rho_w:=|w|^2=\eta^2\rho,
 \quad j_w:=\operatorname{Im}(\overline w\nabla w).
\]
Here $\rho$ and $\rho_w$ are the mass densities of $u$ and $w$.
Then
\begin{equation}\label{eq:im_cutoff_equation}
 iw_t+\Delta w=(-\rho+\rho^2)w+G_\eta,
 \qquad G_\eta:=2\nabla\eta\cdot\nabla u+(\Delta\eta)u.
\end{equation}
Since $\eta$ is supported in $x\cdot e>b$, we have
\begin{equation}\label{eq:im_cutoff_bounds}
 \|w\|_{L_x^2(\R^2)}\leq\delta,\qquad
 \|\nabla w\|_{L_x^2(\R^2)}\leq B_H+C\delta/\ell\leq C_{B_H},\qquad
 \|G_\eta\|_{L_x^2(\R^2)}\leq C(B_H/\ell+\delta/\ell^2).
\end{equation}

For $a(x):=|x|$, define
\[
 \mathcal M_w(t):=2\iint\rho_w(t,y)\nabla a(x-y)\cdot j_w(t,x)
                 \,dx\,dy.
\]
By \eqref{eq:im_cutoff_bounds},
\begin{equation}\label{eq:im_action_bound}
 |\mathcal M_w(t)|\leq2\|w(t)\|_{L_x^2(\R^2)}^3\|\nabla w(t)\|_{L_x^2(\R^2)}
 \leq C_{B_H}\delta^3.
\end{equation}
Set $V_w:=\nabla a*\rho_w$ and
\[
 \mathcal P(\rho):=-\tfrac12\rho^2+\tfrac23\rho^3.
\]
Differentiating $\mathcal M_w$ and using
\eqref{eq:im_cutoff_equation} gives the interaction Morawetz
identity \cite{CGT,PlanchonVega}
\begin{equation}\label{eq:im_localized_identity}
\begin{aligned}
 \mathcal M_w'(t)
 &=c_2\|D^{1/2}\rho_w\|_{L_x^2(\R^2)}^2+\mathcal K_w
     +2\int(\Delta a*\rho_w)\eta^2\mathcal P(\rho)\,dx\\
 &\quad+2\int V_w\cdot\nabla(\eta^2)\mathcal P(\rho)\,dx
     +\mathcal R_\eta,
\end{aligned}
\end{equation}
where $c_2>0$ is the  constant in
\[
 -\iint\rho_w(y)\Delta^2a(x-y)\rho_w(x)\,dx\,dy
 =c_2\|D^{1/2}\rho_w\|_{L_x^2(\R^2)}^2,
\]
and
\[
 \mathcal K_w:=4\sum_{j,k=1}^2\iint\partial_{jk}a(x-y)
 \left[\rho_w(y)\operatorname{Re}
   (\partial_j\overline w(x)\partial_kw(x))
      -j_{w,j}(x)j_{w,k}(y)\right]dx\,dy.
\]
The forcing contribution is
\begin{equation}\label{eq:im_forcing_term}
\begin{aligned}
 \mathcal R_\eta:={}&4\iint\operatorname{Im}(G_\eta(y)\overline w(y))
       \nabla a(x-y)\cdot j_w(x)\,dx\,dy\\
 &+2\operatorname{Re}\int V_w\cdot
       (G_\eta\nabla\overline w-w\nabla\overline{G_\eta})\,dx.
\end{aligned}
\end{equation}
The term $\mathcal K_w$ is nonnegative by symmetrization in $x,y$
and positivity of the Hessian of $a$. The nonlinear terms in
\eqref{eq:im_localized_identity} follow from
\[
 -\rho_w\nabla(-\rho+\rho^2)
 =-\nabla\bigl(\eta^2\mathcal P(\rho)\bigr)
                  +\mathcal P(\rho)\nabla(\eta^2).
\]
Since $\Delta a(x)=|x|^{-1}$, the first nonlinear integral in
\eqref{eq:im_localized_identity} equals
\[
 -\iint\frac{\rho_w(y)\eta(x)^2\rho(x)^2}{|x-y|}\,dx\,dy
 +\frac43\iint\frac{\rho_w(y)\eta(x)^2\rho(x)^3}{|x-y|}\,dx\,dy.
\]
The second term is nonnegative. Since $\|V_w\|_{L_x^\infty(\R^2)}\leq\delta^2$,
the cutoff term in \eqref{eq:im_localized_identity} is bounded by
\begin{equation}\label{eq:im_pressure_boundary}
 C\frac{\delta^2}{\ell}
       \int_{x\cdot e>b}(|u|^4+|u|^6)\,dx
 \leq C_{B_H}\frac{\delta^4}{\ell},
\end{equation}
where we used the Gagliardo--Nirenberg inequality \eqref{eq:appendix_gagliardo_nirenberg} on the
halfspace $x\cdot e>b$.

It remains to estimate $\mathcal R_\eta$ in \eqref{eq:im_forcing_term}.
Integrating the term containing $\nabla\overline{G_\eta}$ by parts gives
\begin{equation}\label{eq:im_forcing_ibp}
 \begin{aligned}
 &2\operatorname{Re}\int V_w\cdot
       (G_\eta\nabla\overline w-w\nabla\overline{G_\eta})\,dx\\
 &\qquad=4\operatorname{Re}\int V_w\cdot G_\eta\nabla\overline w\,dx
       +2\operatorname{Re}\int(\Delta a*\rho_w)G_\eta\overline w\,dx.
 \end{aligned}
\end{equation}
The Hardy--Littlewood--Sobolev and Gagliardo--Nirenberg
inequalities give
\[
 \|w(|x|^{-1}*\rho_w)\|_{L_x^2(\R^2)}
 \leq C\|w\|_{L_x^4(\R^2)}\|\rho_w\|_{L_x^{4/3}(\R^2)}
 \leq C\|w\|_{L_x^2(\R^2)}^2\|\nabla w\|_{L_x^2(\R^2)}.
\]
Thus Cauchy--Schwarz and \eqref{eq:im_forcing_ibp} imply
\begin{equation}\label{eq:im_forcing_L2}
 |\mathcal R_\eta|\leq C\|w\|_{L_x^2(\R^2)}^2\|\nabla w\|_{L_x^2(\R^2)}\|G_\eta\|_{L_x^2(\R^2)}.
\end{equation}
In particular, \eqref{eq:im_cutoff_bounds} yields
\begin{equation}\label{eq:im_linear_boundary}
 \int_I|\mathcal R_\eta(t)|\,dt
       \leq C_{B_H}\frac{\delta^2H}{\ell}.
\end{equation}

\proofstep{3}{Iteration over spatial cutoffs}
Choose an absolute $C_0$ sufficiently large for both
\eqref{eq:im_nested_recurrence} and
\eqref{eq:im_regularized_recurrence}, and set
\[
 q:=C_0\delta^2,\qquad \theta:=\sqrt q.
\]
Choose $\delta_*$ so that $2\sqrt{C_0}\delta_*\leq1/4$.
Let $J$ be the first positive integer such that $L\theta^J\leq1$,
and define
\[
 \ell_j:=L\theta^j\ (0\leq j<J),\qquad \ell_J:=1,\qquad
 r_0:=b+L,\quad r_{j+1}:=r_j-\ell_{j+1}\ (0\leq j<J).
\]
Set $\chi_j(x):=\vartheta((x\cdot e-r_j)/\ell_j)$,
$w_j:=\chi_ju$, and
\[
 \rho_j:=|w_j|^2,\qquad
 A_j:=\int_I\|D^{1/2}\rho_j(t)\|_{L_x^2(\R^2)}^2\,dt.
\]
Since $\sum_{j=1}^J\ell_j\leq L\theta/(1-\theta)+1<L$,
all cutoffs are supported in $x\cdot e>b$.
Moreover, $r_{j+1}+\ell_{j+1}=r_j$, so $\chi_{j+1}=1$ on
$\operatorname{supp}\chi_j$. Consequently,
$0\leq\rho_j\leq\rho_{j+1}$ and
$\chi_j^2\rho^2\leq\rho_{j+1}^2$. By
\eqref{eq:im_pressure_absorption},
\begin{equation}\label{eq:im_cubic_pressure_control}
 \begin{aligned}
 \iint\frac{\rho_j(y)\chi_j(x)^2\rho(x)^2}{|x-y|}\,dx\,dy
 &\leq\int\rho_{j+1}^2(|x|^{-1}*\rho_{j+1})\,dx\\
 &\leq C\delta^2\|D^{1/2}\rho_{j+1}\|_{L_x^2(\R^2)}^2.
 \end{aligned}
\end{equation}
Integrating \eqref{eq:im_localized_identity} on $I$, discarding
$\mathcal K_{w_j}$ and the nonnegative quintic term, and using
\eqref{eq:im_action_bound}, \eqref{eq:im_pressure_boundary},
\eqref{eq:im_linear_boundary}, and
\eqref{eq:im_cubic_pressure_control}, we obtain
\begin{equation}\label{eq:im_nested_recurrence}
 A_j\leq qA_{j+1}
       +C_{B_H}\left(\delta^3+\frac{\delta^2H}{\ell_j}\right),
 \qquad 0\leq j<J.
\end{equation}
For every $0\leq j\leq J$, the fractional product rule,
Sobolev embedding, and the Gagliardo--Nirenberg inequality \eqref{eq:appendix_gagliardo_nirenberg}
give
\begin{equation}\label{eq:im_density_apriori}
 \|D^{1/2}|w_j|^2\|_{L_x^2(\R^2)}
 \lesssim\|w_j\|_{L_x^4(\R^2)}\|D^{1/2}w_j\|_{L_x^4(\R^2)}
 \lesssim\delta^{1/2}\|\nabla w_j\|_{L_x^2(\R^2)}^{3/2},\qquad
 A_j\leq C_{B_H}\delta H.
\end{equation}
In particular, $A_J\leq C_{B_H}\delta H$.
Since $q/\theta=\theta<1$ and $q^J=\theta^{2J}\leq L^{-2}$,
iteration of \eqref{eq:im_nested_recurrence} yields
\begin{equation}\label{eq:im_iterated_bound}
\begin{aligned}
 A_0&\leq C_{B_H}\delta^3\sum_{j=0}^{J-1}q^j
       +C_{B_H}\frac{\delta^2H}{L}
                    \sum_{j=0}^{J-1}(q/\theta)^j+q^JA_J\\
 &\leq C_{B_H}\left(\delta^3+\frac{\delta^2H}{L}
                              +\frac{\delta H}{L^2}\right).
\end{aligned}
\end{equation}
This proves \eqref{eq:exterior_im_finite} for smooth solutions.
Sobolev embedding applied to $|w_0|^2$ gives
\eqref{eq:exterior_im_L4L8}. Finally,
the Gagliardo--Nirenberg inequality \eqref{eq:appendix_gagliardo_nirenberg} on the halfspace gives
\begin{equation}\label{eq:exterior_im_uniform_L8}
 \sup_{t\in I}\|u(t)\|_{L_x^8(\{x:x\cdot e>b+2L\})}^4
 \leq C\delta B_H^3.
\end{equation}
Multiplying \eqref{eq:exterior_im_uniform_L8} by \eqref{eq:exterior_im_L4L8}
proves \eqref{eq:exterior_im_quintic_critical}.

\proofstep{4}{Approximation in $H^1$}
Fix $\delta<\delta'<2\delta_*$. By
Lemma~\ref{lem:local_small_data}, there exist smooth solutions
$u^{(n)}\to u$ in $C(I;H^1)$ such that, for all sufficiently
large $n$,
\[
 \sup_I\|u^{(n)}\|_{H^1}\leq B_H+1,\qquad
 \sup_I\|u^{(n)}\|_{L_x^2(\{x:x\cdot e>b\})}\leq\delta'.
\]
The constants below are uniform in $n$ and depend only on $B_H$;
we continue to denote them by $C_{B_H}$.
Construct the cutoffs in Step~3 using
$q=C_0(\delta')^2$ and $\theta=\sqrt q$.

For each fixed $n$, use
$a_\varepsilon(x):=\sqrt{|x|^2+\varepsilon^2}$ in the identity
\eqref{eq:im_localized_identity}, with the first term replaced by
$-\iint\rho_w(y)\Delta^2a_\varepsilon(x-y)\rho_w(x)\,dx\,dy$. Its Hessian is positive semidefinite, and
\[
 |\nabla a_\varepsilon|\leq1,\qquad
 0\leq\Delta a_\varepsilon(x)\leq2|x|^{-1}\quad(x\ne0).
\]
Thus \eqref{eq:im_action_bound},
\eqref{eq:im_pressure_boundary},
\eqref{eq:im_linear_boundary}, and
\eqref{eq:im_cubic_pressure_control} hold uniformly in
$\varepsilon$, with $\delta'$ in place of $\delta$.
For $0\leq j\leq J$, set
\[
 \rho_j^{(n)}:=|\chi_ju^{(n)}|^2,\qquad
 A_j^{(n)}:=\int_I\|D^{1/2}\rho_j^{(n)}(t)\|_{L_x^2(\R^2)}^2\,dt.
\]
By \eqref{eq:im_density_apriori},
$A_j^{(n)}\leq C_{B_H}\delta'H<\infty$.
Since $\Delta a_\varepsilon(x)\leq2|x|^{-1}$,
\eqref{eq:im_cubic_pressure_control} bounds the integrated cubic
term by $C(\delta')^2A_{j+1}^{(n)}$, uniformly in $\varepsilon$.
The convolution operator with kernel $-\Delta^2a_\varepsilon$ has
multiplier $c_2|\xi|(1+\varepsilon|\xi|)e^{-\varepsilon|\xi|}$.
Thus, after integrating \eqref{eq:im_localized_identity} on $I$
with $a_\varepsilon$ in place of $a$ and discarding the nonnegative
kinetic and quintic terms, we obtain
\begin{equation}\label{eq:im_regularized_recurrence}
\begin{aligned}
 &\int_I\int_{\R^2}|\xi|(1+\varepsilon|\xi|)e^{-\varepsilon|\xi|}
       |\widehat{\rho_j^{(n)}}(t,\xi)|^2\,d\xi\,dt\\
 &\qquad\leq qA_{j+1}^{(n)}
  +C_{B_H}\left((\delta')^3+
             \frac{(\delta')^2H}{\ell_j}\right),\qquad 0\leq j<J.
\end{aligned}
\end{equation}
The multiplier on the left increases to $|\xi|$ as
$\varepsilon\downarrow0$. Monotone convergence therefore gives
\begin{equation}\label{eq:im_approximant_recurrence}
 A_j^{(n)}\leq qA_{j+1}^{(n)}
  +C_{B_H}\left((\delta')^3+
             \frac{(\delta')^2H}{\ell_j}\right).
\end{equation}
Iterating \eqref{eq:im_approximant_recurrence} as in the derivation of \eqref{eq:im_iterated_bound}
yields
\begin{equation}\label{eq:im_smooth_uniform_bound}
 A_0^{(n)}\leq C_{B_H}\left((\delta')^3+
        \frac{(\delta')^2H}{L}+\frac{\delta'H}{L^2}\right).
\end{equation}
To let $n\to\infty$, use the fractional product rule and
the Sobolev embeddings \eqref{eq:appendix_H1_embeddings} to obtain
\[
 \|D^{1/2}(|f|^2-|g|^2)\|_{L_x^2(\R^2)}
 \leq C(\|f\|_{H^1}+\|g\|_{H^1})\|f-g\|_{H^1}.
\]
Since $\chi_0u^{(n)}\to\chi_0u$ in $C(I;H^1)$, this gives
\[
 \sup_{t\in I}\bigl\|D^{1/2}
       (|\chi_0u^{(n)}(t)|^2-|\chi_0u(t)|^2)\bigr\|_{L_x^2(\R^2)}
       \longrightarrow0,
\]
and hence $A_0^{(n)}\to A_0$. The cutoff $\chi_0$ is independent
of $\delta'$. Letting first $n\to\infty$ and then
$\delta'\downarrow\delta$ in \eqref{eq:im_smooth_uniform_bound}
proves \eqref{eq:exterior_im_finite}. Sobolev embedding
then gives \eqref{eq:exterior_im_L4L8}; combining \eqref{eq:exterior_im_L4L8} with
the Gagliardo--Nirenberg inequality \eqref{eq:appendix_gagliardo_nirenberg} gives
\eqref{eq:exterior_im_quintic_critical}.
\end{proof}

\begin{corollary}[Decay of the exterior $L_t^4L_x^8$ norm]
\label{cor:exterior_im_doubling}
Let $u\in C([0,\infty);H^1(\R^2))$ solve \eqref{eq:cq}, with
$\sup_{t\geq0}\|u(t)\|_{H^1}\leq B_H$. Suppose that, for every
$c>0$,
\[
 \sup_{t\in[T,2T]}\|u(t)\|_{L_x^2(\{x:|x|>cT\})}\longrightarrow0
 \qquad(T\to\infty).
\]
Then, for every $c>0$,
\begin{equation}\label{eq:exterior_im_doubling}
 \int_T^{2T}\|u(t)\|_{L_x^8(\{x:|x|>ct\})}^4\,dt\longrightarrow0
 \qquad(T\to\infty).
\end{equation}
The exterior $L_{t,x}^8$ norm on
$\{(t,x):T\leq t\leq2T,\ |x|>ct\}$ also tends to zero.
\end{corollary}

\begin{proof}
Fix $0<c_1<c/(3\sqrt2)$ and set
$\delta_T:=\sup_{t\in[T,2T]}\|u(t)\|_{L_x^2(\{x:|x|>c_1T\})}$.
Let $e_1:=(1,0)$ and $e_2:=(0,1)$ be the standard coordinate vectors.
Apply Proposition~\ref{prop:exterior_im_finite} in each of
the four directions $e=e_1,-e_1,e_2,-e_2$ with $H=T$ and
$L=b=c_1T$. For sufficiently large $T$, we have
$\delta_T\leq\delta_*$ and $L\geq4$, so
\eqref{eq:exterior_im_L4L8} gives
\[
 \int_T^{2T}\|u(t)\|_{L_x^8(\{x:x\cdot e>3c_1T\})}^4\,dt
 \leq C_{B_H}\left(\delta_T^3+
       \frac{\delta_T^2}{c_1}+\frac{\delta_T}{c_1^2T}\right)
 \longrightarrow0.
\]
For $T\leq t\leq2T$, the inequality
$\max(|x_1|,|x_2|)\geq|x|/\sqrt2$ shows that
$\{|x|>ct\}$ is contained in the union of these four halfspaces.
Summing their bounds proves \eqref{eq:exterior_im_doubling}.
The $L_{t,x}^8$ conclusion follows from
\eqref{eq:exterior_im_quintic_critical} on the same halfspaces.
\end{proof}

\subsection{Exterior Strichartz estimates}
\label{subsec:cubic_window_closure}

We next estimate the exterior $L_{t,x}^4$ norm by partitioning
time so that the cubic term can be absorbed on each subinterval.
We then combine this estimate with
Corollaries~\ref{cor:ballistic_evacuation} and
\ref{cor:exterior_im_doubling} to prove decay of the exterior
Strichartz norms on $[T,2T]$.

\begin{proposition}[Exterior $L_{t,x}^4$ estimate]
\label{prop:cubic_window_closure}
Fix $0<a<b$, let $T>0$, and set $I_T:=[T,2T]$. Let
$u\in C(I_T;H^1(\R^2))$ solve \eqref{eq:cq}, with
$\sup_{t\in I_T}\|u(t)\|_{H^1}\leq B_H$.
Define $\delta,k,\varepsilon_T\geq0$ by
\[
 \delta:=\sup_{t\in I_T}\|u(t)\|_{L_x^2(\{x:|x|>aT\})},\qquad
 k^2:=\frac1T\int_{I_T}\int_{aT<|x|<bT}|\nabla u|^2\,dx\,dt,
\]
and
\[
 \varepsilon_T:=\|\mathbf1_{\{|x|>aT\}}|u|^4u\|_{L_{t,x}^{4/3}(I_T\times\R^2)}.
\]
There exist $\delta_*>0$ and $T_*\geq1$, depending at most on
$B_H,a,b$, such that, if $\delta\leq\delta_*$ and $T\geq T_*$, then
\begin{equation}\label{eq:cubic_window_closure}
 \|u\|_{L_{t,x}^4(I_T\times\{x:|x|>bT\})}
                   \leq C_{B_H,a,b}(\delta+k+\varepsilon_T).
\end{equation}
\end{proposition}

\begin{proof}
The case $\delta=0$ is immediate. Assume $\delta>0$.

We first derive an estimate for two smooth spatial cutoffs
$\chi_j,\chi_{j+1}$ taking values in $[0,1]$ and satisfying
\begin{equation}\label{eq:cubic_window_cutoff_support}
 \begin{gathered}
 \operatorname{supp}\chi_j\cup\operatorname{supp}\chi_{j+1}
       \subset\{|x|>aT\},\\
 \chi_{j+1}=1\quad\text{on }\operatorname{supp}\chi_j.
 \end{gathered}
\end{equation}
Set
\[
 w_j:=\chi_ju,\qquad Y_j:=\|w_j\|_{L_{t,x}^4(I_T\times\R^2)}^4,\qquad
 G_j:=2\nabla\chi_j\cdot\nabla u+(\Delta\chi_j)u,\qquad
 f_j:=\chi_j|u|^4u.
\]
Then
\begin{equation}\label{eq:critical_window_actual_equation}
 (i\partial_t+\Delta)w_j=-|u|^2w_j+G_j+f_j.
\end{equation}
Choose $\eta>0$ so that the Strichartz constant in
\eqref{eq:appendix_strichartz} for the pair $(4,4)$ satisfies
$C\eta^2\leq1/2$. Starting at $T$, end each subinterval when
the integral of $|w_{j+1}|^4$ reaches $\eta^4$, or at $2T$ if
the remaining integral is smaller. Since $Y_{j+1}<\infty$ by
the Sobolev embeddings \eqref{eq:appendix_H1_embeddings}, this gives a finite partition
$I_T=\bigcup_{m=1}^{K_j}I_{j,m}$ satisfying
\begin{equation}\label{eq:cubic_window_time_partition}
 \|w_{j+1}\|_{L_{t,x}^4(I_{j,m}\times\R^2)}\leq\eta,\qquad
 K_j\leq1+\eta^{-4}Y_{j+1}.
\end{equation}
At the left endpoint $t_{j,m}$ of each interval,
$\|w_j(t_{j,m})\|_{L_x^2(\R^2)}\leq\delta$. Moreover, $u=w_{j+1}$ on
$\operatorname{supp}w_j$, so H\"older's inequality gives
\[
 \||u|^2w_j\|_{L_{t,x}^{4/3}(I_{j,m}\times\R^2)}
 \leq\eta^2\|w_j\|_{L_{t,x}^4(I_{j,m}\times\R^2)}.
\]
Applying the Strichartz estimate \eqref{eq:appendix_strichartz} to
\eqref{eq:critical_window_actual_equation} and absorbing
the cubic term, we obtain
\begin{equation}\label{eq:cubic_window_local_estimate}
 \|w_j\|_{L_{t,x}^4(I_{j,m}\times\R^2)}
 \leq C\left(\delta+\|G_j\|_{L_t^1(I_{j,m};L_x^2(\R^2))}
                  +\|f_j\|_{L_{t,x}^{4/3}(I_{j,m}\times\R^2)}\right).
\end{equation}
For $G_j$ and $f_j$, we have
\begin{align*}
 \sum_m\|G_j\|_{L_t^1(I_{j,m};L_x^2(\R^2))}^4
 &\leq\|G_j\|_{L_t^1(I_T;L_x^2(\R^2))}^4,\\
 \sum_m\|f_j\|_{L_{t,x}^{4/3}(I_{j,m}\times\R^2)}^4
 &\leq\left(\sum_m\|f_j\|_{L_{t,x}^{4/3}(I_{j,m}\times\R^2)}^{4/3}\right)^3
 =\|f_j\|_{L_{t,x}^{4/3}(I_T\times\R^2)}^4.
\end{align*}
Taking fourth powers in \eqref{eq:cubic_window_local_estimate}
and summing over $m$, only the initial-data term contributes
$K_j\delta^4$. Thus \eqref{eq:cubic_window_time_partition} gives
\begin{equation}\label{eq:critical_window_recurrence}
 Y_j\leq C\delta^4(1+Y_{j+1})
      +C\|G_j\|_{L_t^1(I_T;L_x^2(\R^2))}^4
      +C\|f_j\|_{L_{t,x}^{4/3}(I_T\times\R^2)}^4.
\end{equation}

We next construct cutoffs satisfying
\eqref{eq:cubic_window_cutoff_support} for which
\eqref{eq:critical_window_recurrence} can be iterated. Fix $C_0$ larger than the constant in \eqref{eq:critical_window_recurrence},
and put
\[
 q:=C_0\delta^4,\qquad \theta:=q^{1/16},\qquad L:=\frac{b-a}{2}T.
\]
Choose $\delta_*$ sufficiently small and $T_*$ sufficiently large
so that $\theta\leq1/4$ and $L\geq4$. Let $J$ be the first positive
integer such that $L\theta^J\leq1$, and define
\[
 \ell_j:=L\theta^j\ (0\leq j<J),\quad \ell_J:=1,\quad
 r_0:=aT+L,\quad r_{j+1}:=r_j-\ell_{j+1}\ (0\leq j<J).
\]
Using the function $\vartheta$ fixed in
Section~\ref{subsec:exterior_im_finite}, set
$\chi_j(x):=\vartheta((|x|-r_j)/\ell_j)$.
Since $\sum_{j=1}^J\ell_j\leq L\theta/(1-\theta)+1<L$, we have
$r_j>aT$. Also, $r_{j+1}+\ell_{j+1}=r_j$ and $r_0+\ell_0=bT$.
Thus \eqref{eq:cubic_window_cutoff_support} holds,
$\chi_0=1$ on $|x|\geq bT$, and all cutoff derivatives vanish
outside $aT<|x|<bT$.

For these radial cutoffs,
\[
 \|\nabla\chi_j\|_{L_x^\infty(\R^2)}\leq C\ell_j^{-1},\qquad
 \|\Delta\chi_j\|_{L_x^\infty(\R^2)}
 \leq C\left(\ell_j^{-2}+(aT\ell_j)^{-1}\right)
 \leq C_{a,b}\ell_j^{-2},
\]
where the last inequality uses $1\leq\ell_j\leq L$.
Cauchy--Schwarz inequality in time gives
\[
 \int_{I_T}\|\nabla u(t)\|_{L_x^2(\{x:aT<|x|<bT\})}\,dt
 \leq T^{1/2}\left(\int_{I_T}\int_{aT<|x|<bT}
                       |\nabla u|^2\,dx\,dt\right)^{1/2}=Tk.
\]
Consequently,
\begin{equation}\label{eq:critical_window_G_bound}
 \|G_j\|_{L_t^1(I_T;L_x^2(\R^2))}
 \leq C_{a,b}\left(\frac{Tk}{\ell_j}
                         +\frac{T\delta}{\ell_j^2}\right).
\end{equation}
Since $0\leq\chi_j\leq1$ and $\chi_j=0$ on $|x|\leq aT$,
\begin{equation}\label{eq:critical_window_Q_uniform}
 \|f_j\|_{L_{t,x}^{4/3}(I_T\times\R^2)}\leq \varepsilon_T,
 \qquad 0\leq j\leq J.
\end{equation}

It remains to iterate the estimate. Combining
\eqref{eq:critical_window_recurrence},
\eqref{eq:critical_window_G_bound}, and
\eqref{eq:critical_window_Q_uniform}, we obtain
\begin{equation}\label{eq:cubic_window_spatial_recurrence}
 Y_j\leq qY_{j+1}+C(\delta^4+\varepsilon_T^4)
       +C_{a,b}\frac{T^4k^4}{\ell_j^4}
       +C_{a,b}\frac{T^4\delta^4}{\ell_j^8},\qquad 0\leq j<J.
\end{equation}
At the last cutoff, $\ell_J=1$,
$\|w_J\|_{L_x^2(\R^2)}\leq\delta$, and $\|\nabla w_J\|_{L_x^2(\R^2)}\leq B_H+C\delta$.
By the Gagliardo--Nirenberg inequality \eqref{eq:appendix_gagliardo_nirenberg},
\[
 Y_J\leq C\int_{I_T}\|w_J(t)\|_{L_x^2(\R^2)}^2\|\nabla w_J(t)\|_{L_x^2(\R^2)}^2\,dt
 \leq C_{B_H}T\delta^2.
\]
The definition of $J$ gives $q^J\leq L^{-16}$, while $J\geq1$
gives $q^J\leq q$. Multiplying these inequalities gives
$q^{2J}\leq qL^{-16}$, hence $q^J\leq\sqrt q\,L^{-8}$. Therefore
\begin{equation}\label{eq:critical_window_terminal}
 q^JY_J\leq C_{B_H}T\delta^2\sqrt q\,L^{-8}
          \leq C_{B_H,a,b}\delta^4T^{-7}.
\end{equation}
Since $q/\theta^4=q^{3/4}$ and $q/\theta^8=q^{1/2}$ are
uniformly less than 1, iteration of
\eqref{eq:cubic_window_spatial_recurrence} yields
\begin{align*}
 Y_0&\leq C(\delta^4+\varepsilon_T^4)\sum_{j=0}^{J-1}q^j
       +C_{a,b}\frac{T^4k^4}{L^4}
                    \sum_{j=0}^{J-1}(q/\theta^4)^j\\
 &\quad+C_{a,b}\frac{T^4\delta^4}{L^8}
                    \sum_{j=0}^{J-1}(q/\theta^8)^j+q^JY_J\\
 &\leq C_{B_H,a,b}(\delta^4+k^4+\varepsilon_T^4),
\end{align*}
where we used \eqref{eq:critical_window_terminal},
$L=(b-a)T/2$, and $T\geq1$ in the last line.
Since $\chi_0=1$ on $|x|>bT$, taking fourth roots proves
\eqref{eq:cubic_window_closure}.
\end{proof}

\begin{corollary}[Decay of the exterior Strichartz norms]
\label{cor:cubic_window_unconditional}
Let $u$ be a minimal-mass nonscattering solution.
Then, for every $c>0$,
\begin{equation}\label{eq:cubic_window_unconditional}
 \|u\|_{L_{t,x}^4(\{(t,x):T\leq t\leq2T,\ |x|>ct\})}\longrightarrow0
 \qquad(T\to\infty).
\end{equation}
Moreover, for every fixed smooth radial cutoff $\chi$ as in
\eqref{eq:spatial_cutoff} and every fixed Schr\"odinger-admissible
pair $(q,r)$, with $(q,r)\ne(2,\infty)$,
\begin{equation}\label{eq:exterior_all_S0}
 \|\chi_Tu\|_{L_t^\infty([T,2T];L_x^2(\R^2))}
 +\|\chi_Tu\|_{L_t^q([T,2T];L_x^r(\R^2))}\longrightarrow0
 \qquad(T\to\infty).
\end{equation}
\end{corollary}

\begin{proof}
We first verify that the three quantities on the right-hand side
of \eqref{eq:cubic_window_closure} tend to zero.
Fix $0<a<b$ and set $I_T:=[T,2T]$. Write $\delta_T,k_T$ for the
quantities $\delta,k$ in Proposition~\ref{prop:cubic_window_closure}.
Applying \eqref{eq:ballistic_mass} and
\eqref{eq:ballistic_evacuation} with terminal time $2T$ gives
\begin{equation}\label{eq:cubic_window_mass_kinetic_decay}
 \begin{aligned}
 \delta_T^2&\leq\sup_{0\leq t\leq2T}
                  \int_{|x|>aT}|u(t,x)|^2\,dx\longrightarrow0,\\
 k_T^2&\leq2\left[\frac1{2T}\int_0^{2T}
                  \int_{|x|<bT}|\nabla u(t,x)|^2\,dx\,dt\right]
                  \longrightarrow0,
 \end{aligned}
\end{equation}
where the corresponding speeds are $a/2$ and $b/2$.

To estimate $\varepsilon_T$, choose a smooth radial cutoff
$\eta:\R^2\to[0,1]$ which is zero on $|x|\leq a_0$ and one
on $|x|\geq a$, where $0<a_0<a$. Set $v_T:=\eta(x/T)u$ and
$\widetilde\delta_T:=\|v_T\|_{L_t^\infty(I_T;L_x^2(\R^2))}$.
By \eqref{eq:ballistic_mass}, $\widetilde\delta_T\to0$.
For $t\in I_T$, we have $v_T(t,x)=0$ on $|x|\leq a_0T$ and
\[
 \{|x|>a_0T\}\subset\{|x|>(a_0/2)t\}.
\]
Hence \eqref{eq:exterior_im_doubling} gives
\begin{equation}\label{eq:cubic_window_L4L8_decay}
 A_T:=\int_{I_T}\|v_T(t)\|_{L_x^8(\R^2)}^4\,dt\longrightarrow0.
\end{equation}
The functions $v_T$ are uniformly bounded in $L_t^\infty(I_T;H_x^1(\R^2))$.
Thus the Gagliardo--Nirenberg inequality \eqref{eq:appendix_gagliardo_nirenberg} yields
\[
 \|v_T(t)\|_{L_x^{16/3}(\R^2)}^2
 \leq C\|v_T(t)\|_{L_x^2(\R^2)}^{3/4}\|\nabla v_T(t)\|_{L_x^2(\R^2)}^{5/4}
 \leq C_{B_H,\eta}\widetilde\delta_T^{3/4}.
\]
As $v_T=u$ on $|x|>aT$, H\"older's inequality gives
\begin{equation}\label{eq:quintic_N0_uniform}
 \varepsilon_T\leq\|v_T\|_{L_t^4(I_T;L_x^8(\R^2))}^3
              \|v_T\|_{L_t^\infty(I_T;L_x^{16/3}(\R^2))}^2
 \leq C_{B_H,\eta}A_T^{3/4}\widetilde\delta_T^{3/4}
 \longrightarrow0.
\end{equation}
By \eqref{eq:cubic_window_mass_kinetic_decay},
\eqref{eq:quintic_N0_uniform}, and
Proposition~\ref{prop:cubic_window_closure},
$\|u\|_{L_{t,x}^4(I_T\times\{x:|x|>bT\})}\to0$.
Taking $b=c$ and using $t\geq T$ proves
\eqref{eq:cubic_window_unconditional}.

We next prove \eqref{eq:exterior_all_S0}.
First choose a smooth radial cutoff $\psi:\R^2\to[0,1]$
which is zero on $|x|\leq a$ and one on $|x|\geq b$, and set
$\psi_T(x):=\psi(x/T)$. Then $\|\psi_Tu(T)\|_{L_x^2(\R^2)}\leq\delta_T$.
Since $\|\nabla\psi_T\|_{L_x^\infty(\R^2)}\leq C_\psi/T$ and
$\|\Delta\psi_T\|_{L_x^\infty(\R^2)}\leq C_\psi/T^2$, the same
Cauchy--Schwarz estimate as in \eqref{eq:critical_window_G_bound}
gives
\begin{equation}\label{eq:cubic_window_commutator_decay}
 \|[\Delta,\psi_T]u\|_{L_t^1(I_T;L_x^2(\R^2))}
 \leq C_\psi\left(k_T+\frac{\delta_T}{T}\right)\longrightarrow0.
\end{equation}
For $t\in I_T$, we have
$\{|x|>aT\}\subset\{|x|>(a/2)t\}$.
Therefore \eqref{eq:cubic_window_unconditional}, with
$c=a/2$, and \eqref{eq:quintic_N0_uniform} imply
\begin{equation}\label{eq:cubic_window_nonlinearity_decay}
 \|\psi_TF(u)\|_{L_{t,x}^{4/3}(I_T\times\R^2)}
 \leq\|u\|_{L_{t,x}^4(I_T\times\{x:|x|>aT\})}^3+\varepsilon_T
 \longrightarrow0.
\end{equation}
Applying the Strichartz estimate \eqref{eq:appendix_strichartz} to the equation for
$\psi_Tu$, with initial time $T$, and using
\eqref{eq:cubic_window_mass_kinetic_decay},
\eqref{eq:cubic_window_commutator_decay}, and
\eqref{eq:cubic_window_nonlinearity_decay}, proves
\eqref{eq:exterior_all_S0} with $\psi_T$ in place of
$\chi_T$, for each fixed admissible pair.
For a cutoff $\chi$ in the statement, choose $0<a<b$ such that
$\chi=0$ on $|x|\leq b$, and choose $\psi$ as above.
Then $\psi=1$ on $\operatorname{supp}\chi$, and
$|\chi_Tu|\leq|\psi_Tu|$ proves \eqref{eq:exterior_all_S0}.
\end{proof}

\subsection{Derivative estimates for the exterior nonlinearity}
\label{subsec:source}

We now combine the exterior Strichartz estimates in
\eqref{eq:exterior_all_S0} with the uniform $H^1$ bound
to control the localized nonlinearity with one spatial derivative.

\begin{corollary}[Exterior nonlinear estimates]
\label{cor:actual_exterior_cubic_N1}
Let $u$ be a minimal-mass nonscattering solution, and let
$\chi$ be a fixed smooth cutoff as in \eqref{eq:spatial_cutoff}.
Then, as $T\to\infty$,
\begin{equation}\label{eq:actual_exterior_cubic_N1}
 \|\chi_T|u|^2u\|_{N^1([T,2T])}
                                               \longrightarrow0.
\end{equation}
Moreover,
\begin{equation}\label{eq:actual_exterior_full_N1}
 \|\chi_T|u|^4u\|_{N^1([T,2T])}
 +\|\chi_TF(u)\|_{N^1([T,2T])}
                                               \longrightarrow0.
\end{equation}
\end{corollary}

\begin{proof}
Set $I_T:=[T,2T]$. Choose $a>0$ such that $\chi=0$ on
$|x|\leq a$, and a smooth radial cutoff $\eta:\R^2\to[0,1]$
which is zero on $|x|\leq a/2$ and one on $|x|\geq a$.
Then $\eta_T(x):=\eta(x/T)$ equals 1 on
$\operatorname{supp}\chi_T\cup\operatorname{supp}\nabla\chi_T$.
By \eqref{eq:exterior_all_S0}, with the admissible pair
$(8/3,8)$, and \eqref{eq:uniform_H1_bound},
\[
 \alpha_T:=\|\eta_Tu\|_{L_t^{8/3}(I_T;L_x^8(\R^2))}\longrightarrow0,
 \qquad \|\eta_Tu\|_{L_t^\infty(I_T;H_x^1(\R^2))}\leq C_\eta B_H
 \quad(T\geq1).
\]
The product rule and $\|\nabla\chi_T\|_{L_x^\infty(\R^2)}
=T^{-1}\|\nabla\chi\|_{L_x^\infty(\R^2)}$ give
\begin{equation}\label{eq:exterior_nonlinearity_product_rule}
 \begin{aligned}
 &\bigl|\chi_T|u|^m u\bigr|
       +\bigl|\nabla(\chi_T|u|^m u)\bigr|\\
 &\qquad\leq C_\chi(1+T^{-1})|\eta_Tu|^m
                       (|u|+|\nabla u|),\qquad m=2,4.
 \end{aligned}
\end{equation}

For the cubic term, place the two undifferentiated factors in
$L_t^{8/3}(I_T;L_x^8(\R^2))$ and use the $L_t^\infty(I_T;L_x^2(\R^2))$ bounds
for $u$ and $\nabla u$. H\"older's inequality and
\eqref{eq:exterior_nonlinearity_product_rule} yield
\[
 \|\chi_T|u|^2u\|_{N^1(I_T)}
 \leq C_\chi(1+T^{-1})\alpha_T^2
                 \|u\|_{L_t^\infty(I_T;H_x^1(\R^2))}
 \longrightarrow0,
\]
which proves \eqref{eq:actual_exterior_cubic_N1}.

For the quintic term, we first use
the Gagliardo--Nirenberg inequality \eqref{eq:appendix_gagliardo_nirenberg} with $q=8$ and $p=16$:
\[
 \|\eta_Tu(t)\|_{L_x^{16}(\R^2)}^4
 \leq C\|\nabla(\eta_Tu(t))\|_{L_x^2(\R^2)}^2\|\eta_Tu(t)\|_{L_x^8(\R^2)}^2.
\]
Taking the $L_t^{4/3}(I_T)$ norm gives
\begin{equation}\label{eq:exterior_quintic_factor_bound}
 \begin{aligned}
 \|\eta_Tu\|_{L_t^{16/3}(I_T;L_x^{16}(\R^2))}^4
 &\leq C\|\nabla(\eta_Tu)\|_{L_t^\infty(I_T;L_x^2(\R^2))}^2
                 \|\eta_Tu\|_{L_t^{8/3}(I_T;L_x^8(\R^2))}^2\\
 &\leq C_{B_H,\eta}\alpha_T^2\longrightarrow0.
 \end{aligned}
\end{equation}
Thus, H\"older's inequality
and \eqref{eq:exterior_nonlinearity_product_rule}--
\eqref{eq:exterior_quintic_factor_bound} give
\[
 \|\chi_T|u|^4u\|_{N^1(I_T)}
 \leq C_\chi(1+T^{-1})
       \|\eta_Tu\|_{L_t^{16/3}(I_T;L_x^{16}(\R^2))}^4
       \|u\|_{L_t^\infty(I_T;H_x^1(\R^2))}
 \longrightarrow0.
\]
Together with \eqref{eq:actual_exterior_cubic_N1} and the
triangle inequality, this proves \eqref{eq:actual_exterior_full_N1}.
\end{proof}

\subsection{Global exterior spacetime bounds}
\label{subsec:dyadic_nested_cone_S0}
\label{sec:dyadic_cone_completion}
\label{subsec:cone_summation}

We now extend the estimates on \([T,2T]\) to the positive time half-line.
To this end, we first state two estimates for Duhamel integrals over separated time
intervals and between moving half-spaces.

\begin{lemma}
\label{lem:dyadic_offdiagonal_S0}
Let $T_0\geq1$, $T_j:=2^jT_0$, and $I_j:=[T_j,2T_j]$ for
$j\geq0$. For $k\geq0$, let $G_k\in N^0(I_k)$ be a function
on $I_k\times\R^2$, extended by zero outside $I_k$.
If $|j-k|\geq2$, then
\begin{equation}\label{eq:dyadic_offdiagonal_bound}
 \left\|\int_{I_k}e^{i(t-s)\Delta}G_k(s)\,ds\right\|_{L_{t,x}^4(I_j\times\R^2)}
 \leq C2^{-|j-k|/4}\|G_k\|_{N^0(I_k)}.
\end{equation}
For $|j-k|\leq1$, the integrals over $I_k\cap\{s>t\}$ and
$I_k\cap\{s<t\}$ each satisfy \eqref{eq:dyadic_offdiagonal_bound}.
For a nonnegative sequence $b=(b_k)_{k\geq0}$, define each
component $(\mathcal Tb)_j\in[0,+\infty]$ by
\begin{equation}\label{eq:dyadic_sequence_operator}
 (\mathcal Tb)_j:=\sum_{k\geq0}2^{-|j-k|/4}b_k,\qquad
 \|\mathcal T\|_{\ell^p\to\ell^p}
 \leq \sum_{m\in\mathbb Z}2^{-|m|/4}
 \quad(1\leq p\leq\infty).
\end{equation}
Now let $G$ be a function on $[T_0,\infty)\times\R^2$ with
$G|_{I_k}\in N^0(I_k)$ for every $k\geq0$, and set
\[
 G_k:=\mathbf1_{I_k}G,\qquad
 b_k:=\|G_k\|_{N^0(I_k)}=\|G\|_{N^0(I_k)}.
\]
For each $j$ such that
$(\mathcal Tb)_j<\infty$, the limit
\begin{equation}\label{eq:dyadic_terminal_integral}
 \int_t^\infty e^{i(t-s)\Delta}G(s)\,ds
 :=\lim_{T\to\infty}\int_t^T e^{i(t-s)\Delta}G(s)\,ds
\end{equation}
exists strongly in $L_{t,x}^4(I_j\times\R^2)$, and
\begin{equation}\label{eq:dyadic_duhamel_output}
 \begin{aligned}
 &\left\|\int_{T_0}^t e^{i(t-s)\Delta}G(s)\,ds
       \right\|_{L_{t,x}^4(I_j\times\R^2)}\\
 &\quad+\left\|\int_t^\infty e^{i(t-s)\Delta}G(s)\,ds
       \right\|_{L_{t,x}^4(I_j\times\R^2)}
 \leq C(\mathcal Tb)_j.
 \end{aligned}
\end{equation}
In particular, \eqref{eq:dyadic_terminal_integral} and
\eqref{eq:dyadic_duhamel_output} hold for every $j$ when
$b\in\ell^\infty$.
\end{lemma}

\begin{lemma}
\label{lem:all_frequency_gap_crossing_S0}
Fix $0<A<B$, let $\delta:=B-A\leq1$, and let $T_0\geq1$.
For $e\in\mathbb S^1$, use the notation in
\eqref{eq:directional_notation}. Let $\varphi(t,x)$ and $\psi(t,x)$ be smooth
cutoff functions satisfying, for every $t\geq T_0$,
\begin{equation}\label{eq:gap_crossing_cutoffs}
 \begin{gathered}
 \operatorname{supp}\varphi(t)\subset\{x\cdot e\leq At\},\qquad
 \operatorname{supp}\psi(t)\subset\{x\cdot e\geq Bt\},\\
 |\varphi|+|\psi|\leq C,\qquad
 |\nabla_x\varphi(t)|\leq\frac C{\delta t},\qquad
 |\nabla_x\psi(t)|\leq\frac C{\delta t}.
 \end{gathered}
\end{equation}
Choose a fixed smooth $\eta:\R\to[0,1]$ such that $\eta(v)=0$ for
$v\leq-1/8$ and $\eta(v)=1$ for $v\geq1/8$, and define
\begin{equation}\label{eq:gap_crossing_multipliers}
 \mu:=\frac{A+B}{4},\qquad
 P^+:=\eta\!\left(\frac{D_e-\mu}{\delta}\right),\qquad P^-:=I-P^+.
\end{equation}
Here $\eta((D_e-\mu)/\delta)$ denotes the Fourier multiplier defined on
$L_x^2(\R^2)$ by
\[
 \widehat{P^+f}(\xi):=\eta\!\left(\frac{\xi\cdot e-\mu}{\delta}\right)
                         \widehat f(\xi).
\]
Then, for every fixed $1<p<\infty$, these complementary multipliers
are bounded on $L_x^{p}(\R^2)$ independently of $\delta,\mu$, and $e$.

For $\mathcal G\in L_t^\infty([T_0,\infty);H_x^{3/4}(\R^2))$, set
\begin{align*}
 C^+(t)&:=\lim_{T\to\infty}\psi(t)\int_t^T
                    e^{i(t-s)\Delta}P^+\varphi(s)\mathcal G(s)\,ds,\\
 C^-(t)&:=\psi(t)\int_{T_0}^t
                    e^{i(t-s)\Delta}P^-\varphi(s)\mathcal G(s)\,ds.
\end{align*}
The limit $C^+(t)$ exists strongly in $H^{3/4}$ for
each fixed $t\geq T_0$, and $C^-(t)\in H^{3/4}$. Moreover,
\begin{equation}\label{eq:all_frequency_gap_crossing_S0}
 \|C^+\|_{L_{t,x}^4([T_0,\infty)\times\R^2)}
 +\|C^-\|_{L_{t,x}^4([T_0,\infty)\times\R^2)}
 \leq C \delta^{-4}\|\mathcal G\|_{L_t^\infty([T_0,\infty);H_x^{3/4}(\R^2))}.
\end{equation}
The constant depends only on $\eta$ and the constant in
\eqref{eq:gap_crossing_cutoffs}.
\end{lemma}

\begin{remark}\label{rem:gap_crossing_velocities}
The parameters $A$ and $B$ are the speeds at which the lines
$x\cdot e=At$ and $x\cdot e=Bt$ move in direction $e$.
Their distance at time $t$ is $\delta t$, where $\delta=B-A$.
The parameter $\delta$ also determines the transition width of the
multipliers in \eqref{eq:gap_crossing_multipliers}.
The symbols of $P^+$ and $P^-$ are supported in
\[
 \{2\xi_e-A\geq \delta/4\}
 \quad\text{and}\quad
 \{B-2\xi_e\geq \delta/4\},
\]
respectively. The free Schr\"odinger evolution at frequency
$\xi$ has velocity $2\xi_e$ in direction $e$.
These inequalities make the derivative of the phase with
respect to $\xi_e$ nonzero on the spatial supports in
\eqref{eq:gap_crossing_cutoffs}: the first applies to $C^+$
for $s>t$, and the second to $C^-$ for $s<t$.
They are used in the integration by parts in
Appendix~\ref{app:propagation_proofs}.
\end{remark}

The proofs of Lemmas~\ref{lem:dyadic_offdiagonal_S0} and
\ref{lem:all_frequency_gap_crossing_S0} are given in
Appendix~\ref{app:propagation_proofs}. We now combine these
estimates with the reduced Duhamel formula to prove the global
exterior $L_{t,x}^4$ bound.

\begin{proposition}[Global exterior $L_{t,x}^4$ bound]
\label{prop:conditional_exterior_global_L4}
Let $u$ be a minimal-mass nonscattering solution satisfying
\eqref{eq:reduced_future_duhamel}. Assume that, for every $c>0$,
\begin{equation}\label{eq:local_exterior_L4_L8_assumption}
 \|u\|_{L_{t,x}^4(\{(t,x):T<t<2T,\ |x|>ct\})}
 +\|u\|_{L_{t,x}^8(\{(t,x):T<t<2T,\ |x|>ct\})}
                  \longrightarrow0\qquad(T\to\infty).
\end{equation}
Then, for every $c>0$ and sufficiently large $T_0$,
\begin{equation}\label{eq:global_exterior_L4}
 \|u\|_{L_{t,x}^4(\{(t,x):t>T_0,\ |x|>ct\})}<\infty.
\end{equation}
\end{proposition}
\begin{proof}
Fix $e\in\mathbb S^1$ and $a_0>0$. We first prove the bound
on $x\cdot e>a_0t$. Choose $0<h<\min(a_0,1)$, put
$a_*:=a_0-h>0$, and define
\begin{equation}\label{eq:nested_cone_speeds}
 a_n:=a_*+h2^{-n},\qquad
 g_n:=a_n-a_{n+1}=h2^{-n-1},\qquad n\geq0.
\end{equation}
Choose smooth cutoffs $0\leq\varphi_n,\psi_n\leq1$ such that
\begin{equation}\label{eq:nested_cone_cutoffs}
 \begin{aligned}
 \varphi_n(t,x)&=\begin{cases}
 0,&x\cdot e\leq a_{n+1}t,\\
 1,&x\cdot e\geq(a_{n+1}+g_n/4)t,
 \end{cases}\\
 \psi_n(t,x)&=\begin{cases}
 0,&x\cdot e\leq(a_n-g_n/4)t,\\
 1,&x\cdot e\geq a_nt,
 \end{cases}
 \end{aligned}
\end{equation}
with first spatial derivatives bounded by $C/(g_nt)$. Split
\begin{equation}\label{eq:nested_cone_decomposition}
 G_n:=\varphi_nF(u),\qquad K_n:=(1-\varphi_n)F(u).
\end{equation}
The cutoffs $1-\varphi_n$ and $\psi_n$ satisfy
\eqref{eq:gap_crossing_cutoffs} with
\begin{equation}\label{eq:nested_cone_gap}
 A_n:=a_{n+1}+g_n/4,\qquad B_n:=a_n-g_n/4,\qquad
 \delta_n:=B_n-A_n=g_n/2.
\end{equation}
Let $P_n^\pm$ be the multipliers in
\eqref{eq:gap_crossing_multipliers} for these speeds.

For $T_j:=2^jT_0$ and $I_j:=[T_j,2T_j]$, set
\begin{equation}\label{eq:nested_cone_norms}
 \begin{aligned}
 X_{n,j}&:=\|u\|_{L_{t,x}^4(\{(t,x):t\in I_j,\ x\cdot e>a_nt\})},\\
 \varepsilon&:=\sup_{j\geq0}\left(
 \|u\|_{L_{t,x}^4(\{(t,x):t\in I_j,\ x\cdot e>a_*t\})}^2
 +\|u\|_{L_{t,x}^8(\{(t,x):t\in I_j,\ x\cdot e>a_*t\})}^4\right).
 \end{aligned}
\end{equation}
Fix $\varepsilon_*>0$, to be chosen sufficiently small below.
 By
\eqref{eq:local_exterior_L4_L8_assumption} with $c=a_*$,
there exists $\overline T\geq1$ such that
\[
 \|u\|_{L_{t,x}^4(\{(t,x):S<t<2S,\ |x|>a_*t\})}^2
 +\|u\|_{L_{t,x}^8(\{(t,x):S<t<2S,\ |x|>a_*t\})}^4
 \leq\varepsilon_*\qquad(S\geq\overline T).
\]
Choose $T_0\geq\overline T$. For every $n,j\geq0$, we have
$T_j\geq\overline T$ and, for $t\in I_j$,
\[
 \{x\cdot e>a_nt\}\subset\{x\cdot e>a_*t\}
                         \subset\{|x|>a_*t\}.
\]
Hence \eqref{eq:nested_cone_norms} gives
\begin{equation}\label{eq:nested_cone_uniform_smallness}
 \varepsilon\leq\varepsilon_*,\qquad
 \sup_{n\geq0}\|X_n\|_{\ell^\infty}\leq\varepsilon^{1/2}.
\end{equation}
The lower bound $T_0\geq\overline T$ is determined by the
estimate at speed $a_*$ and the target $\varepsilon_*$,
independently of $n,j$. For every $n\geq0$,
\eqref{eq:nested_cone_cutoffs} and \eqref{eq:nested_cone_gap}
verify \eqref{eq:gap_crossing_cutoffs} on $[T_0,\infty)$ with
$\varphi=1-\varphi_n$, $\psi=\psi_n$, $A=A_n$, and $B=B_n$.
Since $0<A_n<B_n$ and $0<\delta_n\leq1$, the cutoff and speed
conditions in Lemma~\ref{lem:all_frequency_gap_crossing_S0}
are satisfied for every $n\geq0$.

Where $G_n\ne0$, we have $x\cdot e>a_{n+1}t>a_*t$. H\"older's inequality gives
\begin{equation}\label{eq:nested_cone_source_coefficient}
 \begin{aligned}
 \|G_n\|_{N^0(I_j)}
 &\leq\left(
 \|u\|_{L_{t,x}^4(\{(t,x):t\in I_j,\ x\cdot e>a_*t\})}^2
 +\|u\|_{L_{t,x}^8(\{(t,x):t\in I_j,\ x\cdot e>a_*t\})}^4\right)X_{n+1,j}\\
 &\leq\varepsilon X_{n+1,j}.
 \end{aligned}
\end{equation}
We next estimate the contributions of $K_n$.
The fractional product estimate
\eqref{eq:appendix_fractional_product} and the embeddings
\eqref{eq:appendix_H1_embeddings} give
\begin{equation}\label{eq:nonlinearity_uniform_H34}
 \begin{aligned}
 \|F(u(t))\|_{H^{3/4}}
 &\leq C\bigl(\|u(t)\|_{L_x^{16}(\R^2)}^2+\|u(t)\|_{L_x^{32}(\R^2)}^4\bigr)
                    \|u(t)\|_{W^{3/4,8/3}}\\
 &\leq C_{B_H}.
 \end{aligned}
\end{equation}
Let $C_n^\pm$ denote the two fields in
Lemma~\ref{lem:all_frequency_gap_crossing_S0}, with
$\psi=\psi_n$, $\varphi=1-\varphi_n$, and $\mathcal G=F(u)$.
By \eqref{eq:all_frequency_gap_crossing_S0} and
\eqref{eq:nonlinearity_uniform_H34},
\begin{equation}\label{eq:nested_cone_crossing_bound}
 \|C_n^+\|_{L_{t,x}^4([T_0,\infty)\times\R^2)}
 +\|C_n^-\|_{L_{t,x}^4([T_0,\infty)\times\R^2)}
 \leq C_{B_H}\delta_n^{-4}.
\end{equation}

We next derive the Duhamel formula
\eqref{eq:nested_cone_duhamel_identity} for $\psi_nu$.
Fix $n$ and a compact interval $J=[t_1,t_2]\subset(T_0,\infty)$.
For $T>t_2$, the finite Duhamel formula gives
\begin{equation}\label{eq:nested_cone_finite_duhamel}
 \begin{aligned}
 \psi_nP_n^+u(t)
 ={}&\psi_ne^{i(t-T)\Delta}P_n^+u(T)\\
 &+i\psi_n(t)\int_t^T e^{i(t-s)\Delta}P_n^+G_n(s)\,ds\\
 &+i\psi_n(t)\int_t^T e^{i(t-s)\Delta}P_n^+K_n(s)\,ds.
 \end{aligned}
\end{equation}
We pass to the limit in these three terms separately.

By \eqref{eq:nested_cone_source_coefficient},
\eqref{eq:nested_cone_uniform_smallness}, and the uniform
$L_x^{4/3}(\R^2)$ bound for $P_n^+$,
\[
 \sup_{k\geq0}\|P_n^+G_n\|_{N^0(I_k)}
 \leq C\varepsilon^{3/2}<\infty.
\]
Lemma~\ref{lem:dyadic_offdiagonal_S0} shows that the
integral of $P_n^+G_n$ in \eqref{eq:nested_cone_finite_duhamel}
converges strongly in $L_{t,x}^4(J\times\R^2)$ as $T\to\infty$.

For the term containing $K_n$, use
\eqref{eq:gap_crossing_H34} with decay exponent $2$ and
\eqref{eq:nonlinearity_uniform_H34}. For $T'>T>t_2$, they give
\begin{equation}\label{eq:nested_cone_terminal_tail}
 \begin{aligned}
 &\sup_{t\in J}\left\|
 \psi_n(t)\int_T^{T'}e^{i(t-s)\Delta}P_n^+K_n(s)\,ds
                         \right\|_{H^{3/4}}\\
 &\quad\leq C_{B_H}(1+\delta_n^{-1})^{3/4}
                       \int_T^{T'}(1+\delta_n^2s)^{-2}\,ds\\
 &\quad\leq C_{B_H}\delta_n^{-2}(1+\delta_n^{-1})^{3/4}
                       (1+\delta_n^2T)^{-1}\longrightarrow0.
 \end{aligned}
\end{equation}
Here $n$ is fixed. Thus the third term in
\eqref{eq:nested_cone_finite_duhamel} converges to $iC_n^+$
in $L_t^\infty(J;H_x^{3/4}(\R^2))$.

 Let
$\theta\in C_c^\infty((t_1,t_2)\times\R^2)$. In the $L_x^2(\R^2)$
pairing, using \eqref{eq:linear_future_weak_zero} we have
\[
 \begin{aligned}
 &\left\langle\psi_ne^{i(t-T)\Delta}P_n^+u(T),\theta(t)\right\rangle\\
 &\qquad=\left\langle e^{i(t-T)\Delta}u(T),
                   (P_n^+)^*(\psi_n\theta(t))\right\rangle
 \longrightarrow0
 \end{aligned}
\]
for each fixed $t\in J$.
The test function on the right is in $L_x^2(\R^2)$ and is independent
of $T$, so this follows from the stated weak $H^1$ convergence.
The absolute value of the pairing is bounded by \(CB_H\|\theta(t)\|_{L_x^2(\R^2)}\).
Dominated convergence therefore yields convergence to zero upon integration over J.

We may now take the limit in \eqref{eq:nested_cone_finite_duhamel} against \(\theta\).
Combining the limiting form of \eqref{eq:nested_cone_finite_duhamel} with
the Duhamel formula \eqref{eq:appendix_duhamel} for $P_n^-u$ with initial time
$T_0$, and using $P_n^++P_n^-=I$, gives
\begin{equation}\label{eq:nested_cone_duhamel_identity}
 \begin{aligned}
 \psi_nu={}&\psi_ne^{i(t-T_0)\Delta}P_n^-u(T_0)
       -i\psi_n\int_{T_0}^t e^{i(t-s)\Delta}P_n^-G_n(s)\,ds\\
 &+i\psi_n\int_t^\infty e^{i(t-s)\Delta}P_n^+G_n(s)\,ds
       +iC_n^+-iC_n^-.
 \end{aligned}
\end{equation}
The identity holds in spacetime distributions on
$(T_0,\infty)\times\R^2$, since $J$ was arbitrary.

The free term in \eqref{eq:nested_cone_duhamel_identity} has
global $L_{t,x}^4([T_0,\infty)\times\R^2)$ norm at most $CB_H$ by
the Strichartz estimate \eqref{eq:appendix_strichartz}. Define $f_{n,j}$ as the sum
of its $L_{t,x}^4(I_j\times\R^2)$ norm and those of $C_n^+,C_n^-$.
Then \eqref{eq:nested_cone_gap} and
\eqref{eq:nested_cone_crossing_bound} give
\begin{equation}\label{eq:nested_cone_f_bound}
 f_n\geq0,\qquad
 \|f_n\|_{\ell^4}\leq C_{B_H}(1+\delta_n^{-4})
                      \leq C_{B_H,h}2^{4n}.
\end{equation}
Since $\psi_n=1$ on $x\cdot e>a_nt$, the uniform
$L_x^{4/3}(\R^2)$ bounds of $P_n^\pm$,
\eqref{eq:nested_cone_source_coefficient}, and
Lemma~\ref{lem:dyadic_offdiagonal_S0} yield
\begin{equation}\label{eq:nested_cone_sequence}
 X_n\leq f_n+\rho\mathcal TX_{n+1},\qquad
 \rho:=C\varepsilon.
\end{equation}
The constant in $\rho$ is independent of $n,T_0$.

Choose $\varepsilon_*$, and then $T_0$ as in the construction of \eqref{eq:nested_cone_uniform_smallness}, so that
\begin{equation}\label{eq:nested_cone_contraction}
 16\rho\sum_{\ell\in\mathbb Z}2^{-|\ell|/4}<1.
\end{equation}
Iterating \eqref{eq:nested_cone_sequence} gives
\[
 X_0\leq\sum_{n=0}^{m-1}\rho^n\mathcal T^nf_n
                  +\rho^m\mathcal T^mX_m.
\]
First use \eqref{eq:dyadic_sequence_operator} in $\ell^\infty$:
\[
 \|\rho^m\mathcal T^mX_m\|_{\ell^\infty}
 \leq\left(\rho\sum_{\ell\in\mathbb Z}2^{-|\ell|/4}\right)^m
             \sup_n\|X_n\|_{\ell^\infty}\longrightarrow0.
\]
The remaining series converges in $\ell^4$, since
\eqref{eq:nested_cone_f_bound} and
\eqref{eq:nested_cone_contraction} imply
\[
 \sum_{n\geq0}\|\rho^n\mathcal T^nf_n\|_{\ell^4}
 \leq C_{B_H,h}\sum_{n\geq0}
                  \left(16\rho\sum_{\ell\in\mathbb Z}2^{-|\ell|/4}\right)^n<\infty.
\]
Its terms are nonnegative, so its sum dominates $X_0$
componentwise. Thus $X_0\in\ell^4$, and
\begin{equation}\label{eq:global_halfspace_L4}
 \int_{T_0}^\infty\int_{x\cdot e>a_0t}|u(t,x)|^4\,dx\,dt
 =\sum_{j\geq0}X_{0,j}^4<\infty.
\end{equation}

Fix $c>0$, set $a_0=c/2$, and choose
\[
 e_1:=(1,0),\qquad e_2:=(-1,0),\qquad
 e_3:=(0,1),\qquad e_4:=(0,-1).
\]
If $|x|>ct$, then
$\max_{1\leq k\leq4}x\cdot e_k
 =\max(|x_1|,|x_2|)\geq|x|/\sqrt2>a_0t$. Hence
\begin{equation}\label{eq:radial_halfspace_cover}
 \{x:|x|>ct\}\subset
       \bigcup_{k=1}^4\{x:x\cdot e_k>a_0t\},\qquad t>0.
\end{equation}
For each $k$, \eqref{eq:global_halfspace_L4} gives a starting
time $T^{(k)}\geq1$ for which the integral in direction $e_k$
is finite. Put $T_*:=\max_{1\leq k\leq4}T^{(k)}$. Then
\[
 \begin{aligned}
 \int_{T_*}^\infty\int_{|x|>ct}|u|^4
 &\leq\sum_{k=1}^4\int_{T_*}^\infty
                              \int_{x\cdot e_k>a_0t}|u|^4\\
 &\leq\sum_{k=1}^4\int_{T^{(k)}}^\infty
                              \int_{x\cdot e_k>a_0t}|u|^4<\infty.
 \end{aligned}
\]
This proves \eqref{eq:global_exterior_L4} for every
$T_0\geq T_*$.
\end{proof}

\begin{remark}
\label{rem:global_exterior_application}
The $L_{t,x}^4$ limit in \eqref{eq:local_exterior_L4_L8_assumption}
is \eqref{eq:cubic_window_unconditional}. To obtain
the $L_{t,x}^8$ limit, fix $c>0$, set $a=c/2$, and use the directions
$e_1,\ldots,e_4$ from \eqref{eq:radial_halfspace_cover}.
Write $E_{e,a}(t):=\{x:x\cdot e>at\}$ and $I_T:=[T,2T]$.
Applying the Gagliardo--Nirenberg inequality \eqref{eq:appendix_GN_L8} to $u(t)$ on each halfspace,
using $\|\nabla u(t)\|_{L_x^2(\R^2)}\leq B_H$, and integrating in time gives
\begin{equation}\label{eq:dyadic_halfspace_L8}
 \begin{aligned}
 \int_{I_T}\int_{|x|>ct}|u|^8
 &\leq\sum_{k=1}^4\int_{I_T}\int_{E_{e_k,a}(t)}|u|^8\\
 &\leq CB_H^4\sum_{k=1}^4
                       \int_{I_T}\int_{E_{e_k,a}(t)}|u|^4\\
 &\leq4CB_H^4\int_{I_T}\int_{|x|>at}|u|^4
 \longrightarrow0.
 \end{aligned}
\end{equation}
Here each $E_{e_k,a}(t)$ is contained in $\{|x|>at\}$, and
the last limit follows from \eqref{eq:cubic_window_unconditional}
at speed $a$. The limits \eqref{eq:cubic_window_unconditional}
and \eqref{eq:dyadic_halfspace_L8}, together with
\eqref{eq:reduced_future_duhamel}, verify the hypotheses of
Proposition~\ref{prop:conditional_exterior_global_L4}.

\end{remark}

\section{Exterior kinetic energy decay}\label{sec:gradient_rigidity}

In this section, we prove that
$\|\nabla u(t)\|_{L_x^2(\{x:|x|>ct\})}\longrightarrow0$ as $t\to\infty$
for every fixed $c>0$. We first prove decay of the normal derivative on $x\cdot e>Bt$,
assuming that a suitable cutoff of $F(u)$ belongs to
$N^1([T_0,\infty))$. We then verify this assumption using the
global exterior $L_{t,x}^4$ bound \eqref{eq:global_exterior_L4}
and the derivative estimate \eqref{eq:actual_exterior_full_N1}
on $[T,2T]$. Finally, we combine the estimates in finitely many
directions to obtain decay of the full gradient on $|x|>ct$.

\subsection{Estimates for the exterior normal derivative}
\label{subsec:normal_crossing_tail_rigidity}
\label{subsec:incoming}

We begin with two propagation estimates and use them to prove
the conditional decay statement in
Proposition~\ref{cor:global_exterior_source_normal_rigidity}.
For $e\in\mathbb S^1$, use \eqref{eq:directional_notation} and define
$\mathcal H_e^{1/4}:=L^2_{x_e}(\R;H^{1/4}_{x_\perp}(\R))$, with norm
\[
 \|f\|_{\mathcal H_e^{1/4}}
 :=\left(\int_{\R}
       \|f(x_e,\cdot)\|_{H^{1/4}_{x_\perp}(\R)}^2\,dx_e\right)^{1/2}.
\]
Plancherel's theorem gives
$\|f\|_{\mathcal H_e^{1/4}}\leq C\|f\|_{H^{1/4}(\R^2)}$.
Together with \eqref{eq:nonlinearity_uniform_H34}, this implies,
for every uniformly $H^1$-bounded solution and every $T_0>0$,
\begin{equation}\label{eq:normal_nonlinearity_bound}
 \sup_{t\geq T_0}\|F(u(t))\|_{\mathcal H_e^{1/4}}
 \leq C\sup_{t\geq T_0}\|F(u(t))\|_{H^{3/4}}
 \leq C_{B_H}.
\end{equation}
We will use \eqref{eq:normal_nonlinearity_bound} to apply
Lemma~\ref{lem:normal_crossing_full_frequency} with
$\mathcal G=F(u)$.

\begin{lemma}
\label{lem:normal_crossing_full_frequency}
Let $0<A<B$, $\delta:=B-A$, and $T_0\geq1$.
Choose $\eta$ as in Lemma~\ref{lem:all_frequency_gap_crossing_S0}
and define $\mu,P^\pm$ by \eqref{eq:gap_crossing_multipliers}.
Let $\varphi(t,x_e)$ be measurable, with
$|\varphi|\leq1$ and
$\operatorname{supp}\varphi(t)\subset\{x_e\leq At\}$
for $t\geq T_0$, and let
$\mathcal G\in L_t^\infty([T_0,\infty);L_x^2(\R^2))$.
For $m\in\{0,1\}$, $t\geq T_0$, and $T>t$, set
$\Omega_t:=\{x:x_e>Bt\}$ and write
\begin{align*}
 C_{m,T}^+(t)&:=\mathbf1_{\Omega_t}
       \int_t^T\partial_e^m e^{i(t-s)\Delta}P^+
                         \varphi(s)\mathcal G(s)\,ds,\\
 C_m^-(t)&:=\mathbf1_{\Omega_t}
       \int_{T_0}^t\partial_e^m e^{i(t-s)\Delta}P^-
                         \varphi(s)\mathcal G(s)\,ds.
\end{align*}
Here each Duhamel integral is taken in $H^{-m}(\R^2)$.
On $\Omega_t$, it agrees with an $L_x^2(\Omega_t)$ function.
We define $C_{m,T}^+(t)$ and $C_m^-(t)$ by setting the respective
functions equal to 0 outside $\Omega_t$.
For each fixed $t\geq T_0$, the limit
\[
 C_m^+(t):=\lim_{T\to\infty}C_{m,T}^+(t)
\]
exists strongly in $L_x^2(\R^2)$, and
\begin{equation}\label{eq:normal_crossing_pointwise_L2}
 \|C_m^\pm(t)\|_{L_x^2(\R^2)}
 \leq C(1+A+B)^m(1+\delta^{-1})^6t^{-2}
             \|\mathcal G\|_{L_t^\infty([T_0,\infty);L_x^2(\R^2))}.
\end{equation}
If $\mathcal G\in L_t^\infty([T_0,\infty);\mathcal H_e^{1/4})$,
the same limit is strong in $L_x^4(\R^2)$, and
\begin{equation}\label{eq:normal_crossing_pointwise_L4}
 \|C_m^\pm(t)\|_{L_x^4(\R^2)}
 \leq C(1+A+B)^m(1+\delta^{-1})^7t^{-9/4}
             \|\mathcal G\|_{L_t^\infty([T_0,\infty);\mathcal H_e^{1/4})}.
\end{equation}
Consequently,
\begin{equation}\label{eq:normal_crossing_global_L4}
 \|C_m^\pm\|_{L_{t,x}^4([T_0,\infty)\times\R^2)}
 \leq C(1+A+B)^m(1+\delta^{-1})^7T_0^{-2}
             \|\mathcal G\|_{L_t^\infty([T_0,\infty);\mathcal H_e^{1/4})}.
\end{equation}
The bounds \eqref{eq:normal_crossing_pointwise_L2},
\eqref{eq:normal_crossing_pointwise_L4}, and
\eqref{eq:normal_crossing_global_L4} also hold for $C_{m,T}^+$
on $T_0\leq t<T$, uniformly in $T$. The constant $C$ depends only on $\eta$.
\end{lemma}

The proof is given in Appendix~\ref{app:propagation_proofs}.
We also need the following estimate for the free evolution.

\begin{lemma}
\label{lem:incoming_free_velocity}
Let $B>0$ and $\varepsilon>0$. If $f\in H^1(\R^2)$ satisfies
$\operatorname{supp}\widehat f\subset\{\xi_e\leq B/2-\varepsilon\}$,
then
\[
 \|e^{it\Delta}f\|_{L_x^2(\{x:x_e>Bt\})}+
 \|e^{it\Delta}\partial_ef\|_{L_x^2(\{x:x_e>Bt\})}
 \longrightarrow0\qquad(t\to\infty).
\]
\end{lemma}
\begin{proof}
For $g\in L_x^2(\R^2)$ and $t>0$, the factorization of the free
propagator gives, with a unitary Fourier transform,
\[
 \|e^{it\Delta}g\|_{L_x^2(\{x:x_e>Bt\})}
 =\left\|\mathbf1_{\{\xi_e>B/2\}}
       \mathcal F\bigl(e^{i|x|^2/(4t)}g\bigr)\right\|_{L_\xi^2(\R^2)}.
\]
If $\widehat g$ vanishes on $\{\xi_e>B/2\}$, Plancherel's
theorem and dominated convergence yield
\[
 \|e^{it\Delta}g\|_{L_x^2(\{x:x_e>Bt\})}
 \leq\|(e^{i|x|^2/(4t)}-1)g\|_{L_x^2(\R^2)}\longrightarrow0.
\]
Apply this to $g=f$ and $g=\partial_ef$.
\end{proof}

We now combine the two lemmas with the reduced Duhamel formula.

\begin{proposition}[Decay of the exterior normal derivative]
\label{cor:global_exterior_source_normal_rigidity}
Let $u$ be a uniformly $H^1$-bounded solution of \eqref{eq:cq} satisfying
\eqref{eq:reduced_future_duhamel}. Fix $T_0>0$ and $0<a<A<B$, and choose
$\vartheta\in C^\infty(\R;[0,1])$ equal to 0 on
$(-\infty,a]$ and 1 on $[A,\infty)$. Set
$\zeta(t,x):=\vartheta(x_e/t)$, and assume
\begin{equation}\label{eq:normal_exterior_source_assumption}
 G:=\zeta F(u)\in N^1([T_0,\infty)).
\end{equation}
Then
\begin{equation}\label{eq:normal_exterior_energy_pointwise_decay}
 \|\partial_eu(t)\|_{L_x^2(\{x:x_e>Bt\})}\longrightarrow0
 \qquad(t\to\infty).
\end{equation}
\end{proposition}
\begin{proof}
By increasing \(T_0\) if necessary, we may assume \(T_0\geq 1\).
Set $\delta:=B-A$, $\Omega_t:=\{x:x_e>Bt\}$, and use the multipliers
in \eqref{eq:gap_crossing_multipliers}. Write
$F(u)=G+K$, where $K:=(1-\zeta)F(u)$. Since
$\operatorname{supp}(1-\zeta(t))\subset\{x_e\leq At\}$,
Lemma~\ref{lem:normal_crossing_full_frequency} applies
with $\varphi=1-\zeta$ and $\mathcal G=F(u)$. By
\eqref{eq:normal_crossing_pointwise_L2} and
\eqref{eq:normal_nonlinearity_bound}, its functions $C_1^\pm$ satisfy
\begin{equation}\label{eq:normal_interior_duhamel_decay}
 \|C_1^+(t)\|_{L_x^2(\R^2)}+\|C_1^-(t)\|_{L_x^2(\R^2)}
 \leq C_{B_H}(1+A+B)(1+\delta^{-1})^6t^{-2}.
\end{equation}
For the integrals containing $G$, the dual Strichartz estimate
\eqref{eq:appendix_dual_strichartz} gives
\begin{equation}\label{eq:normal_exterior_duhamel_tail}
 \left\|\int_t^\infty e^{i(t-s)\Delta}P^\pm G(s)\,ds\right\|_{H^1}
 \leq C\|G\|_{N^1([t,\infty))}\longrightarrow0.
\end{equation}
These integrals converge strongly in $H^1$.

We first obtain a formula for $\partial_eP^+u$ on $\Omega_t$.
For fixed $t\geq T_0$ and $T>t$, the finite Duhamel formula gives
\[
 \begin{aligned}
 \partial_eP^+u(t)
 ={}&\partial_e e^{i(t-T)\Delta}P^+u(T)
       +i\int_t^T e^{i(t-s)\Delta}P^+\partial_eG(s)\,ds\\
    &+i\int_t^T\partial_e e^{i(t-s)\Delta}P^+K(s)\,ds.
 \end{aligned}
\]
The first term tends weakly to zero in $L_x^2(\R^2)$ by
\eqref{eq:reduced_future_duhamel} and the finite Duhamel formula.
The second converges strongly in $L_x^2(\R^2)$ by
\eqref{eq:normal_exterior_duhamel_tail}; the third converges
strongly in $L_x^2(\Omega_t)$ by
Lemma~\ref{lem:normal_crossing_full_frequency}.
Testing against $\phi\in C_c^\infty(\Omega_t)$ and taking
$T\to\infty$, we obtain the $L_x^2(\Omega_t)$ identity
\begin{equation}\label{eq:normal_outgoing_identity}
 \mathbf1_{\Omega_t}\partial_eP^+u(t)
 =i\mathbf1_{\Omega_t}\int_t^\infty
       e^{i(t-s)\Delta}P^+\partial_eG(s)\,ds+iC_1^+(t).
\end{equation}

For $P^-u$, define
\[
 f_-:=e^{-iT_0\Delta}P^-u(T_0)
       -i\int_{T_0}^\infty e^{-is\Delta}P^-G(s)\,ds.
\]
The integral converges strongly in $H^1$ by the dual
Strichartz estimate \eqref{eq:appendix_dual_strichartz}, and
\begin{equation}\label{eq:normal_incoming_fourier_support}
 \operatorname{supp}\widehat f_-
 \subset\{\xi_e\leq\mu+\delta/8\}
 =\{\xi_e\leq B/2-\delta/8\}.
\end{equation}
The Duhamel formula with initial time $T_0$ now gives
\begin{equation}\label{eq:normal_incoming_identity}
 \begin{aligned}
 \mathbf1_{\Omega_t}\partial_eP^-u(t)
 ={}&\mathbf1_{\Omega_t}e^{it\Delta}\partial_ef_-\\
    &+i\mathbf1_{\Omega_t}\int_t^\infty
       e^{i(t-s)\Delta}P^-\partial_eG(s)\,ds-iC_1^-(t).
 \end{aligned}
\end{equation}
Lemma~\ref{lem:incoming_free_velocity} shows that the first term tends to zero
in $L_x^2(\R^2)$. Combining \eqref{eq:normal_outgoing_identity} and
\eqref{eq:normal_incoming_identity} with $P^++P^-=1$ yields
\[
 \begin{aligned}
 \|\partial_eu(t)\|_{L_x^2(\Omega_t)}
 \leq{}&\|e^{it\Delta}\partial_ef_-\|_{L_x^2(\Omega_t)}
       +C\|G\|_{N^1([t,\infty))}\\
       &+C_{B_H}(1+A+B)(1+\delta^{-1})^6t^{-2}
       \longrightarrow0,
 \end{aligned}
\]
which proves \eqref{eq:normal_exterior_energy_pointwise_decay}.
\end{proof}

\subsection{Exterior gradient bounds and kinetic energy decay}
\label{subsec:gradient_from_exterior_tail}
\label{subsec:pointwise_exterior_kinetic_decay}

We now prove global exterior $L_{t,x}^4$ bounds for $\nabla u$
and then deduce decay of the exterior kinetic energy.
We first bound the gradient on each dyadic time interval and
sum these estimates using the global exterior $L_{t,x}^4$ bound for $u$.
This gives the $N^1$ estimate in
\eqref{eq:normal_exterior_source_assumption}, allowing us to apply
Proposition~\ref{cor:global_exterior_source_normal_rigidity}
in finitely many directions.

For the next proposition, fix $T_0\geq1$, and write
\[
 I_n:=[2^nT_0,2^{n+1}T_0],\qquad
 \mathcal E_{e,c}:=\{(t,x):t\geq T_0,\ x\cdot e>ct\},
 \qquad n\geq0.
\]

\begin{proposition}[Global exterior gradient bounds]
\label{prop:gradient_from_exterior_tail}
Let $u$ be a solution of \eqref{eq:cq} satisfying
\eqref{eq:uniform_H1_bound} and the reduced Duhamel formula
\eqref{eq:reduced_future_duhamel}. Assume that
\begin{equation}\label{eq:gradient_exterior_L4_assumption}
 u\in L_{t,x}^4(\{(t,x):t\geq T_0,\ |x|>ct\})
 \qquad\text{for every }c>0,
\end{equation}
and that, for every fixed smooth radial cutoff $\chi$ as in
\eqref{eq:spatial_cutoff},
\begin{equation}\label{eq:gradient_dyadic_N1_assumption}
 \|\chi(x/T)F(u)\|_{N^1([T,2T])}\longrightarrow0
 \qquad(T\to\infty).
\end{equation}
Then, for every $c>0$ and $e\in\mathbb S^1$,
\begin{equation}\label{eq:global_directional_gradient_L4}
 \partial_eu\in L_{t,x}^4(\mathcal E_{e,c}).
\end{equation}
Consequently, for every $c>0$,
\begin{equation}\label{eq:global_exterior_gradient_L4}
 \nabla u\in L_{t,x}^4(\{(t,x):t\geq T_0,\ |x|>ct\}).
\end{equation}
\end{proposition}
\begin{proof}
We first prove a uniform bound on the dyadic time intervals.
Fix $e\in\mathbb S^1$ and $0<a<b$, and set $A=(a+b)/2$.
Choose a smooth function $\vartheta:\R\to[0,1]$ which is 0
on $(-\infty,a]$ and 1 on $[A,\infty)$, and set
$\zeta(t,x):=\vartheta(x\cdot e/t)$. Then
\begin{equation}\label{eq:gradient_cutoff_supports}
 \begin{gathered}
 \operatorname{supp}\zeta(t)\subset\{x\cdot e\geq at\},\qquad
 \operatorname{supp}(1-\zeta(t))\subset\{x\cdot e\leq At\},\\
 |\partial_e\zeta(t,x)|\leq C_{a,b}/t.
 \end{gathered}
\end{equation}
Take the multipliers in \eqref{eq:gap_crossing_multipliers}
with $B=b$, so that $\delta=b-A$ and $\mu=(A+b)/4$.
Let $C_1^\pm$ be the functions in
Lemma~\ref{lem:normal_crossing_full_frequency} with
$\varphi=1-\zeta$ and $\mathcal G=F(u)$. By
\eqref{eq:normal_nonlinearity_bound} and
\eqref{eq:normal_crossing_global_L4},
\begin{equation}\label{eq:normal_cross_input}
 K_{e,a,b}:=\|C_1^+\|_{L_{t,x}^4([T_0,\infty)\times\R^2)}
              +\|C_1^-\|_{L_{t,x}^4([T_0,\infty)\times\R^2)}
 \leq C_{a,b,B_H}T_0^{-2}.
\end{equation}
These functions vanish outside $\mathcal E_{e,b}$.

Set $f:=\partial_e(\zeta F(u))$. For $T=2^nT_0$ and
$t\in I_n$, \eqref{eq:gradient_cutoff_supports} places the
supports of $\zeta$ and $\partial_e\zeta$ in $|x|\geq aT$.
Choose a fixed smooth radial cutoff $\chi$ which is 1 on
$|x|\geq a$ and 0 on a smaller ball, and write
$\chi_T(x):=\chi(x/T)$. Thus
\begin{equation}\label{eq:gradient_localized_source_bound}
 f=\zeta\,\partial_e(\chi_TF(u))
                +(\partial_e\zeta)\chi_TF(u),\qquad
 \|f\|_{N^0(I_n)}
 \leq C_{a,b}(1+T^{-1})\|\chi_TF(u)\|_{N^1(I_n)}.
\end{equation}
For each fixed $n$, the uniform $H^1$ bound gives
$\|u\|_{L_{t,x}^4(I_n\times\R^2)}\leq C B_H|I_n|^{1/4}<\infty$.
Thus the spacetime estimate \eqref{eq:H1_scattering_criterion}
gives $F(u)\in N^1(I_n)$ for every $n$. Thus
\eqref{eq:gradient_localized_source_bound} makes each $\beta_n$ below finite,
and \eqref{eq:gradient_localized_source_bound} together with
\eqref{eq:gradient_dyadic_N1_assumption} yields
\begin{equation}\label{eq:uniform_gradient_source_blocks}
 \beta_n:=\|f\|_{N^0(I_n)},\qquad
 \sup_{n\geq0}\beta_n<\infty,\qquad
 \beta_n\longrightarrow0.
\end{equation}

We next prove the identity
\begin{equation}\label{eq:normal_gradient_identity}
 \begin{aligned}
 \partial_eu(t)
 ={}&e^{i(t-T_0)\Delta}P^-\partial_eu(T_0)
           +iC_1^+(t)-iC_1^-(t)\\
    &+i\int_t^\infty e^{i(t-s)\Delta}P^+f(s)\,ds
       -i\int_{T_0}^t e^{i(t-s)\Delta}P^-f(s)\,ds
 \end{aligned}
\end{equation}
as a spacetime distribution on $\{t>T_0,\ x\cdot e>bt\}$.
Apply $\partial_eP^+$ and $\partial_eP^-$ to the Duhamel
formula \eqref{eq:appendix_duhamel} with terminal time $T>t$
and initial time $T_0$, respectively. On $x\cdot e>bt$, this gives
\begin{equation}\label{eq:gradient_finite_terminal_duhamel}
 \partial_eP^+u(t)
 =\partial_e e^{i(t-T)\Delta}P^+u(T)
   +i\int_t^T e^{i(t-s)\Delta}P^+f(s)\,ds+iC_{1,T}^+(t),
\end{equation}
and
\begin{equation}\label{eq:gradient_initial_duhamel}
 \partial_eP^-u(t)
 =e^{i(t-T_0)\Delta}P^-\partial_eu(T_0)
  -i\int_{T_0}^t e^{i(t-s)\Delta}P^-f(s)\,ds-iC_1^-(t).
\end{equation}
Here $C_{1,T}^+$ is the finite integral from
Lemma~\ref{lem:normal_crossing_full_frequency}, with
$B=b$, $\varphi=1-\zeta$, and $\mathcal G=F(u)$.
To prove \eqref{eq:normal_gradient_identity},
it suffices to let $T\to\infty$ in
\eqref{eq:gradient_finite_terminal_duhamel} and then add
\eqref{eq:gradient_initial_duhamel}.

Fix a compact interval $J\subset(T_0,\infty)$.
The multiplier argument in the proof of
Lemma~\ref{lem:all_frequency_gap_crossing_S0} applies
for every $\delta>0$ and gives
\[
 \|P^\pm f\|_{N^0(I_n)}\leq C\|f\|_{N^0(I_n)}=C\beta_n.
\]
Thus \eqref{eq:uniform_gradient_source_blocks} allows us
to apply Lemma~\ref{lem:dyadic_offdiagonal_S0} to
$G=P^+f$, which gives
\begin{equation}\label{eq:gradient_projected_integral_limit}
 \int_t^T e^{i(t-s)\Delta}P^+f(s)\,ds
 \longrightarrow
 \int_t^\infty e^{i(t-s)\Delta}P^+f(s)\,ds
 \quad\text{strongly in }L_{t,x}^4(J\times\R^2).
\end{equation}
For $C_{1,T}^+$, let $T'\to\infty$ in
\eqref{eq:crossing_source_tails} and use the limit $C_1^+$
from Lemma~\ref{lem:normal_crossing_full_frequency}.
For $T>\sup J$, this yields
\begin{equation}\label{eq:gradient_crossing_terminal_limit}
 \|C_{1,T}^+-C_1^+\|_{L_{t,x}^4(J\times\R^2)}
 \leq C_{a,b,B_H}|J|^{1/4}T^{-9/4}\longrightarrow0.
\end{equation}

We now show that $\partial_e e^{i(t-T)\Delta}P^+u(T)$
converges weakly to zero in $L_x^2(\R^2)$ as $T\to\infty$
for each fixed $t>T_0$.
For each fixed $t$, the Duhamel
formula \eqref{eq:appendix_duhamel} and
\eqref{eq:reduced_future_duhamel} give
\[
 e^{i(t-T)\Delta}u(T)
 =u(t)-i\int_t^T e^{i(t-s)\Delta}F(u(s))\,ds
 \rightharpoonup0\quad\text{in }H^1(\R^2).
\]
Since $\partial_eP^+:H^1\to L_x^2(\R^2)$ is bounded, the free term in
\eqref{eq:gradient_finite_terminal_duhamel} tends weakly to
zero in $L_x^2(\R^2)$ as $T\to\infty$.
Let $\phi\in C_c^\infty(\{t>T_0,\ x\cdot e>bt\})$, and
choose $J$ to contain its time support. The bound
\[
 \left|\int_{\R^2}
   \partial_e e^{i(t-T)\Delta}P^+u(T)\,
                        \overline{\phi(t,x)}\,dx\right|
 \leq C B_H\|\phi(t)\|_{L_x^2(\R^2)}
\]
has an integrable right-hand side on $J$. Dominated
convergence therefore gives
\begin{equation}\label{eq:gradient_free_term_distribution_limit}
 \lim_{T\to\infty}\int_J\int_{\R^2}
   \partial_e e^{i(t-T)\Delta}P^+u(T)\,
                   \overline{\phi(t,x)}\,dx\,dt=0.
\end{equation}

Multiply \eqref{eq:gradient_finite_terminal_duhamel} by
$\overline\phi$ and integrate over $J\times\R^2$.
The limit of the term containing $u(T)$ is zero by
\eqref{eq:gradient_free_term_distribution_limit}.
For the integral containing $P^+f$ and for $C_{1,T}^+$,
use \eqref{eq:gradient_projected_integral_limit} and
\eqref{eq:gradient_crossing_terminal_limit}, respectively;
H\"older's inequality applies since $\phi\in L_{t,x}^{4/3}(J\times\R^2)$.
Thus, as a spacetime distribution on $\{t>T_0,\ x\cdot e>bt\}$,
\begin{equation}\label{eq:gradient_outgoing_duhamel_identity}
 \partial_eP^+u(t)
 =i\int_t^\infty e^{i(t-s)\Delta}P^+f(s)\,ds+iC_1^+(t).
\end{equation}
Adding \eqref{eq:gradient_initial_duhamel} to \eqref{eq:gradient_outgoing_duhamel_identity}
and using $P^++P^-=I$ proves
\eqref{eq:normal_gradient_identity}.

For each $j\geq0$, set
\[
 \begin{aligned}
 q_j&:=\|\partial_eu\|_{L_{t,x}^4(\{(t,x):t\in I_j,\ x\cdot e>bt\})},\\
 h_j&:=\|e^{i(t-T_0)\Delta}P^-\partial_eu(T_0)\|_{L_{t,x}^4(I_j\times\R^2)}
       +\|C_1^+\|_{L_{t,x}^4(I_j\times\R^2)}
       +\|C_1^-\|_{L_{t,x}^4(I_j\times\R^2)}.
 \end{aligned}
\]
The Strichartz estimate \eqref{eq:appendix_strichartz} and
\eqref{eq:normal_cross_input} give
$\|h\|_{\ell^4}\leq C B_H+K_{e,a,b}$. By
\eqref{eq:dyadic_duhamel_output}, the identity
\eqref{eq:normal_gradient_identity} implies
\begin{equation}\label{eq:gradient_dyadic_output}
 q_j\leq h_j+C(\mathcal T\beta)_j
       =h_j+C\sum_{k\geq0}2^{-|j-k|/4}\beta_k.
\end{equation}
Thus \eqref{eq:dyadic_sequence_operator}, with $p=\infty$, gives
\begin{equation}\label{eq:uniform_directional_gradient_blocks}
 M_{e,b}:=\sup_{j\geq0}q_j
 \leq\|h\|_{\ell^4}+C\|\beta\|_{\ell^\infty}<\infty.
\end{equation}
Repeating this argument with $a,b$ replaced by $a/2,a$
gives $M_{e,a}<\infty$.

We now use \eqref{eq:gradient_exterior_L4_assumption} to sum
\eqref{eq:gradient_dyadic_output}. Since
$\mathcal E_{e,a}\subset\{t\geq T_0,\ |x|>at\}$, the quantities
$U_p:=\|u\|_{L_{t,x}^p(\mathcal E_{e,a})}$ satisfy
\begin{equation}\label{eq:gradient_halfspace_spacetime_input}
 U_4<\infty,\qquad U_8^8\leq C B_H^4U_4^4<\infty,
\end{equation}
where the second inequality follows from the halfspace $L_{t,x}^8$
estimate \eqref{eq:halfspace_L8_from_L4}.
Set
\[
 \begin{aligned}
 A_n&:=\|u\|_{L_{t,x}^4(\{(t,x):t\in I_n,\ x\cdot e>at\})}^2
          +\|u\|_{L_{t,x}^8(\{(t,x):t\in I_n,\ x\cdot e>at\})}^4,\\
 g_n&:=\|\partial_eu\|_{L_{t,x}^4(\{(t,x):t\in I_n,\ x\cdot e>at\})}.
 \end{aligned}
\]
The intervals $I_n$ are disjoint up to their endpoints. Hence
\begin{equation}\label{eq:gradient_coefficient_sequence}
 \begin{aligned}
 \|A\|_{\ell^2}
 &\leq\left(\sum_n\|u\|_{L_{t,x}^4(\{(t,x):t\in I_n,\ x\cdot e>at\})}^4\right)^{1/2}
       +\left(\sum_n\|u\|_{L_{t,x}^8(\{(t,x):t\in I_n,\ x\cdot e>at\})}^8\right)^{1/2}\\
 &=U_4^2+U_8^4,\qquad \|g\|_{\ell^\infty}=M_{e,a}.
 \end{aligned}
\end{equation}
The pointwise inequality
$|\partial_eF(u)|\leq C(|u|^2+|u|^4)|\partial_eu|$
and H\"older's inequality give
\[
 \|\zeta\partial_eF(u)\|_{N^0(I_n)}\leq C A_ng_n.
\]
For the term containing $\partial_e\zeta$, put
$r_n:=\|(\partial_e\zeta)F(u)\|_{N^0(I_n)}$. By
\eqref{eq:gradient_cutoff_supports} and H\"older's inequality,
\begin{equation}\label{eq:gradient_cutoff_derivative_sum}
 \begin{aligned}
 \|r\|_{\ell^{4/3}}
 &=\|(\partial_e\zeta)F(u)\|_{N^0([T_0,\infty))}\\
 &\leq C_{a,b}T_0^{-1}(U_4^3+U_8^4U_4)<\infty.
 \end{aligned}
\end{equation}
Since $f=\zeta\partial_eF(u)+(\partial_e\zeta)F(u)$,
we have $\beta_n\leq CA_ng_n+r_n$. Combining
\eqref{eq:gradient_coefficient_sequence} and
\eqref{eq:gradient_cutoff_derivative_sum} with
$\ell^2\subset\ell^4$ and $\ell^{4/3}\subset\ell^4$ yields
\begin{equation}\label{eq:gradient_source_sequence_l4}
 \begin{aligned}
 \|\beta\|_{\ell^4}
 &\leq C\|A g\|_{\ell^2}+\|r\|_{\ell^{4/3}}\\
 &\leq CM_{e,a}(U_4^2+U_8^4)
        +C_{a,b}T_0^{-1}(U_4^3+U_8^4U_4)<\infty.
 \end{aligned}
\end{equation}
For every integer $N\geq0$, \eqref{eq:gradient_dyadic_output}
and \eqref{eq:dyadic_sequence_operator}, now with $p=4$, give
\begin{equation}\label{eq:gradient_partial_sum_bound}
 \left(\sum_{j=0}^N q_j^4\right)^{1/4}
 \leq\|h\|_{\ell^4}+C\|\beta\|_{\ell^4}.
\end{equation}
The right-hand side is finite and independent of $N$.
Let $N\to\infty$ in \eqref{eq:gradient_partial_sum_bound}.
Monotone convergence, \eqref{eq:gradient_source_sequence_l4},
and $\|h\|_{\ell^4}\leq CB_H+K_{e,a,b}$ give
\begin{equation}\label{eq:gradient_quantitative_closure}
 \begin{aligned}
 \|\partial_eu\|_{L_{t,x}^4(\mathcal E_{e,b})}
 &=\left(\sum_{j\geq0}q_j^4\right)^{1/4}\\
 &\leq C_{a,b}\bigl[B_H+K_{e,a,b}
       +M_{e,a}(U_4^2+U_8^4)
       +T_0^{-1}(U_4^3+U_8^4U_4)\bigr]<\infty.
 \end{aligned}
\end{equation}
This proves \eqref{eq:global_directional_gradient_L4}.

To estimate the full gradient, fix $c>0$. For $k=0,\ldots,7$,
let $\Sigma_k$ be the sector centered at angle $k\pi/4$ with
half-angle $\pi/8$, that is,
\[
 \Sigma_k:=\left\{r(\cos\theta,\sin\theta):r>0,\quad
            \frac{k\pi}{4}-\frac\pi8\leq\theta
                         \leq\frac{k\pi}{4}+\frac\pi8\right\}.
\]
Set
\[
 e_k^\pm:=\bigl(\cos(k\pi/4\pm\pi/4),
                    \sin(k\pi/4\pm\pi/4)\bigr).
\]
The sectors cover $\R^2\setminus\{0\}$, and $e_k^+,e_k^-$ are
orthonormal. For $x\in\Sigma_k$,
$x\cdot e_k^\pm\geq|x|\cos(3\pi/8)$, so
\begin{equation}\label{eq:gradient_sector_cover}
 \Sigma_k\cap\{|x|>ct\}
 \subset\{x:x\cdot e_k^+>(c/4)t,\ x\cdot e_k^->(c/4)t\},
 \qquad t>0.
\end{equation}
On each sector,
\[
 |\nabla u|^2=|\partial_{e_k^+}u|^2+|\partial_{e_k^-}u|^2,
 \qquad
 |\nabla u|^4\leq2\bigl(|\partial_{e_k^+}u|^4
                            +|\partial_{e_k^-}u|^4\bigr).
\]
Since the sectors cover $\{|x|>ct\}$, we integrate this
pointwise bound over $\Sigma_k\cap\{|x|>ct\}$ and
$t\in[T_0,\infty)$, and sum over $k$. The inclusion
\eqref{eq:gradient_sector_cover} then gives
\[
 \begin{aligned}
 &\int_{T_0}^\infty\int_{|x|>ct}|\nabla u|^4\,dx\,dt\\
 &\quad\leq\sum_{k=0}^7\int_{T_0}^\infty
            \int_{\Sigma_k\cap\{|x|>ct\}}|\nabla u|^4\,dx\,dt\\
 &\quad\leq2\sum_{k=0}^7\sum_{\pm}
      \|\partial_{e_k^\pm}u\|_{L_{t,x}^4(\mathcal E_{e_k^\pm,c/4})}^4
 <\infty.
 \end{aligned}
\]
Applying \eqref{eq:global_directional_gradient_L4}
with $e=e_k^\pm$ and $c$ replaced by $c/4$ shows that
each term in the last sum is finite.
This proves \eqref{eq:global_exterior_gradient_L4}.
\end{proof}

\begin{corollary}[Global $N^1$ bounds for the localized nonlinearity]
\label{cor:global_exterior_N1_from_gradient}
Under the hypotheses of
Proposition~\ref{prop:gradient_from_exterior_tail}, fix $c>0$.
Let $\zeta(t,x)$ be a smooth cutoff satisfying
\begin{equation}\label{eq:global_nonlinearity_cutoff_conditions}
 \operatorname{supp}\zeta(t)\subset\{|x|>ct\},\qquad
 \|\zeta\|_{L_{t,x}^\infty([T_0,\infty)\times\R^2)}<\infty,\qquad
 |\nabla_x\zeta(t,x)|\leq C/t\quad(t\geq T_0).
\end{equation}
Then
\begin{equation}\label{eq:global_exterior_N1}
 \zeta F(u)\in N^1([T_0,\infty)).
\end{equation}
\end{corollary}
\begin{proof}
Write $\Omega:=\{(t,x):t\geq T_0,\ |x|>ct\}$ and set
\[
 U_4:=\|u\|_{L_{t,x}^4(\Omega)},\qquad
 U_8:=\|u\|_{L_{t,x}^8(\Omega)},\qquad
 D_4:=\|\nabla u\|_{L_{t,x}^4(\Omega)}.
\]
The first and third norms are finite by
\eqref{eq:gradient_exterior_L4_assumption} and
\eqref{eq:global_exterior_gradient_L4}. To bound $U_8$,
use the four halfspaces in \eqref{eq:radial_halfspace_cover}
with $a_0=c/2$. The halfspace estimate
\eqref{eq:halfspace_L8_from_L4} gives
\[
 \begin{aligned}
 U_8^8&\leq CB_H^4\sum_{k=1}^4
                  \|u\|_{L_{t,x}^4(\mathcal E_{e_k,c/2})}^4\\
 &\leq4CB_H^4
       \int_{T_0}^\infty\int_{|x|>(c/2)t}|u|^4\,dx\,dt<\infty,
 \end{aligned}
\]
where finiteness follows from
\eqref{eq:gradient_exterior_L4_assumption} with $c$ replaced by $c/2$.
The product rule gives, on $\Omega$,
\[
 |\zeta F(u)|+|\nabla(\zeta F(u))|
 \leq C_\zeta(|u|^2+|u|^4)
                   \bigl(|\nabla u|+(1+T_0^{-1})|u|\bigr).
\]
Outside $\Omega$ the left-hand side vanishes. H\"older's
inequality yields
\[
 \|\zeta F(u)\|_{N^1([T_0,\infty))}
 \leq C_\zeta(U_4^2+U_8^4)
                  \bigl(D_4+(1+T_0^{-1})U_4\bigr)<\infty.
\]
\end{proof}

\begin{remark}
The uniform $H^1$ bound also gives, for every $0<\sigma<1$,
\begin{equation}\label{eq:F_uniform_fractional}
 \sup_{t\in\R}\|F(u(t))\|_{H^\sigma}\leq C_{B_H,\sigma}.
\end{equation}
Indeed, the fractional product estimate
\eqref{eq:appendix_fractional_product} gives
\[
 \|F(u(t))\|_{H^\sigma}
 \leq C_\sigma\left(\|u(t)\|_{L_x^{4/(1-\sigma)}(\R^2)}^2
               +\|u(t)\|_{L_x^{8/(1-\sigma)}(\R^2)}^4\right)
               \|u(t)\|_{W^{\sigma,2/\sigma}}
 \leq C_{B_H,\sigma},
\]
where $\sigma/2+(1-\sigma)/2=1/2$.
\end{remark}

Finally, we apply
Proposition~\ref{cor:global_exterior_source_normal_rigidity}
in finitely many directions to obtain decay of the exterior kinetic energy.

\begin{corollary}[Decay of the exterior kinetic energy]
\label{cor:full_exterior_gradient_pointwise}
Let $u$ be a uniformly $H^1$-bounded solution of \eqref{eq:cq}
satisfying \eqref{eq:reduced_future_duhamel}, and fix $T_0>0$.
Suppose that, for every $c>0$ and every smooth cutoff $\zeta$ with
$\operatorname{supp}\zeta(t)\subset\{|x|>ct\}$,
$\|\zeta\|_{L_{t,x}^\infty([T_0,\infty)\times\R^2)}<\infty$, and
$|\nabla_x\zeta(t,x)|\leq C/t$ for $t\geq T_0$, one has
\begin{equation}\label{eq:pointwise_exterior_source_assumption}
 \zeta F(u)\in N^1([T_0,\infty)).
\end{equation}
Then, for each fixed $c>0$,
\begin{equation}\label{eq:pointwise_exterior_kinetic_vanishes}
 \|\nabla u(t)\|_{L_x^2(\{x:|x|>ct\})}\longrightarrow0
 \qquad(t\to\infty).
\end{equation}
\end{corollary}
\begin{proof}
Fix $c>0$, set $b=c/4$, and choose $0<a<A<b$.
Let $\vartheta\in C^\infty(\R;[0,1])$ be 0 on
$(-\infty,a]$ and 1 on $[A,\infty)$. For each fixed
$e\in\mathbb S^1$, set $\zeta_e(t,x):=\vartheta(x\cdot e/t)$.
Then
\[
 \operatorname{supp}\zeta_e(t)
 \subset\{x:x\cdot e\geq at\}
 \subset\{|x|>(a/2)t\},\qquad
 |\nabla_x\zeta_e(t,x)|\leq\|\vartheta'\|_{L_v^\infty(\R)}/t.
\]
Thus \eqref{eq:pointwise_exterior_source_assumption}, with speed
$a/2$, gives $\zeta_eF(u)\in N^1([T_0,\infty))$.
Proposition~\ref{cor:global_exterior_source_normal_rigidity},
with $B=b$, gives
\eqref{eq:normal_exterior_energy_pointwise_decay}
for every fixed $e$.

Use the sectors $\Sigma_k$ and orthonormal pairs $e_k^\pm$
constructed in the proof of
Proposition~\ref{prop:gradient_from_exterior_tail}.
Since $|\nabla u|^2=|\partial_{e_k^+}u|^2+|\partial_{e_k^-}u|^2$,
the inclusion \eqref{eq:gradient_sector_cover} gives
\[
 \begin{aligned}
 \|\nabla u(t)\|_{L_x^2(\{x:|x|>ct\})}^2
 &\leq\sum_{k=0}^7
       \int_{\Sigma_k\cap\{|x|>ct\}}|\nabla u(t,x)|^2\,dx\\
 &\leq\sum_{k=0}^7\sum_{\pm}
       \|\partial_{e_k^\pm}u(t)\|_{L_x^2(\{x:x\cdot e_k^\pm>bt\})}^2.
 \end{aligned}
\]
Applying \eqref{eq:normal_exterior_energy_pointwise_decay}
with $e=e_k^\pm$ and $B=b$ gives
$\|\partial_{e_k^\pm}u(t)\|_{L_x^2(\{x:x\cdot e_k^\pm>bt\})}^2
\longrightarrow0$ as $t\to\infty$.
The finite sum therefore tends to zero, which proves
\eqref{eq:pointwise_exterior_kinetic_vanishes}.
\end{proof}

\section{Proof of Theorem~\ref{thm:threshold_scattering}}
\label{subsec:normal_cone_rigidity_closure}

The preceding sections establish several properties of a
minimal-mass nonscattering solution. We now rule out such a solution
by contradiction and prove
Theorem~\ref{thm:threshold_scattering}.

Suppose that Theorem~\ref{thm:threshold_scattering}
is false. Then there exists a minimal-mass nonscattering solution
$u$ satisfying \eqref{eq:minimal_mass_nonscattering}.
Conservation of mass and energy, \eqref{eq:coercive_energy}, and
the continuation criterion in Lemma~\ref{lem:local_small_data} give
\[
 u\in C(\R;H^1(\R^2)),\qquad
 B_H=\sup_{t\in\R}\|u(t)\|_{H^1}<\infty,
\]
where $B_H$ is the quantity defined in \eqref{eq:uniform_H1_bound}.
Here $u$ denotes the solution after the Galilean transform in
Section~\ref{subsec:conservation}, which sets $P(u)=0$ and
preserves the mass and $L_{t,x}^4$ norm.
The condition $P(u)=0$ is used in
\eqref{eq:truncated_center_motion} and
\eqref{eq:translated_action_bound}.
We use this transformed solution throughout the contradiction
argument, with conserved energy $E=E(u)$.
Proposition~\ref{prop:reduced_future_duhamel} supplies the reduced
Duhamel formula \eqref{eq:reduced_future_duhamel}, with convergence
weakly in $H^1$.

Fix $T_0=1$. The exterior $L_{t,x}^4$ limit
\eqref{eq:cubic_window_unconditional} and the $L_{t,x}^8$ estimate
\eqref{eq:dyadic_halfspace_L8}
verify \eqref{eq:local_exterior_L4_L8_assumption}.
For each fixed $c>0$,
Proposition~\ref{prop:conditional_exterior_global_L4}
gives $T_1=T_1(c)\geq T_0$ such that
\[
 \int_{T_1}^{\infty}\int_{|x|>ct}|u|^4\,dx\,dt<\infty.
\]
By the uniform $H^1$ bound and the Sobolev embedding
\eqref{eq:appendix_H1_embeddings},
\[
 \begin{aligned}
 \int_{T_0}^{\infty}\int_{|x|>ct}|u|^4\,dx\,dt
 &\leq\int_{T_0}^{T_1}\|u(t)\|_{L_x^4(\R^2)}^4\,dt
       +\int_{T_1}^{\infty}\int_{|x|>ct}|u|^4\,dx\,dt\\
 &\leq CB_H^4(T_1-T_0)
       +\int_{T_1}^{\infty}\int_{|x|>ct}|u|^4\,dx\,dt<\infty.
 \end{aligned}
\]
Thus \eqref{eq:gradient_exterior_L4_assumption} holds with
this fixed $T_0$ for every $c>0$.
Moreover, \eqref{eq:actual_exterior_full_N1} gives
\[
 \|\chi(x/T)F(u)\|_{N^1([T,2T])}\longrightarrow0
 \qquad(T\to\infty)
\]
for every fixed smooth radial cutoff $\chi$ satisfying
\eqref{eq:spatial_cutoff}.
These are the hypotheses
\eqref{eq:gradient_exterior_L4_assumption} and
\eqref{eq:gradient_dyadic_N1_assumption} of
Proposition~\ref{prop:gradient_from_exterior_tail}. Hence
\[
 \|\nabla u\|_{L_{t,x}^4(\{(t,x):t\geq T_0,\ |x|>ct\})}<\infty
 \qquad\text{for every }c>0.
\]
Corollary~\ref{cor:global_exterior_N1_from_gradient}
then gives
\[
 \zeta F(u)\in N^1([T_0,\infty))
\]
for every $c>0$ and every smooth cutoff $\zeta$ satisfying
\eqref{eq:global_nonlinearity_cutoff_conditions}.

Fix $c>0$. Corollary~\ref{cor:full_exterior_gradient_pointwise}
and \eqref{eq:uniform_H1_bound} give
\[
 0\leq\|\nabla u(t)\|_{L_x^2(\{x:|x|>ct\})}^2\leq B_H^2,
 \qquad
 \|\nabla u(t)\|_{L_x^2(\{x:|x|>ct\})}^2\longrightarrow0
 \quad(t\to\infty).
\]
Given $\varepsilon>0$, choose $T_\varepsilon\geq T_0$ such that
$\|\nabla u(t)\|_{L_x^2(\{x:|x|>ct\})}^2\leq\varepsilon$ for
$t\geq T_\varepsilon$. For $R\geq T_\varepsilon$, splitting
the time integral at $T_\varepsilon$ gives
\[
 0\leq\frac1R\int_{T_0}^R\int_{|x|>ct}
       |\nabla u(t,x)|^2\,dx\,dt
 \leq\frac{(T_\varepsilon-T_0)B_H^2}{R}+\varepsilon.
\]
Letting $R\to\infty$ and then $\varepsilon\downarrow0$, we obtain
\begin{equation}\label{eq:exterior_kinetic_average_zero}
 \frac1R\int_{T_0}^R\int_{|x|>ct}
       |\nabla u(t,x)|^2\,dx\,dt
 \longrightarrow0
 \qquad(R\to\infty).
\end{equation}

For $T_0\leq t\leq R$, we have
$\{|x|\leq ct\}\subset\{|x|\leq cR\}$.
Thus \eqref{eq:ballistic_evacuation}, with $z=0$ and $v=c$, gives
\begin{equation}\label{eq:interior_kinetic_average_zero}
 0\leq\frac1R\int_{T_0}^R\int_{|x|\leq ct}
              |\nabla u(t,x)|^2\,dx\,dt
 \leq\frac1R\int_0^R\int_{|x|\leq cR}
              |\nabla u(t,x)|^2\,dx\,dt
 \longrightarrow0.
\end{equation}
Also note that
\begin{equation}\label{eq:initial_kinetic_average_zero}
 0\leq\frac1R\int_0^{T_0}\|\nabla u(t)\|_{L_x^2(\R^2)}^2\,dt
 \leq\frac{T_0B_H^2}{R}\longrightarrow0.
\end{equation}
Combining \eqref{eq:initial_kinetic_average_zero} with \eqref{eq:totalavg} and
\eqref{eq:interior_kinetic_average_zero}, we find
\[
 \begin{aligned}
 \frac1R\int_{T_0}^R\int_{|x|>ct}
       |\nabla u(t,x)|^2\,dx\,dt
 &=\frac1R\int_0^R\|\nabla u(t)\|_{L_x^2(\R^2)}^2\,dt
   -\frac1R\int_0^{T_0}\|\nabla u(t)\|_{L_x^2(\R^2)}^2\,dt\\
 &\quad-\frac1R\int_{T_0}^R\int_{|x|\leq ct}
       |\nabla u(t,x)|^2\,dx\,dt\\
 &\longrightarrow2E>0\qquad(R\to\infty),
 \end{aligned}
\]
contradicting \eqref{eq:exterior_kinetic_average_zero}.

Now let $u$ be any $H^1$ solution with mass $m_Q$.
If its forward $L_{t,x}^4$ norm were infinite, the Galilean
transform in Section~\ref{subsec:conservation} would preserve
this norm and set its momentum to zero, producing a solution
satisfying \eqref{eq:minimal_mass_nonscattering}.
The contradiction above therefore gives
\[
 \|u\|_{L_{t,x}^4([0,\infty)\times\R^2)}<\infty.
\]
Time reversal combined with complex conjugation gives another
solution $\overline{u(-t,x)}$ with the same mass. Applying the
forward bound to this solution yields
\[
 \begin{aligned}
 \|u\|_{L_{t,x}^4(\R\times\R^2)}^4
 &=\|u\|_{L_{t,x}^4((-\infty,0]\times\R^2)}^4
   +\|u\|_{L_{t,x}^4([0,\infty)\times\R^2)}^4
 <\infty.
 \end{aligned}
\]
The scattering criterion in
Lemma~\ref{lem:H1_scattering_criterion} now gives
$u_\pm\in H^1(\R^2)$ such that
\[
 \|u(t)-e^{it\Delta}u_\pm\|_{H^1}\longrightarrow0
 \qquad(t\to\pm\infty).
\]

\appendix
\section{Analytic tools and auxiliary results}
\label{app:local_theory}

\subsection{Linear and functional estimates}

\begin{lemma}[Strichartz estimates {\cite{Strichartz,KT}}]
\label{lem:strichartz_estimates}
Let $I\subset\R$ be an interval, let $t_0\in I$, and let
$0\leq s\leq1$. Fix two Schr\"odinger-admissible pairs
$(q,r)$ and $(\widetilde q,\widetilde r)$ in dimension two.
In particular, neither pair is $(2,\infty)$.
For $f\in H^s(\R^2)$ and
$G\in L_t^{\widetilde q'}(I;W_x^{s,\widetilde r'}(\R^2))$, we have
\begin{equation}\label{eq:appendix_strichartz}
\begin{aligned}
 \|e^{i(t-t_0)\Delta}f\|_{L_t^q(I;W_x^{s,r}(\R^2))}
 &\leq C\|f\|_{H^s},\\
 \left\|\int_{t_0}^t e^{i(t-\tau)\Delta}G(\tau)\,d\tau
       \right\|_{L_t^q(I;W_x^{s,r}(\R^2))}
 &\leq C\|G\|_{L_t^{\widetilde q'}(I;W_x^{s,\widetilde r'}(\R^2))}.
\end{aligned}
\end{equation}
The corresponding dual estimate is
\begin{equation}\label{eq:appendix_dual_strichartz}
 \left\|\int_I e^{-it\Delta}G(t)\,dt\right\|_{H^s}
 \leq C\|G\|_{L_t^{\widetilde q'}(I;W_x^{s,\widetilde r'}(\R^2))}.
\end{equation}
The constants depend only on $s$ and the admissible pairs.
\end{lemma}

\begin{lemma}[Fractional product estimates {\cite{KatoPonce,GrafakosOh}}]
\label{lem:appendix_fractional_product}
Let $0\leq s\leq1$ and $1<r,p_1,p_2,q_1,q_2<\infty$ satisfy
$1/r=1/p_1+1/q_1=1/p_2+1/q_2$. Then
\begin{equation}\label{eq:appendix_fractional_product}
 \|D^s(fg)\|_{L_x^r(\R^2)}
 \leq C\bigl(\|D^sf\|_{L_x^{p_1}(\R^2)}\|g\|_{L_x^{q_1}(\R^2)}
              +\|f\|_{L_x^{p_2}(\R^2)}\|D^sg\|_{L_x^{q_2}(\R^2)}\bigr).
\end{equation}
The estimate \eqref{eq:appendix_fractional_product} also holds with $J^s$ in place of $D^s$.
\end{lemma}

\begin{lemma}[Sobolev embedding and interpolation]
\label{lem:appendix_sobolev_interpolation}
Let $0<s<2$, $1<p<2/s$, and $1/q=1/p-s/2$. Then
\begin{equation}\label{eq:appendix_sobolev}
 \|f\|_{L_x^q(\R^2)}\leq C_{s,p}\|D^sf\|_{L_x^p(\R^2)},
 \qquad f\in\mathcal S(\R^2).
\end{equation}
We also have
\begin{equation}\label{eq:appendix_H1_embeddings}
 H^1(\R^2)\hookrightarrow L_x^p(\R^2)\quad(2\leq p<\infty),
 \qquad
 H^1(\R^2)\hookrightarrow W^{\sigma,2/\sigma}(\R^2)
 \quad(0<\sigma<1).
\end{equation}
For $0\leq\theta\leq1$ and $1<p<\infty$, complex interpolation
\cite{Calderon1964} gives
\begin{equation}\label{eq:appendix_bessel_interpolation}
 \|J^\theta f\|_{L_x^p(\R^2)}
 \leq C_{p,\theta}\|f\|_{L_x^p(\R^2)}^{1-\theta}\|Jf\|_{L_x^p(\R^2)}^\theta,
 \qquad f\in W^{1,p}(\R^2).
\end{equation}
For any interval $I$ and $1\leq a\leq\infty$,
\[
 \|G\|_{L_t^a(I;W_x^{\theta,p}(\R^2))}
 \leq C_{p,\theta}
       \|G\|_{L_t^a(I;L_x^p(\R^2))}^{1-\theta}
       \|G\|_{L_t^a(I;W_x^{1,p}(\R^2))}^\theta.
\]
The constants are independent of $I$.
\end{lemma}

\begin{lemma}[Gagliardo--Nirenberg inequalities
{\cite{Nirenberg1959}}]
\label{lem:appendix_gagliardo_nirenberg}
\label{lem:exterior_global_L8}
Let $\Omega=\R^2$ or let $\Omega$ be an open halfspace in $\R^2$.
If $2\leq q\leq p<\infty$ and $f\in H^1(\Omega)\cap L_x^q(\Omega)$,
then
\begin{equation}\label{eq:appendix_gagliardo_nirenberg}
 \|f\|_{L_x^p(\Omega)}
 \leq C_{p,q}\|\nabla f\|_{L_x^2(\Omega)}^{1-q/p}
                    \|f\|_{L_x^q(\Omega)}^{q/p}.
\end{equation}
 In particular,
\begin{equation}\label{eq:appendix_GN_L8}
 \|f\|_{L_x^8(\Omega)}^8
 \leq C\|\nabla f\|_{L_x^2(\Omega)}^4
        \|f\|_{L_x^4(\Omega)}^4.
\end{equation}
Consequently, if $u\in L_t^\infty([T_0,\infty);H_x^1(\R^2))$ with
$\sup_{t\geq T_0}\|u(t)\|_{H^1}\leq B_H$, then for $|e|=1$ and $c>0$,
\begin{equation}\label{eq:halfspace_L8_from_L4}
 \int_{T_0}^\infty\int_{x\cdot e>ct}|u|^8
 \leq C B_H^4\int_{T_0}^\infty\int_{x\cdot e>ct}|u|^4.
\end{equation}
\end{lemma}

\subsection{Local theory, scattering, and stability}

\begin{lemma}[$H^1$ local theory and small data
{\cite{GinibreVelo,CazenaveWeisslerH1}}]
\label{lem:local_small_data}
Let $t_0\in\R$ and $f\in H^1(\R^2)$. There is a unique
maximal-lifespan solution
$u\in C(I_{\max};H^1(\R^2))$ of \eqref{eq:cq}, with
$t_0\in I_{\max}$ and $u(t_0)=f$, in the class
$S^1_{\mathrm{loc}}(I_{\max})$. This solution satisfies \eqref{eq:appendix_duhamel}.
The following properties hold:
\begin{enumerate}
\item \textbf{Conservation laws.}
The mass, energy, and momentum of $u$ are conserved.

\item \textbf{Spacetime bounds and continuous dependence.}
For every compact $J\subset I_{\max}$, we have $\|u\|_{S^1(J)}<\infty$.
The solution depends continuously on the initial data in
$C(J;H^1)\cap S^1(J)$ on each such interval containing $t_0$.

\item \textbf{Continuation criterion.}
If a finite endpoint of $I_{\max}$ is approached while the $H^1$
norm remains bounded, the solution extends beyond that endpoint.

\item \textbf{Persistence of regularity and smooth approximation.}
If $f\in H^k$ for an integer $k\geq1$, then
$u\in C(I_{\max};H^k)$. In particular, if
$f_m\in C_c^\infty(\R^2)$ and $f_m\to f$ in $H^1$, the corresponding
smooth solutions exist on every compact $J\subset I_{\max}$ containing
$t_0$ for all sufficiently large $m$, and converge to $u$ in
$C(J;H^1)\cap S^1(J)$.

\item \textbf{Small-data theory.}
There are constants $\delta_0>0$ and $C>0$ such that, if
$\|f\|_{H^1}\leq\delta_0$, then $I_{\max}=\R$ and
\begin{equation}\label{eq:appendix_small_data}
 \|u\|_{S^1(\R)}\leq C\|f\|_{H^1}.
\end{equation}
Moreover,
\begin{equation}\label{eq:appendix_small_data_correction}
 \|u-e^{i(t-t_0)\Delta}f\|_{S^1(\R)}
 \leq C\|f\|_{H^1}^3.
\end{equation}
\end{enumerate}
\end{lemma}

The next lemma is a consequence of the local and scattering theory
in \cite{ChengMiaoZhao,Cheng}.

\begin{lemma}[Spacetime criterion for $H^1$ scattering {\cite{ChengMiaoZhao,Cheng}}]
\label{lem:H1_scattering_criterion}
Let $I$ be an interval containing $t_0$, and let
$u\in C(I;H^1)$ solve~\eqref{eq:cq} in the sense of the local
$H^1$ theory. Suppose that
\[
 \|u\|_{L_t^\infty(I;H_x^1(\R^2))}\leq B <\infty,
 \qquad \|u\|_{L_{t,x}^4(I\times\R^2)}\leq L<\infty.
\]
Then
\begin{equation}\label{eq:H1_scattering_criterion}
 \|u\|_{S^1(I)}+
 \|F(u)\|_{N^1(I)}\leq C(B,L),
\end{equation}
where $C(B,L)$ is independent of $I$. If $I$ contains
$[t_0,\infty)$ (respectively, $(-\infty,t_0]$), then $u$ scatters
forward (respectively, backward) in time in $H^1$.
\end{lemma}

\begin{lemma}[Regularity bounds for the cubic equation {\cite{Cazenave}}]
\label{lem:persistence_regularity}
Let $I\subset\R$ be an interval and let
$w\in C(I;L_x^2(\R^2))$ solve \eqref{eq:cubic_nls} in the sense of
the local $L_x^2(\R^2)$ theory. Assume that
$\|w\|_{L_{t,x}^4(I\times\R^2)}\leq L<\infty$.
Fix $0\leq s\leq1$ and $t_0\in I$. If $w(t_0)\in H^s$, then
$w\in C(I;H^s)$ and
\begin{equation}\label{eq:appendix_cubic_regularity}
 \|D^sw\|_{S^0(I)}
 \leq C_s(L)\|D^sw(t_0)\|_{L_x^2(\R^2)}.
\end{equation}
Here $D^0$ is the identity, and the constant is independent of $I$
and $t_0$.

If $\pm\infty$ is an endpoint of $I$ and $w$ has an $L_x^2(\R^2)$
scattering state $f_\pm\in H^s$ , the same conclusion holds with $\|D^sf_\pm\|_{L_x^2(\R^2)}$ on the right-hand
side.
In this case,
\[
 \|w(t)-e^{it\Delta}f_\pm\|_{H^s}\longrightarrow0
 \qquad(t\to\pm\infty).
\]
\end{lemma}

\begin{lemma}[Stability {\cite{CKSTT,TaoVisan,Cheng}}]
\label{lem:stability}
Let $I\subset\R$ be a compact time interval, let $t_0\in I$, and let
$\widetilde v\in C(I;H^1(\R^2))$ solve the approximate equation
\[
 (i\partial_t+\Delta)\widetilde v=F(\widetilde v)+e
\]
on $I\times\R^2$ in the Duhamel sense. Assume that
\[
 \|\widetilde v\|_{L_t^\infty(I;H_x^1(\R^2))}\leq B,
 \qquad
 \|\widetilde v\|_{L_t^4(I;W_x^{1/2,4}(\R^2))}\leq L
\]
for some positive constants $B$ and $L$. Let $f\in H^1(\R^2)$ satisfy
\[
 \|f-\widetilde v(t_0)\|_{H^{1/2}}\leq B'
\]
for some positive constant $B'$. Assume further that
\[
 \bigl\|e^{i(t-t_0)\Delta}(f-\widetilde v(t_0))\bigr\|_{L_t^4(I;W_x^{1/2,4}(\R^2))}
 \leq\varepsilon,
 \qquad
 \|\langle D\rangle^{1/2}e\|_{L_{t,x}^{4/3}(I\times\R^2)}
 \leq\varepsilon
\]
for some $0<\varepsilon<\varepsilon_0(B,B',L)$, where
$\varepsilon_0(B,B',L)>0$ is sufficiently small.
Then there exists a unique solution $v\in C(I;H^1(\R^2))$ of
\eqref{eq:cq} with $v(t_0)=f$ such that
\[
 \|v-\widetilde v\|_{L_t^4(I;W_x^{1/2,4}(\R^2))}
 \leq C(B,B',L)\varepsilon,
 \qquad
 \|v\|_{L_t^4(I;W_x^{1/2,4}(\R^2))}\leq L+C(B,B',L)\varepsilon.
\]
The constants $\varepsilon_0$ and $C$ are independent of $I$ and $t_0$.
\end{lemma}

The same conclusion holds on an unbounded interval when the hypotheses
hold there with the same bounds. Apply the lemma on an increasing
sequence of compact intervals containing $t_0$. Uniqueness identifies
the solutions on overlaps, and monotone convergence gives the stated
spacetime estimates.

\subsection{Scattering and classification for the cubic equation}

\begin{lemma}[Scattering below the ground-state mass
{\cite[Theorem~1.6]{DodsonSub}}]
\label{lem:cubic_subthreshold}
Let $\phi\in L_x^2(\R^2)$ satisfy $M(\phi)<m_Q$.
Then the solution $w$ of \eqref{eq:cubic_nls} with $w(0)=\phi$
is global and satisfies
\[
 \|w\|_{L_{t,x}^4(\R\times\R^2)}<\infty.
\]
Moreover, there exist $w_\pm\in L_x^2(\R^2)$ such that
\[
 \|w(t)-e^{it\Delta}w_\pm\|_{L_x^2(\R^2)}\longrightarrow0
 \qquad(t\to\pm\infty).
\]
\end{lemma}

\begin{lemma}[Classification at the ground-state mass
{\cite[Theorem~3]{Dodson}}]
\label{lem:cubic_threshold_classification}
Let $w\in C(J_{\max};L_x^2(\R^2))$ be a maximal-lifespan solution
of \eqref{eq:cubic_nls} with $M(w)=m_Q$. Suppose that, for some
$s_0\in J_{\max}$, its $L_{t,x}^4$ norm is infinite on at least
one of
\[
 (\inf J_{\max},s_0]\times\R^2,
 \qquad [s_0,\sup J_{\max})\times\R^2.
\]
Then, up to the symmetries of \eqref{eq:cubic_nls}, $w$ is the
soliton $e^{is}Q$ or a pseudoconformal transform of this soliton.
Here the symmetries include phase rotations, space and time
translations, Galilean transformations, mass-preserving dilations,
and complex conjugation together with time reversal.
\end{lemma}

\section{Proofs of the propagation estimates}
\label{app:propagation_proofs}

In this appendix we finish the proofs of
Lemmas~\ref{lem:dyadic_offdiagonal_S0},
\ref{lem:all_frequency_gap_crossing_S0}, and
\ref{lem:normal_crossing_full_frequency}.
We use the intervals, cutoffs, and multipliers in their statements
and the notation in \eqref{eq:directional_notation}.

\begin{lemma}[A uniform one-dimensional kernel estimate]
\label{lem:quadratic_phase_kernel}
Let $M>0$, $\rho,\tau\geq0$, and
$a\in C_c^\infty((1/2,2))$. For each integer $L\geq1$,
\begin{equation}\label{eq:quadratic_phase_kernel}
 \left|M\int_{\R}e^{\pm i(M\rho v+M^2\tau v^2)}a(v)\,dv\right|
 \leq C_LM(1+M\rho+M^2\tau)^{-L}
       \sum_{j=0}^L\|\partial_v^ja\|_{L_v^\infty(\R)}.
\end{equation}
The constant $C_L$ is independent of $M,\rho,\tau,a$, and the sign.
\end{lemma}
\begin{proof}
Set $\Lambda:=M\rho+M^2\tau$ and
$\Phi(v):=\pm(M\rho v+M^2\tau v^2)$.
For $\Lambda\leq1$, the estimate follows by direct integration.
For $\Lambda\geq1$ and $v\in(1/2,2)$,
\[
 |\Phi'(v)|\geq\Lambda,\qquad
 |\Phi''(v)|\leq2\Lambda,\qquad
 \left|\partial_v^k\frac1{\Phi'(v)}\right|
 \leq C_k\Lambda^{-1}\quad(k\geq0).
\]
On $(1/2,2)$, define $a_0:=a$ and
$a_{j+1}:=-\partial_v(a_j/(i\Phi'))$, extending each amplitude
by zero outside this interval.
Integration by parts preserves the support of the amplitudes and gives
\[
 \int_{\R}e^{i\Phi}a\,dv=\int_{\R}e^{i\Phi}a_L\,dv,
 \qquad
 \|a_L\|_{L_v^\infty(\R)}
 \leq C_L\Lambda^{-L}
        \sum_{j=0}^L\|\partial_v^ja\|_{L_v^\infty(\R)}.
\]
The support has length at most $3/2$, which proves
\eqref{eq:quadratic_phase_kernel}.
\end{proof}

\begin{proof}[Proof of Lemma~\ref{lem:dyadic_offdiagonal_S0}]
If $|j-k|\geq2$, then $|t-s|\sim\max(T_j,T_k)$ on
$I_j\times I_k$. The dispersive estimate and H\"older's inequality give
\[
 \begin{aligned}
 \left\|\int_{I_k}e^{i(t-s)\Delta}G_k(s)\,ds\right\|_{L^4_x(\R^2)}
 &\leq C\max(T_j,T_k)^{-1/2}
             \int_{I_k}\|G_k(s)\|_{L^{4/3}_x(\R^2)}\,ds\\
 &\leq C\max(T_j,T_k)^{-1/2}T_k^{1/4}
             \|G_k\|_{N^0(I_k)}.
 \end{aligned}
\]
Taking the $L^4_t(I_j)$ norm, we obtain
\begin{equation}\label{eq:dyadic_separated_appendix}
 \begin{aligned}
 \left\|\int_{I_k}e^{i(t-s)\Delta}G_k(s)\,ds
             \right\|_{L^4_{t,x}(I_j\times\R^2)}
 &\leq CT_j^{1/4}T_k^{1/4}\max(T_j,T_k)^{-1/2}
             \|G_k\|_{N^0(I_k)}\\
 &=C2^{-|j-k|/4}\|G_k\|_{N^0(I_k)}.
 \end{aligned}
\end{equation}

For $|j-k|\leq1$, choose $[t_-,t_+]\supset I_j\cup I_k$ and
extend $G_k$ by zero to $[t_-,t_+]$. For $t\in I_j$, we have
\[
 \begin{aligned}
 \int_{I_k\cap\{s<t\}}e^{i(t-s)\Delta}G_k(s)\,ds
 &=\int_{t_-}^t e^{i(t-s)\Delta}G_k(s)\,ds,\\
 \int_{I_k\cap\{s>t\}}e^{i(t-s)\Delta}G_k(s)\,ds
 &=-\int_{t_+}^t e^{i(t-s)\Delta}G_k(s)\,ds.
 \end{aligned}
\]
Applying the inhomogeneous Strichartz estimate
\eqref{eq:appendix_strichartz} with both pairs equal to $(4,4)$
and with $t_0=t_-$ and $t_0=t_+$, respectively, gives
\begin{equation}\label{eq:dyadic_adjacent_appendix}
 \begin{aligned}
 \left\|\int_{I_k\cap\{s<t\}}e^{i(t-s)\Delta}G_k(s)\,ds
             \right\|_{L^4_{t,x}(I_j\times\R^2)}
 &\leq C_1\|G_k\|_{N^0(I_k)}
 \leq C2^{-|j-k|/4}\|G_k\|_{N^0(I_k)},\\
 \left\|\int_{I_k\cap\{s>t\}}e^{i(t-s)\Delta}G_k(s)\,ds
             \right\|_{L^4_{t,x}(I_j\times\R^2)}
 &\leq C_1\|G_k\|_{N^0(I_k)}
 \leq C2^{-|j-k|/4}\|G_k\|_{N^0(I_k)}.
 \end{aligned}
\end{equation}
Here, we take
$C\geq2^{1/4}C_1$ since $2^{-|j-k|/4}\geq2^{-1/4}$.
This completes the proof of \eqref{eq:dyadic_offdiagonal_bound}.

We next prove \eqref{eq:dyadic_sequence_operator}.
Let $\widetilde b$ be the sequence on $\mathbb Z$ equal to $b_k$
for $k\geq0$ and 0 for $k<0$. Define the sequence
$w:=(w_m)_{m\in\mathbb Z}$ by $w_m:=2^{-|m|/4}$.
Then $\|w\|_{\ell^1(\mathbb Z)}=\sum_{m\in\mathbb Z}2^{-|m|/4}<\infty$ and
\[
 (w*\widetilde b)_j
 =\sum_{k\in\mathbb Z}w_{j-k}\widetilde b_k
 =\sum_{k\geq0}2^{-|j-k|/4}b_k
 =(\mathcal Tb)_j\qquad(j\geq0).
\]
Young's inequality gives
\[
 \begin{aligned}
 \|\mathcal Tb\|_{\ell^p}
 &\leq \|w*\widetilde b\|_{\ell^p(\mathbb Z)}
 \leq \|w\|_{\ell^1(\mathbb Z)}\|\widetilde b\|_{\ell^p(\mathbb Z)}\\
 &=\left(\sum_{m\in\mathbb Z}2^{-|m|/4}\right)\|b\|_{\ell^p}
 \qquad(1\leq p\leq\infty).
 \end{aligned}
\]
This proves \eqref{eq:dyadic_sequence_operator}.

To prove \eqref{eq:dyadic_terminal_integral}, fix $j$ with
$(\mathcal Tb)_j<\infty$. For finite $T\geq2T_j$, the integral
$\int_t^T e^{i(t-s)\Delta}G(s)\,ds$ belongs to $L^4_{t,x}(I_j\times\R^2)$ by
\eqref{eq:dyadic_separated_appendix} and \eqref{eq:dyadic_adjacent_appendix}.
For $K\geq j+2$ and
$T'>T\geq T_K$, apply \eqref{eq:dyadic_separated_appendix} to
$\mathbf1_{[T,T']}G_k$ to obtain
\[
 \left\|\int_T^{T'}e^{i(t-s)\Delta}G(s)\,ds
       \right\|_{L^4_{t,x}(I_j\times\R^2)}
 \leq C\sum_{k\geq K}2^{-|j-k|/4}b_k\longrightarrow0
 \qquad(K\to\infty).
\]
This proves \eqref{eq:dyadic_terminal_integral}.

To prove \eqref{eq:dyadic_duhamel_output}, fix $T\geq2T_j$ and apply
\eqref{eq:dyadic_separated_appendix} and
\eqref{eq:dyadic_adjacent_appendix} to
$\mathbf1_{[T_0,T]}G_k=\mathbf1_{I_k\cap[T_0,T]}G$. Since
$\|\mathbf1_{[T_0,T]}G_k\|_{N^0(I_k)}\leq b_k$, summing over $k$ gives
\begin{equation}\label{eq:dyadic_finite_terminal_output}
 \begin{aligned}
 &\left\|\int_{T_0}^t e^{i(t-s)\Delta}G(s)\,ds
                       \right\|_{L^4_{t,x}(I_j\times\R^2)}\\
 &\quad+\left\|\int_t^T e^{i(t-s)\Delta}G(s)\,ds
                       \right\|_{L^4_{t,x}(I_j\times\R^2)}
 \leq C\sum_{k\geq0}2^{-|j-k|/4}b_k
 =C(\mathcal Tb)_j.
 \end{aligned}
\end{equation}
Letting $T\to\infty$ in \eqref{eq:dyadic_finite_terminal_output}
and using the strong convergence in \eqref{eq:dyadic_terminal_integral}
gives \eqref{eq:dyadic_duhamel_output}.

\end{proof}

\begin{proof}[Proof of Lemma~\ref{lem:all_frequency_gap_crossing_S0}]
We first prove the $L^p$ bounds. For $g\in\mathcal S(\R)$, the
Hilbert transform is defined by
\[
 Hg(z):=\frac1\pi\lim_{\varepsilon\downarrow0}
             \int_{|y|>\varepsilon}\frac{g(z-y)}y\,dy.
\]
Its Fourier transform satisfies
$\widehat{Hg}(\xi)=-i\operatorname{sgn}(\xi)\widehat g(\xi)$.
By the theorem of M.~Riesz,
$H$ extends to a bounded operator on $L^p_z(\R)$ for $1<p<\infty$.
Define
\[
 \begin{aligned}
 H_e:\ L^p_x(\R^2)&\longrightarrow L^p_x(\R^2),\qquad f\longmapsto H_ef,\\
 (H_ef)(x_e,x_\perp)&:=H[f(\cdot,x_\perp)](x_e).
 \end{aligned}
\]
Fubini's theorem gives
\[
 \|H_ef\|_{L^p_x(\R^2)}^p
 =\int_{\R}\|H[f(\cdot,x_\perp)]\|_{L^p_{x_e}(\R)}^p\,dx_\perp
 \leq C_p^p\|f\|_{L^p_x(\R^2)}^p.
\]
Since
\[
 \mathbf1_{\{D_e>\alpha\}}f
 =e^{i\alpha x_e}\frac{I+iH_e}{2}(e^{-i\alpha x_e}f),
 \qquad
 \|\mathbf1_{\{D_e>\alpha\}}f\|_{L^p_x(\R^2)}\leq C_p\|f\|_{L^p_x(\R^2)},
\]
the bound for $\mathbf1_{\{D_e>\alpha\}}$ is uniform in $\alpha\in\R$ and $e\in\mathbb S^1$.
For the multipliers in \eqref{eq:gap_crossing_multipliers}, we have
\[
 \begin{aligned}
 \widehat{P^\pm f}(\xi)&=m_\pm(\xi_e)\widehat f(\xi),
 \qquad m_-(\xi_e):=1-m_+(\xi_e),\\
 m_+(\xi_e)&:=\eta\!\left(\frac{\xi_e-\mu}{\delta}\right)
 =\int_{\R}\eta'(v)\mathbf1_{\{\xi_e>\mu+\delta v\}}\,dv .
 \end{aligned}
\]
Thus, for $f\in\mathcal S(\R^2)$, Fourier inversion and Minkowski's
inequality yield
\begin{equation}\label{eq:gap_crossing_multiplier_Lp}
 \begin{aligned}
 P^+f&=\int_{\R}\eta'(v)\,
                    \mathbf1_{\{D_e>\mu+\delta v\}}f\,dv,\\
 \|P^+f\|_{L^p_x(\R^2)}
 &\leq\int_{\R}|\eta'(v)|\,
                \|\mathbf1_{\{D_e>\mu+\delta v\}}f\|_{L^p_x(\R^2)}\,dv
 \leq C_p\|\eta'\|_{L^1_v(\R)}\|f\|_{L^p_x(\R^2)}.
 \end{aligned}
\end{equation}
By density, \eqref{eq:gap_crossing_multiplier_Lp} holds for every
$f\in L^p_x(\R^2)$.
The bound for $P^-=I-P^+$ follows.

For $s,t\geq T_0$, define
\[
 T_{t,s}^\pm f:=\psi(t)e^{i(t-s)\Delta}P^\pm(\varphi(s)f),
 \qquad f\in L^2(\R^2).
\]
Thus
\begin{equation}\label{eq:gap_crossing_operators}
 \begin{aligned}
 \psi(t)\int_t^S e^{i(t-s)\Delta}P^+\varphi(s)\mathcal G(s)\,ds
   &=\int_t^S T_{t,s}^+\mathcal G(s)\,ds\qquad(S>t),\\
 C^-(t)&=\int_{T_0}^t T_{t,s}^-\mathcal G(s)\,ds.
 \end{aligned}
\end{equation}
We will prove that $\int_t^S T_{t,s}^+\mathcal G(s)\,ds$
converges strongly in $H^{3/4}$ as $S\to\infty$; its limit is $C^+(t)$.
We first estimate the $L^2$ operator norms of $T_{t,s}^+$ for
$s\geq t$ and $T_{t,s}^-$ for $s\leq t$.

To decompose the multipliers in the frequency $\xi_e$, set
\begin{equation}\label{eq:gap_crossing_frequency_changes}
 \gamma_+(\xi_e):=\xi_e-A/2,
 \qquad \gamma_-(\xi_e):=B/2-\xi_e.
\end{equation}
Define, for $\gamma\in\R$,
\begin{equation}\label{eq:gap_crossing_shifted_symbols}
 \begin{aligned}
 q_+(\gamma)&:=m_+(A/2+\gamma)
              =\eta\!\left(\frac\gamma\delta-\frac14\right),\\
 q_-(\gamma)&:=m_-(B/2-\gamma)
              =1-\eta\!\left(\frac14-\frac\gamma\delta\right).
 \end{aligned}
\end{equation}
Thus $m_\pm(\xi_e)=q_\pm(\gamma_\pm(\xi_e))$.
The functions $q_\pm$ vanish for $\gamma\leq\delta/8$ and equal
one for $\gamma\geq3\delta/8$.
Choose a nonnegative $\chi\in C_c^\infty((1/2,2))$ such that
$\sum_{\ell\in\mathbb Z}\chi(2^{-\ell}v)=1$ for $v>0$.
With $M_0:=\delta/16$, define
\begin{equation}\label{eq:gap_crossing_frequency_partition}
 \begin{gathered}
 \widehat{P_M^\pm f}(\xi):=q_\pm(\gamma_\pm(\xi_e))
       \chi\!\left(\frac{\gamma_\pm(\xi_e)}M\right)\widehat f(\xi),
 \qquad M\in\{2^\ell M_0:\ell=0,1,\ldots\},\\
 P_{\leq R}^\pm:=\sum_{M\leq R}P_M^\pm,
 \qquad R\in[M_0,\infty).
 \end{gathered}
\end{equation}
Since $q_\pm(\gamma)=0$ for $\gamma\leq2M_0$, the terms with
$\ell<0$ vanish. Hence
\[
 \begin{gathered}
 0\leq q_\pm(\gamma_\pm(\xi_e))
       \sum_{M\leq R}\chi\!\left(\frac{\gamma_\pm(\xi_e)}M\right)\leq1,\\
 q_\pm(\gamma_\pm(\xi_e))
       \sum_{M\leq R}\chi\!\left(\frac{\gamma_\pm(\xi_e)}M\right)
       \longrightarrow m_\pm(\xi_e)\qquad(R\to\infty).
 \end{gathered}
\]
Dominated convergence in Fourier space gives
$P_{\leq R}^\pm f\to P^\pm f$ in $H^\sigma$ for every
$\sigma\geq0$ and $f\in H^\sigma$. Therefore, we first prove
the $L^2$ bounds with $P^\pm$ replaced by $P_{\leq R}^\pm$,
uniformly in $R$.
Set
\[
 T_{t,s,R}^\pm f
 :=\psi(t)e^{i(t-s)\Delta}P_{\leq R}^\pm(\varphi(s)f),
 \qquad R\geq M_0.
\]

If $\delta^2\max(t,s)\leq1$, unitarity and
\eqref{eq:gap_crossing_cutoffs} give
\[
 \begin{aligned}
 \|T_{t,s}^\pm\|_{L^2_x(\R^2)\to L^2_x(\R^2)}&\leq C,\\
 \|T_{t,s,R}^\pm\|_{L^2_x(\R^2)\to L^2_x(\R^2)}
   &\leq C\qquad(R\geq M_0).
 \end{aligned}
\]
For every $J\geq1$, we have
$C\leq 2^JC(1+\delta^2\max(t,s))^{-J}$ in this case.
Thus, with $C_J\geq2^JC$,
\begin{equation}\label{eq:gap_crossing_small_parameter_L2}
 \max\bigl\{\|T_{t,s}^\pm\|_{L^2_x(\R^2)\to L^2_x(\R^2)},
             \|T_{t,s,R}^\pm\|_{L^2_x(\R^2)\to L^2_x(\R^2)}\bigr\}
 \leq C_J(1+\delta^2\max(t,s))^{-J}.
\end{equation}

Suppose now that $\delta^2\max(t,s)\geq1$.
We estimate the one-dimensional kernels of
$e^{i(t-s)\partial_{x_e}^2}P_M^\pm$ on $x_e\geq Bt$ and $y_e\leq As$,
as required by \eqref{eq:gap_crossing_cutoffs}. Set
\begin{equation}\label{eq:gap_crossing_spatial_coordinates}
 \tau:=|t-s|,\qquad r:=x_e-Bt\geq0,\qquad h:=As-y_e\geq0,
 \qquad \rho:=r+h+\delta\min(t,s).
\end{equation}
For fixed $s,t$, the one-dimensional kernel of
$e^{i(t-s)\partial_{x_e}^2}P_M^\pm$, evaluated at
$x_e=Bt+r$ and $y_e=As-h$, is
\begin{equation}\label{eq:gap_crossing_frequency_kernel}
 K_M(r,h):=\frac1{2\pi}\int_{\R}
 e^{i[(x_e-y_e)\xi_e-(t-s)\xi_e^2]}
 q_\pm(\gamma_\pm(\xi_e))
 \chi\!\left(\frac{\gamma_\pm(\xi_e)}M\right)\,d\xi_e.
\end{equation}
This integral is absolutely convergent because its integrand has
compact support in $\xi_e$.
Substituting $\xi_e=A/2+\gamma_+$ or $\xi_e=B/2-\gamma_-$
from \eqref{eq:gap_crossing_frequency_changes} gives
\begin{equation}\label{eq:gap_crossing_phase}
 (x_e-y_e)\xi_e-(t-s)\xi_e^2=
 \begin{cases}
 c_++\rho\gamma_++\tau\gamma_+^2,&s\geq t,\\
 c_- -\rho\gamma_- -\tau\gamma_-^2,&s\leq t,
 \end{cases}
\end{equation}
where
\[
 c_+:=\tfrac A2(x_e-y_e)+\tfrac{A^2}4\tau,\qquad
 c_-:=\tfrac B2(x_e-y_e)-\tfrac{B^2}4\tau.
\]
The numbers $c_\pm$ are real and independent of $\xi_e$, so
$|e^{ic_\pm}|=1$. In \eqref{eq:gap_crossing_frequency_kernel},
put $v:=\gamma_\pm(\xi_e)/M$ and define
\[
 a_M^\pm(v):=q_\pm(Mv)\chi(v),\qquad
 \Phi_M^\pm(v):=\pm(M\rho v+M^2\tau v^2).
\]
Then
\begin{equation}\label{eq:gap_crossing_rescaled_kernel}
 K_M(r,h)=\frac{M e^{ic_\pm}}{2\pi}
                 \int_{\R}e^{i\Phi_M^\pm(v)}a_M^\pm(v)\,dv.
\end{equation}
By \eqref{eq:gap_crossing_shifted_symbols}, for every integer $j\geq1$,
\[
 \begin{aligned}
 \partial_v^j[q_+(Mv)]
 &=\left(\frac M\delta\right)^j
      \eta^{(j)}\!\left(\frac{Mv}{\delta}-\frac14\right),\\
 \partial_v^j[q_-(Mv)]
 &=(-1)^{j+1}\left(\frac M\delta\right)^j
      \eta^{(j)}\!\left(\frac14-\frac{Mv}{\delta}\right).
 \end{aligned}
\]
These derivatives vanish unless $\delta/8\leq Mv\leq3\delta/8$.
On $\operatorname{supp}\chi\subset(1/2,2)$, this implies
$\delta/16\leq M\leq3\delta/4$, so $M/\delta$ is uniformly bounded.
Together with $|q_\pm|\leq1$, the product rule therefore gives
\begin{equation}\label{eq:gap_crossing_amplitude_bound}
 \|\partial_v^k a_M^\pm\|_{L^\infty_v(\R)}\leq C_k
 \qquad(k=0,1,\ldots).
\end{equation}

Applying Lemma~\ref{lem:quadratic_phase_kernel} to
\eqref{eq:gap_crossing_rescaled_kernel}, with the amplitude bounds
\eqref{eq:gap_crossing_amplitude_bound}, gives
\begin{equation}\label{eq:gap_crossing_dyadic_kernel}
 |K_M(r,h)|\leq C_LM(1+M\rho+M^2\tau)^{-L}.
\end{equation}

We sum \eqref{eq:gap_crossing_dyadic_kernel} to obtain the $L^2$
bound for $T_{t,s}^\pm$. Set
$Z:=\delta(r+h)+\delta^2\max(t,s)\geq1$.
By \eqref{eq:gap_crossing_spatial_coordinates} and $M\geq M_0=\delta/16$,
\begin{equation}\label{eq:gap_crossing_phase_lower_bound}
 \begin{aligned}
 M\rho+M^2\tau
 &\geq M(\rho+M_0\tau)\\
 &=\frac{M}{16\delta}Z+\frac{15M}{16}
                    \bigl(r+h+\delta\min(t,s)\bigr)
 \geq\frac{M}{16\delta}Z.
 \end{aligned}
\end{equation}
Using \eqref{eq:gap_crossing_phase_lower_bound} to sum
\eqref{eq:gap_crossing_dyadic_kernel} over
$M=2^\ell M_0$, $\ell\geq0$, we obtain, for $L>1$,
\begin{equation}\label{eq:gap_crossing_kernel_sum}
 \sum_{M\geq M_0}|K_M(r,h)|
 \leq C_L(1+Z)^{-L}\sum_{M\geq M_0}M(M/\delta)^{-L}
 \leq C_L\delta(1+Z)^{-L}.
\end{equation}
Take $L=J+1$, with $J\geq1$, and set
\begin{equation}\label{eq:gap_crossing_schur_kernel}
 k(r,h):=\delta(1+Z)^{-J-1}
 =\delta[1+\delta(r+h)+\delta^2\max(t,s)]^{-J-1}.
\end{equation}
Direct integration gives
\begin{equation}\label{eq:gap_crossing_schur_integrals}
 \max\left\{
 \sup_{r\geq0}\int_0^\infty k(r,h)\,dh,
 \sup_{h\geq0}\int_0^\infty k(r,h)\,dr\right\}
 \leq C_J(1+\delta^2\max(t,s))^{-J}.
\end{equation}
For $R\geq M_0$, define
\[
 K_R(r,h):=\sum_{M\leq R}K_M(r,h).
\]
By \eqref{eq:gap_crossing_kernel_sum} with $L=J+1$ and
\eqref{eq:gap_crossing_schur_kernel},
\begin{equation}\label{eq:gap_crossing_partial_sum_kernel}
 |K_R(r,h)|\leq\sum_{M\leq R}|K_M(r,h)|
 \leq C_J k(r,h).
\end{equation}
By \eqref{eq:gap_crossing_frequency_partition} and
\eqref{eq:gap_crossing_frequency_kernel}, for $f\in\mathcal S(\R^2)$
and $r\geq0$,
\[
 \begin{aligned}
 (T_{t,s,R}^\pm f)(Bt+r,\cdot)
 ={}&\psi(t,Bt+r,\cdot)\int_0^\infty K_R(r,h)\,
 e^{i(t-s)\partial_{x_\perp}^2}\\
 &\qquad\qquad
 [\varphi(s,As-h,\cdot)f(As-h,\cdot)]\,dh.
 \end{aligned}
\]
Since $e^{i(t-s)\partial_{x_\perp}^2}$ is unitary on
$L^2_{x_\perp}(\R)$, Minkowski's inequality,
\eqref{eq:gap_crossing_cutoffs}, and
\eqref{eq:gap_crossing_partial_sum_kernel} yield
\begin{equation}\label{eq:gap_crossing_partial_sum_bound}
 \|(T_{t,s,R}^\pm f)(Bt+r,\cdot)\|_{L^2_{x_\perp}(\R)}
 \leq C_J\int_0^\infty k(r,h)
                    \|f(As-h,\cdot)\|_{L^2_{x_\perp}(\R)}\,dh.
\end{equation}
Apply Schur's inequality to the integral operator with kernel $k(r,h)$.
By \eqref{eq:gap_crossing_schur_integrals}, its norm from
$L^2_h([0,\infty))$ to $L^2_r([0,\infty))$ is at most
$C_J(1+\delta^2\max(t,s))^{-J}$. Since $T_{t,s,R}^\pm f$ vanishes
for $x_e<Bt$ and
\[
 \int_0^\infty\|f(As-h,\cdot)\|_{L^2_{x_\perp}(\R)}^2\,dh
 =\|f\|_{L^2_x(\{x\in\R^2:x_e\leq As\})}^2
 \leq\|f\|_{L^2_x(\R^2)}^2,
\]
integrating \eqref{eq:gap_crossing_partial_sum_bound} in $L^2_r([0,\infty))$
and using density gives
\begin{equation}\label{eq:gap_crossing_large_parameter_L2}
 \|T_{t,s,R}^\pm\|_{L^2_x(\R^2)\to L^2_x(\R^2)}
 \leq C_J(1+\delta^2\max(t,s))^{-J}
 \qquad(R\geq M_0).
\end{equation}
Combining \eqref{eq:gap_crossing_small_parameter_L2} and
\eqref{eq:gap_crossing_large_parameter_L2}, then letting
$R\to\infty$ and using
$P_{\leq R}^\pm g\to P^\pm g$ in $L^2$ for each $g\in L^2$, we obtain
\begin{equation}\label{eq:all_frequency_gap_crossing_L2}
 \|T_{t,s}^\pm\|_{L^2_x(\R^2)\to L^2_x(\R^2)}
 \leq C_J(1+\delta^2\max(t,s))^{-J}.
\end{equation}

To obtain the bound on $H^{3/4}$, we also estimate one spatial derivative.
For $f\in H^1$,
\[
 \partial_jT_{t,s}^\pm f
 = (\partial_j\psi(t))e^{i(t-s)\Delta}P^\pm\varphi(s)f
   +T_{t,s}^\pm\partial_jf
   +\psi(t)e^{i(t-s)\Delta}P^\pm(\partial_j\varphi(s))f.
\]
For $\theta=\psi(t),\varphi(s)$ one has
$\operatorname{supp}\partial_j\theta\subset\operatorname{supp}\theta$.
The proof of \eqref{eq:all_frequency_gap_crossing_L2} uses
only the supports and $L^\infty$ bounds of the two cutoff factors.
Repeating that proof with $\psi(t)$ replaced by $\partial_j\psi(t)$,
or with $\varphi(s)$ replaced by $\partial_j\varphi(s)$, and using
\eqref{eq:gap_crossing_cutoffs}, we obtain
\begin{equation}\label{eq:gap_crossing_cutoff_derivative_L2}
 \begin{aligned}
 \|(\partial_j\psi(t))e^{i(t-s)\Delta}P^\pm(\varphi(s)f)\|_{L^2_x(\R^2)}
 &\leq\frac{C_J}{\delta t}
          (1+\delta^2\max(t,s))^{-J}\|f\|_{L^2_x(\R^2)},\\
 \|\psi(t)e^{i(t-s)\Delta}P^\pm((\partial_j\varphi(s))f)\|_{L^2_x(\R^2)}
 &\leq\frac{C_J}{\delta s}
          (1+\delta^2\max(t,s))^{-J}\|f\|_{L^2_x(\R^2)}.
 \end{aligned}
\end{equation}
Combining \eqref{eq:gap_crossing_cutoff_derivative_L2} with
\eqref{eq:all_frequency_gap_crossing_L2} for
$T_{t,s}^\pm f$ and $T_{t,s}^\pm\partial_jf$, and using $s,t\geq1$,
gives
\begin{equation}\label{eq:gap_crossing_H1}
 \|T_{t,s}^\pm\|_{H^1\to H^1}
 \leq C_J(1+\delta^{-1})(1+\delta^2\max(t,s))^{-J}.
\end{equation}
Interpolating the $L^2_x(\R^2)\to L^2_x(\R^2)$ bound
\eqref{eq:all_frequency_gap_crossing_L2} and the
$H^1\to H^1$ bound \eqref{eq:gap_crossing_H1} gives
\begin{equation}\label{eq:gap_crossing_H34}
 \|T_{t,s}^\pm f\|_{H^{3/4}}
 \leq C_J(1+\delta^{-1})^{3/4}
         (1+\delta^2\max(t,s))^{-J}\|f\|_{H^{3/4}}.
\end{equation}

It remains to construct $C^+(t)$ and estimate $C^\pm(t)$.
For fixed $t\geq T_0$ and finite $S>t$,
\eqref{eq:gap_crossing_H34} ensures that
$\int_t^S T_{t,s}^+\mathcal G(s)\,ds$ is well defined in $H^{3/4}$.
Taking $J=2$ in \eqref{eq:gap_crossing_H34} and integrating
over $[S,S']$ gives,
for $S'>S\geq t$,
\begin{equation}\label{eq:gap_crossing_terminal_tail}
 \left\|\int_S^{S'}T_{t,s}^+\mathcal G(s)\,ds\right\|_{H^{3/4}}
 \leq C\delta^{-2}(1+\delta^{-1})^{3/4}(1+\delta^2S)^{-1}
       \|\mathcal G\|_{L^\infty_t([T_0,\infty);H^{3/4}_x(\R^2))}.
\end{equation}
It follows that
\[
 C^+(t)=\lim_{S\to\infty}\int_t^S T_{t,s}^+\mathcal G(s)\,ds
 \quad\text{in }H^{3/4}.
\]
The convergence is uniform for $t$ in compact intervals by
\eqref{eq:gap_crossing_terminal_tail}.
The integral for $C^-(t)$ in \eqref{eq:gap_crossing_operators}
is over the finite interval $[T_0,t]$.
In the two time directions,
\[
 \begin{aligned}
 \int_t^\infty(1+\delta^2s)^{-2}\,ds
 &=\delta^{-2}(1+\delta^2t)^{-1},\\
 \int_{T_0}^t(1+\delta^2t)^{-2}\,ds
 &\leq t(1+\delta^2t)^{-2}\leq \delta^{-2}(1+\delta^2t)^{-1}.
 \end{aligned}
\]
Thus \eqref{eq:gap_crossing_H34} yields
\begin{equation}\label{eq:gap_crossing_time_bound}
 \|C^\pm(t)\|_{H^{3/4}}
 \leq C\delta^{-2}(1+\delta^{-1})^{3/4}(1+\delta^2t)^{-1}
       \|\mathcal G\|_{L^\infty_t([T_0,\infty);H^{3/4}_x(\R^2))}.
\end{equation}
Finally,
\[
 \|(1+\delta^2t)^{-1}\|_{L^4_t([T_0,\infty))}
 =3^{-1/4}\delta^{-1/2}(1+\delta^2T_0)^{-3/4}
 \leq C\delta^{-1/2}.
\]
Together with $H^{3/4}(\R^2)\hookrightarrow L^4(\R^2)$ and
\eqref{eq:gap_crossing_time_bound}, this gives
\[
 \begin{aligned}
 \|C^+\|_{L^4_{t,x}([T_0,\infty)\times\R^2)}
 +\|C^-\|_{L^4_{t,x}([T_0,\infty)\times\R^2)}
 &\leq C\delta^{-5/2}(1+\delta^{-1})^{3/4}
       \|\mathcal G\|_{L^\infty_t([T_0,\infty);H^{3/4}_x(\R^2))}\\
 &\leq C\delta^{-4}\|\mathcal G\|_{L^\infty_t([T_0,\infty);H^{3/4}_x(\R^2))},
 \end{aligned}
\]
since $0<\delta\leq1$. This proves
\eqref{eq:all_frequency_gap_crossing_S0} for every $T_0\geq1$.
\end{proof}

\begin{proof}[Proof of Lemma~\ref{lem:normal_crossing_full_frequency}]
Fix $m\in\{0,1\}$. Use the multipliers $P_M^\pm$ and
$P_{\leq R}^\pm$ in \eqref{eq:gap_crossing_frequency_partition},
with $M_0:=\delta/16$ and $R\geq M_0$. For $s,t\geq T_0$, define
\begin{equation}\label{eq:normal_crossing_truncated_operator}
 T_{m,R}^\pm(t,s)f
 :=\mathbf1_{\Omega_t}\partial_e^m e^{i(t-s)\Delta}
                  P_{\leq R}^\pm(\varphi(s)f),
 \qquad f\in L^2(\R^2).
\end{equation}
Since $|\xi_e|\leq(A+B)/2+2R$ on the frequency support of
$P_{\leq R}^\pm$, the operator $\partial_e^mP_{\leq R}^\pm$
is bounded on $L^2$. Thus, for $S>t$, the integrals
\begin{equation}\label{eq:normal_crossing_finite_integrals}
 \begin{aligned}
 W_{m,R,S}^+(t)&:=\int_t^S\partial_e^m e^{i(t-s)\Delta}
                         P_{\leq R}^+(\varphi(s)\mathcal G(s))\,ds,\\
 W_{m,R}^-(t)&:=\int_{T_0}^t\partial_e^m e^{i(t-s)\Delta}
                         P_{\leq R}^-(\varphi(s)\mathcal G(s))\,ds
 \end{aligned}
\end{equation}
are well defined in $L^2$, and
\begin{equation}\label{eq:normal_crossing_operator_integrals}
 \begin{aligned}
 \mathbf1_{\Omega_t}W_{m,R,S}^+(t)
   &=\int_t^S T_{m,R}^+(t,s)\mathcal G(s)\,ds,\\
 \mathbf1_{\Omega_t}W_{m,R}^-(t)
   &=\int_{T_0}^t T_{m,R}^-(t,s)\mathcal G(s)\,ds.
 \end{aligned}
\end{equation}
The proof follows the argument for
Lemma~\ref{lem:all_frequency_gap_crossing_S0}.
We first prove $L^2$ and $L^4$ bounds for $T_{m,R}^\pm(t,s)$
that are uniform in $R$, with $+$ for $s\geq t$ and $-$ for
$s\leq t$. Integrating these bounds will allow us to take
$R\to\infty$ on each finite time interval and then $S\to\infty$.

Use the coordinates in \eqref{eq:gap_crossing_spatial_coordinates}:
$\tau:=|t-s|$, $x_e=Bt+r$, $y_e=As-h$, and
$\rho:=r+h+\delta\min(t,s)$, where $r,h\geq0$.
The phase identities \eqref{eq:gap_crossing_phase} and the
frequency decomposition \eqref{eq:gap_crossing_frequency_partition}
hold for every $\delta>0$. The Fourier differentiation formula gives
$\widehat{\partial_e^m f}(\xi)=(i\xi_e)^m\widehat f(\xi)$.
Thus the kernel of
$\partial_{x_e}^m e^{i(t-s)\partial_{x_e}^2}P_M^\pm$
is obtained by multiplying the integrand in
\eqref{eq:gap_crossing_frequency_kernel} by $(i\xi_e)^m$.
With the same change of variables $v=\gamma_\pm(\xi_e)/M$
as in \eqref{eq:gap_crossing_rescaled_kernel}, we accordingly define
\begin{equation}\label{eq:normal_crossing_derivative_amplitudes}
 \begin{aligned}
 a_{M,m}^+(v)&:=(i(A/2+Mv))^m a_M^+(v),\\
 a_{M,m}^-(v)&:=(i(B/2-Mv))^m a_M^-(v).
 \end{aligned}
\end{equation}
By \eqref{eq:gap_crossing_rescaled_kernel}, the kernel of
$\partial_{x_e}^m e^{i(t-s)\partial_{x_e}^2}P_M^\pm$,
evaluated at $x_e=Bt+r$ and $y_e=As-h$, is
\begin{equation}\label{eq:normal_crossing_frequency_kernel}
 K_{M,m}(r,h):=\frac{Me^{ic_\pm}}{2\pi}
       \int_{\R}e^{i\Phi_M^\pm(v)}a_{M,m}^\pm(v)\,dv.
\end{equation}
The product rule and \eqref{eq:gap_crossing_amplitude_bound} give, for every integer $J\geq1$,
\begin{equation}\label{eq:normal_crossing_amplitude_bound}
 \max_{0\leq k\leq J}\|\partial_v^k a_{M,m}^\pm\|_{L^\infty_v(\R)}
 \leq C_{m,J}(1+A+B+M)^m.
\end{equation}
Lemma~\ref{lem:quadratic_phase_kernel}, applied with amplitude
$a_{M,m}^\pm$, and \eqref{eq:normal_crossing_amplitude_bound} give
\begin{equation}\label{eq:normal_crossing_dyadic_kernel}
 |K_{M,m}(r,h)|
 \leq C_{m,J}(1+A+B+M)^m M(1+M\rho+M^2\tau)^{-J}.
\end{equation}

For $J>m+1$, we have
\begin{equation}\label{eq:normal_crossing_geometric_sum}
 \begin{gathered}
 (1+A+B+M)^m\leq(1+A+B+M_0)^m(M/M_0)^m,\qquad
 1+M\rho+M^2\tau\geq M(\rho+M_0\tau),\\
 M_0^{-m}\sum_{\ell=0}^\infty(2^\ell M_0)^{m+1-J}
   =\frac{M_0^{1-J}}{1-2^{m+1-J}}.
 \end{gathered}
\end{equation}
Summing \eqref{eq:normal_crossing_dyadic_kernel} by
\eqref{eq:normal_crossing_geometric_sum}, we obtain
\begin{equation}\label{eq:normal_crossing_summed_kernel}
 \sum_{M\geq M_0}|K_{M,m}(r,h)|
 \leq C_{m,J}(1+A+B+M_0)^m M_0^{1-J}
              (\rho+M_0\tau)^{-J}.
\end{equation}
This proves absolute convergence for all $r,h\geq0$.
Set $L:=\delta\min(t,s)+M_0\tau$. Then
\begin{equation}\label{eq:normal_crossing_separation}
 \rho+M_0\tau=r+h+L,\qquad
 L\geq\frac{\delta}{16}\max(t,s)
       \geq\frac{\delta}{32}(t+\tau).
\end{equation}
The last inequality uses $s=t+\tau$ when $s\geq t$ and
$0\leq\tau\leq t$ when $s\leq t$.

We now use \eqref{eq:normal_crossing_summed_kernel} to
estimate the corresponding integral operators. For this purpose,
let $K:[0,\infty)^2\to\mathbb C$ be any measurable kernel satisfying
$|K(r,h)|\leq C_*(r+h+L)^{-J}$, where $L>0$ and $J>1$, and define
\[
 (\mathcal Af)(r):=\int_0^\infty K(r,h)f(h)\,dh,
 \qquad f\in L^2_h([0,\infty)),\quad r\geq0.
\]
Direct integration gives
\begin{equation}\label{eq:normal_crossing_schur_integrals}
 \max\left\{
 \sup_{r\geq0}\int_0^\infty|K(r,h)|\,dh,\,
 \sup_{h\geq0}\int_0^\infty|K(r,h)|\,dr
 \right\}\leq\frac{C_*}{J-1}L^{1-J}.
\end{equation}
Schur's inequality and \eqref{eq:normal_crossing_schur_integrals} give
\begin{equation}\label{eq:normal_crossing_spatial_L2}
 \|\mathcal Af\|_{L^2_r([0,\infty))}
 \leq\frac{C_*}{J-1}L^{1-J}\|f\|_{L^2_h([0,\infty))}.
\end{equation}
By Cauchy--Schwarz,
\begin{equation}\label{eq:normal_crossing_spatial_Linfty}
 \begin{aligned}
 \|\mathcal Af\|_{L^\infty_r([0,\infty))}
 &\leq C_*\sup_{r\geq0}
       \left(\int_0^\infty(r+h+L)^{-2J}\,dh\right)^{1/2}
              \|f\|_{L^2_h([0,\infty))}\\
 &=\frac{C_*}{\sqrt{2J-1}}L^{1/2-J}\|f\|_{L^2_h([0,\infty))}.
 \end{aligned}
\end{equation}
The inequality
\[
 \|\mathcal Af\|_{L^4_r([0,\infty))}^4
 \leq\|\mathcal Af\|_{L^\infty_r([0,\infty))}^2
       \|\mathcal Af\|_{L^2_r([0,\infty))}^2,
\]
together with \eqref{eq:normal_crossing_spatial_L2} and
\eqref{eq:normal_crossing_spatial_Linfty}, yields
\begin{equation}\label{eq:normal_crossing_spatial_bounds}
 \|\mathcal Af\|_{L^4_r([0,\infty))}
 \leq\frac{C_*L^{3/4-J}}{(J-1)^{1/2}(2J-1)^{1/4}}
           \|f\|_{L^2_h([0,\infty))}.
\end{equation}

Define
\[
 K_{m,R}(r,h):=\sum_{M\leq R}K_{M,m}(r,h),\qquad
 (\mathcal A_{m,R}f)(r):=
   \int_0^\infty\sum_{M\leq R}|K_{M,m}(r,h)|f(h)\,dh.
\]
Take $J=4$ in \eqref{eq:normal_crossing_summed_kernel}.
Since $M_0=\delta/16\leq(A+B)/16$, the kernel of
$\mathcal A_{m,R}$ is bounded by
$C(1+A+B)^m\delta^{-3}(r+h+L)^{-4}$.
Applying \eqref{eq:normal_crossing_spatial_L2} and
\eqref{eq:normal_crossing_spatial_bounds}, followed by
\eqref{eq:normal_crossing_separation}, we obtain, uniformly in $R$,
\begin{equation}\label{eq:normal_crossing_scalar_operator_bounds}
 \begin{aligned}
 \|\mathcal A_{m,R}\|_{L^2_h([0,\infty))\to L^2_r([0,\infty))}
 &\leq C(1+A+B)^m\delta^{-6}(t+\tau)^{-3},\\
 \|\mathcal A_{m,R}\|_{L^2_h([0,\infty))\to L^4_r([0,\infty))}
 &\leq C(1+A+B)^m\delta^{-25/4}(t+\tau)^{-13/4}.
 \end{aligned}
\end{equation}

The kernel representation \eqref{eq:normal_crossing_frequency_kernel}
gives, for $r>0$,
\begin{equation}\label{eq:normal_crossing_operator_kernel}
 \begin{aligned}
 [T_{m,R}^\pm(t,s)f](Bt+r,\cdot)
 =\int_0^\infty K_{m,R}(r,h)e^{i(t-s)\partial_{x_\perp}^2}
       [\varphi(s,As-h)f(As-h,\cdot)]\,dh.
 \end{aligned}
\end{equation}
Let $\mathcal H:=L^2_{x_\perp}(\R)$ or
$\mathcal H:=H^{1/4}_{x_\perp}(\R)$.
The group $e^{i(t-s)\partial_{x_\perp}^2}$ is unitary on
$\mathcal H$, and $\varphi(s,As-h)$ is a scalar of modulus
at most one. Taking the $\mathcal H$ norm in
\eqref{eq:normal_crossing_operator_kernel}, we find
\begin{equation}\label{eq:normal_crossing_hilbert_bound}
 \begin{gathered}
 \|[T_{m,R}^\pm(t,s)f](Bt+r,\cdot)\|_{\mathcal H}
 \leq\int_0^\infty\sum_{M\leq R}|K_{M,m}(r,h)|
                      \|f(As-h,\cdot)\|_{\mathcal H}\,dh,\\
 \int_0^\infty\|f(As-h,\cdot)\|_{\mathcal H}^2\,dh
 \leq\|f\|_{L^2_{x_e}(\R;\mathcal H)}^2.
 \end{gathered}
\end{equation}
Apply the $L^2_h([0,\infty))\to L^2_r([0,\infty))$ estimate in
\eqref{eq:normal_crossing_scalar_operator_bounds} with
$\mathcal H=L^2$, and its $L^2_h([0,\infty))\to L^4_r([0,\infty))$ estimate with
$\mathcal H=H^{1/4}$. In the latter case also use
$H^{1/4}_{x_\perp}(\R)\hookrightarrow L^4_{x_\perp}(\R)$.
Together with \eqref{eq:normal_crossing_hilbert_bound}, this gives
\begin{equation}\label{eq:normal_crossing_operator_bounds}
 \begin{aligned}
 \|T_{m,R}^\pm(t,s)f\|_{L^2_x(\R^2)}
 &\leq C(1+A+B)^m\delta^{-6}(t+\tau)^{-3}\|f\|_{L^2_x(\R^2)},\\
 \|T_{m,R}^\pm(t,s)f\|_{L^4_x(\R^2)}
 &\leq C\|[T_{m,R}^\pm(t,s)f](Bt+r,\cdot)
                    \|_{L^4_r([0,\infty);H^{1/4}_{x_\perp}(\R))}\\
 &\leq C(1+A+B)^m\delta^{-25/4}(t+\tau)^{-13/4}
                    \|f\|_{\mathcal H_e^{1/4}}.
 \end{aligned}
\end{equation}

We now integrate \eqref{eq:normal_crossing_operator_bounds} in
\eqref{eq:normal_crossing_operator_integrals}. The ranges for
$\tau$ are $[0,S-t]$ and $[0,t-T_0]$, respectively, and
\begin{equation}\label{eq:normal_crossing_time_integrals}
 \int_0^\infty(t+\tau)^{-3}\,d\tau=\tfrac12t^{-2},
 \qquad
 \int_0^\infty(t+\tau)^{-13/4}\,d\tau=\tfrac49t^{-9/4}.
\end{equation}
Minkowski's inequality and \eqref{eq:normal_crossing_time_integrals}
therefore give, for $W=W_{m,R,S}^+(t)$ or $W=W_{m,R}^-(t)$,
\begin{equation}\label{eq:normal_crossing_finite_integral_bounds}
 \begin{aligned}
 \|\mathbf1_{\Omega_t}W\|_{L^2_x(\R^2)}
 &\leq C(1+A+B)^m\delta^{-6}t^{-2}
       \|\mathcal G\|_{L^\infty_t([T_0,\infty);L^2_x(\R^2))},\\
 \|\mathbf1_{\Omega_t}W\|_{L^4_x(\R^2)}
 &\leq C(1+A+B)^m\delta^{-25/4}t^{-9/4}
       \|\mathcal G\|_{L^\infty_t([T_0,\infty);\mathcal H_e^{1/4})}.
 \end{aligned}
\end{equation}
These estimates are uniform in $R,S$.
For later time integration, note that
$s\mapsto\varphi(s)\mathcal G(s)$ is strongly measurable in
$L^2$, and the integrals in \eqref{eq:normal_crossing_finite_integrals}
are continuous in $t$ with values in $L^2$.
Multiplication by $\mathbf1_{\Omega_t}$ preserves this continuity:
for each $f\in L^2$, dominated convergence gives
$\|(\mathbf1_{\Omega_t}-\mathbf1_{\Omega_{t'}})f\|_{L^2_x(\R^2)}\to0$
as $t'\to t$.
Thus the finite integrals, and their strong $L^2$ limits below,
have measurable spacetime representatives.

We next take $R\to\infty$ for fixed $t\geq T_0$ and $S>t$.
By \eqref{eq:gap_crossing_frequency_partition},
$P_{\leq R}^+\to P^+$ strongly on $L^2$, with operator norms
at most one. Since
$\|\partial_e^m e^{i(t-s)\Delta}\|_{L^2_x(\R^2)\to H^{-m}_x(\R^2)}\leq1$,
dominated convergence in \eqref{eq:normal_crossing_finite_integrals}
gives
\begin{equation}\label{eq:crossing_Hminus_limit}
 W_{m,R,S}^+(t)\longrightarrow
 W_{m,S}^+(t):=\int_t^S\partial_e^m e^{i(t-s)\Delta}
                    P^+(\varphi(s)\mathcal G(s))\,ds
 \quad\hbox{in }H^{-m}(\R^2).
\end{equation}
To obtain convergence on $\Omega_t$, take $R\geq2(1+A+B+M_0)$. For $M>R$,
$(1+A+B+M)^m\leq C M^m$. Summing
\eqref{eq:normal_crossing_dyadic_kernel} with $J=4$ gives
\begin{equation}\label{eq:normal_crossing_frequency_kernel_tail}
 \sum_{M>R}|K_{M,m}(r,h)|
 \leq C(\rho+M_0\tau)^{-4}\sum_{M>R}M^{m-3}
 \leq CR^{m-3}(\rho+M_0\tau)^{-4}.
\end{equation}
Apply \eqref{eq:normal_crossing_spatial_L2} and
\eqref{eq:normal_crossing_spatial_bounds} with
$C_*=CR^{m-3}$ and $J=4$, and repeat
\eqref{eq:normal_crossing_hilbert_bound} for the sum over
$R<M\leq R'$. By \eqref{eq:normal_crossing_separation},
followed by time integration using
\eqref{eq:normal_crossing_time_integrals}, we obtain, uniformly
for $R'>R$ and $S>t$,
\begin{equation}\label{eq:crossing_frequency_tails}
 \begin{aligned}
 \|\mathbf1_{\Omega_t}(W_{m,R',S}^+-W_{m,R,S}^+)(t)\|_{L^2_x(\R^2)}
 &\leq CR^{m-3}\delta^{-3}t^{-2}
             \|\mathcal G\|_{L^\infty_t([T_0,\infty);L^2_x(\R^2))},\\
 \|\mathbf1_{\Omega_t}(W_{m,R',S}^+-W_{m,R,S}^+)(t)\|_{L^4_x(\R^2)}
 &\leq CR^{m-3}\delta^{-13/4}t^{-9/4}
             \|\mathcal G\|_{L^\infty_t([T_0,\infty);\mathcal H_e^{1/4})}.
 \end{aligned}
\end{equation}
Since $m=0,1$, the first estimate in
\eqref{eq:crossing_frequency_tails} proves that
$\mathbf1_{\Omega_t}W_{m,R,S}^+(t)$ converges strongly in $L^2$.
Denote its limit by $C_{m,S}^+(t)$.
For $\phi\in C_c^\infty(\Omega_t)$, its zero extension
$\widetilde\phi$ belongs to $C_c^\infty(\R^2)$.
The $H^{-m}$ convergence in \eqref{eq:crossing_Hminus_limit}
and the strong $L^2$ convergence give
\begin{equation}\label{eq:normal_crossing_restriction_pairing}
 \begin{aligned}
 \langle W_{m,S}^+(t)|_{\Omega_t},\phi\rangle
 &:=\langle W_{m,S}^+(t),\widetilde\phi\rangle\\
 &=\lim_{R\to\infty}\int_{\Omega_t}
             W_{m,R,S}^+(t,x)\overline{\phi(x)}\,dx\\
 &=\int_{\Omega_t}C_{m,S}^+(t,x)\overline{\phi(x)}\,dx.
 \end{aligned}
\end{equation}
This identifies $C_{m,S}^+(t)$ with the function specified in
the lemma. Under the additional assumption on $\mathcal G$,
the second estimate in \eqref{eq:crossing_frequency_tails}
gives a strong $L^4$ limit. Testing this limit against
$\phi\in C_c^\infty(\Omega_t)$ in
\eqref{eq:normal_crossing_restriction_pairing} shows that it
equals $C_{m,S}^+(t)$ almost everywhere on $\Omega_t$.
Both limits vanish outside $\Omega_t$.
Letting $R'\to\infty$ in \eqref{eq:crossing_frequency_tails}, we obtain
\begin{equation}\label{eq:normal_crossing_uniform_frequency_limit}
 \begin{aligned}
 \sup_{S>t}\|\mathbf1_{\Omega_t}W_{m,R,S}^+(t)-C_{m,S}^+(t)\|_{L^2_x(\R^2)}
 &\leq CR^{m-3}\delta^{-3}t^{-2}
       \|\mathcal G\|_{L^\infty_t([T_0,\infty);L^2_x(\R^2))}
       \longrightarrow0,\\
 \sup_{S>t}\|\mathbf1_{\Omega_t}W_{m,R,S}^+(t)-C_{m,S}^+(t)\|_{L^4_x(\R^2)}
 &\leq CR^{m-3}\delta^{-13/4}t^{-9/4}
       \|\mathcal G\|_{L^\infty_t([T_0,\infty);\mathcal H_e^{1/4})}
       \longrightarrow0
 \end{aligned}
\end{equation}
as $R\to\infty$.

For $W_{m,R}^-(t)$ defined in
\eqref{eq:normal_crossing_finite_integrals}, dominated convergence
on $[T_0,t]$ gives
\begin{equation}\label{eq:normal_crossing_initial_Hminus_limit}
 W_{m,R}^-(t)\longrightarrow
 W_m^-(t):=\int_{T_0}^t\partial_e^m e^{i(t-s)\Delta}
                       P^-(\varphi(s)\mathcal G(s))\,ds
 \quad\hbox{in }H^{-m}(\R^2).
\end{equation}
Using \eqref{eq:normal_crossing_frequency_kernel_tail} in
\eqref{eq:normal_crossing_spatial_L2},
\eqref{eq:normal_crossing_spatial_bounds}, and
\eqref{eq:normal_crossing_hilbert_bound}, and integrating over
$0\leq\tau\leq t-T_0$ by \eqref{eq:normal_crossing_time_integrals},
gives the two bounds in \eqref{eq:crossing_frequency_tails}
with $W_{m,R,S}^+$ replaced by $W_{m,R}^-$.
The test-function argument in
\eqref{eq:normal_crossing_restriction_pairing}, now using
\eqref{eq:normal_crossing_initial_Hminus_limit}, identifies
the strong limits as $C_m^-(t)$. Thus
\[
 \lim_{R\to\infty}
 \|\mathbf1_{\Omega_t}W_{m,R}^-(t)-C_m^-(t)\|_{L^2_x(\R^2)}=0,
 \qquad
 \lim_{R\to\infty}
 \|\mathbf1_{\Omega_t}W_{m,R}^-(t)-C_m^-(t)\|_{L^4_x(\R^2)}=0.
\]

It remains to take $S\to\infty$ in $C_{m,S}^+(t)$.
For $S'>S>t$, integrate \eqref{eq:normal_crossing_operator_bounds}
on $[S,S']$, where $t+\tau=s$, and then take $R\to\infty$
as justified by \eqref{eq:crossing_frequency_tails}.
Since $\int_S^\infty s^{-3}\,ds=\tfrac12S^{-2}$ and
$\int_S^\infty s^{-13/4}\,ds=\tfrac49S^{-9/4}$, we obtain
\begin{equation}\label{eq:crossing_source_tails}
 \begin{aligned}
 \|C_{m,S'}^+(t)-C_{m,S}^+(t)\|_{L^2_x(\R^2)}
 &\leq C(1+A+B)^m\delta^{-6}S^{-2}
             \|\mathcal G\|_{L^\infty_t([T_0,\infty);L^2_x(\R^2))},\\
 \|C_{m,S'}^+(t)-C_{m,S}^+(t)\|_{L^4_x(\R^2)}
 &\leq C(1+A+B)^m\delta^{-25/4}S^{-9/4}
             \|\mathcal G\|_{L^\infty_t([T_0,\infty);\mathcal H_e^{1/4})}.
 \end{aligned}
\end{equation}
Hence $C_{m,S}^+(t)$ converges strongly in $L^2$ as
$S\to\infty$ to a limit denoted by $C_m^+(t)$.
Under the additional assumption on $\mathcal G$, it also
converges strongly in $L^4$ to $C_m^+(t)$, by the
test-function argument in \eqref{eq:normal_crossing_restriction_pairing}.
For a compact interval $I\subset[T_0,\infty)$ and $S>\sup I$,
let $S'\to\infty$ in the second estimate of
\eqref{eq:crossing_source_tails} and take the $L^4_t(I)$ norm.
This gives
\begin{equation}\label{eq:normal_crossing_local_time_limit}
 \|C_{m,S}^+-C_m^+\|_{L^4_{t,x}(I\times\R^2)}
 \leq C|I|^{1/4}(1+A+B)^m\delta^{-25/4}S^{-9/4}
       \|\mathcal G\|_{L^\infty_t([T_0,\infty);\mathcal H_e^{1/4})}.
\end{equation}
Combining \eqref{eq:normal_crossing_uniform_frequency_limit}
with \eqref{eq:crossing_source_tails} and the triangle inequality
yields, for each fixed $t$,
\[
 \lim_{R,S\to\infty}
 \|\mathbf1_{\Omega_t}W_{m,R,S}^+(t)-C_m^+(t)\|_{L^p_x(\R^2)}=0,
 \qquad p=2,4.
\]

Finally, take $R\to\infty$ in
\eqref{eq:normal_crossing_finite_integral_bounds}, and then
$S\to\infty$ for $C_m^+$. Since
$\delta^{-6}\leq(1+\delta^{-1})^6$ and
$\delta^{-25/4}\leq(1+\delta^{-1})^7$, we obtain
\eqref{eq:normal_crossing_pointwise_L2} and
\eqref{eq:normal_crossing_pointwise_L4}.
Keeping $S=T$ finite gives the same pointwise estimates for
$C_{m,T}^+$ when $T_0\leq t<T$. Also note that
\begin{equation}\label{eq:normal_crossing_time_weight}
 \begin{aligned}
 \|t^{-9/4}\|_{L^4_t([T_0,T))}
 &\leq\|t^{-9/4}\|_{L^4_t([T_0,\infty))}
 =8^{-1/4}T_0^{-2}.
 \end{aligned}
\end{equation}
Taking the $L^4_t([T_0,\infty))$ norm in
\eqref{eq:normal_crossing_pointwise_L4} and using
\eqref{eq:normal_crossing_time_weight} gives
\eqref{eq:normal_crossing_global_L4}.
Taking the $L^4_t([T_0,T))$ norm of the corresponding estimate
for $C_{m,T}^+$ gives the same spacetime bound, uniformly in $T$.
\end{proof}

\section*{Acknowledgements}
The authors used AI tools for mathematical discussions and language
editing during the preparation of this manuscript. ChatGPT provided
key assistance with the proof of
Lemma~\ref{lem:all_frequency_gap_crossing_S0}, particularly in
formulating the uniform one-dimensional nonstationary-phase
kernel estimate used there. The authors independently checked and
revised all mathematical statements, proofs, and results, verified
the final arguments, and take full responsibility for the content
of this article.

\end{document}